\documentclass{article}
\usepackage{graphicx}
\usepackage{amsmath}
\usepackage{amsfonts}
\usepackage{amssymb}
\usepackage{ifthen}
\usepackage{bm,bbm}

\newtheorem{theorem}{Theorem}[section]
\newtheorem{lemma}[theorem]{Lemma}

\newtheorem{remark}[theorem]{Remark}
\newtheorem{definition}[theorem]{Definition}

\newcommand{\DEFINED}[1]{{\bf #1}}

\newcommand{\SET}[1]{\mathcal{#1}}

\newcommand{\BISTOCHASTIC}{bistochastic}

\newcommand{\Ii}{{\bf i}}

\newcommand{\bO}{\bullet}

\newcommand{\ZEROvect}{\mathbf{0}}

\newcommand{\ONESvect}{\mathbf{e}}

\newcommand{\PHASE}[1]{e^{\Ii #1}}

\newcommand{\STbasis}[1]{e_{#1}}

\newcommand{\REALS}{\mathbbm{R}}

\newcommand{\COMPLEX}{\mathbbm{C}}

\newcommand{\NATURAL}{\mathbbm{N}}

\newcommand{\BISTOCHSPACE}{\mathcal{B}}
\newcommand{\UNITARY}{\mathcal{U}}
\newcommand{\DIAGONAL}{\mathcal{D}}

\newcommand{\IDENTITY}{\mathbf{id}}

\newcommand{\CONJ}[1]{\overline{#1}}

\newcommand{\diag}{\mathbf{diag}}

\newcommand{\RELofEQUI}{\simeq}

\renewcommand{\Re}{\mathbf{Re}}
\renewcommand{\Im}{\mathbf{Im}}

\newcommand{\DIFFERENTIAL}[2]{{D \! #1}_{#2}}

\newcommand{\TANGENTspace}[2]{\mathbf{T}_{#2} #1}

\newcommand{\VECTORform}[1]{{\mathbf{vec}}_{#1}}
\newcommand{\VECR}{\VECTORform{\mathbbm{R}}}
\newcommand{\VECC}{\VECTORform{\mathbbm{C}}}
\newcommand{\VEC}{\VECTORform{}}

\newcommand{\SPANNEDspace}[1]{{\mathbf{span}}_{#1}}
\newcommand{\SPANR}{\SPANNEDspace{\mathbbm{R}}}
\newcommand{\SPANC}{\SPANNEDspace{\mathbbm{C}}}

\newcommand{\NULLspace}[1]{{\mathbf{null}}_{#1}}
\newcommand{\NULLR}{\NULLspace{\mathbbm{R}}}
\newcommand{\NULLC}{\NULLspace{\mathbbm{C}}}

\newcommand{\DEFECT}{\mathbf{d}}

\newcommand{\RANK}{\mathbf{rank}}

\newcommand{\lcm}{\mathbf{lcm}}

\newcommand{\NUMBERofIN}[2]{\sharp_{#2} #1}

\newcommand{\DMPbound}{\mathbf{b}}

\newcommand{\HADprod}{\circ}

\newcommand{\EXPentrywise}{\mathbf{EXP}}

\newcommand{\ELEMENTof}[3]{{\left[ #1 \right]}_{#2,#3}}

\newcommand{\DELTA}[2]
{
  \ifthenelse{\equal{#1}{} \or \equal{#2}{}}
    {\Delta^{#1#2}}
    {\Delta^{#1,#2}}
}

\newcommand{\PCM}[3]{#1^{#2,#3}}
\newcommand{\pcmELEMENTof}[3]{\left[ #1 \right]^{#2,#3}}

\newcommand{\PROOFstart}{{\bf Proof}\\}
\newcommand{\PROOFend}{$\blacksquare$}

\begin{document}

\title{Defect of a unitary matrix}

\author{
  Wojciech Tadej$^{1}$  and  Karol \.Zyczkowski$^{2,3}$\\
  \smallskip
 $^1${\small Faculty of Mathematics and Natural Sciences, College of Sciences,}          \\
     {\small Cardinal Stefan Wyszy{\'n}ski University, Warsaw, Poland}                   \\
                                                                                            \\
  $^2${\small Institute of Physics, Jagiellonian University, Krak{\'o}w, Poland}         \\
  $^3${\small Center for Theoretical Physics, Polish Academy of Sciences, Warsaw, Poland}\\
  \smallskip
  {\small e-mail: wtadej@wp.pl               \quad \quad
                  karol@tatry.if.uj.edu.pl   }\\
}
\date{December 17, 2007}

\maketitle

\begin{abstract}
 We analyze properties of a map $f$
 sending a unitary matrix $U$ of size $N$
 into a doubly stochastic matrix $B=f(U)$ defined by $B_{i,j}=|U_{i,j}|^2$. For any $U$
 we define its {\sl defect}, determined by the dimension of the image
    $\DIFFERENTIAL{f}{}( \TANGENTspace{\UNITARY}{U} )$
 of the space
    $\TANGENTspace{\UNITARY}{U}$
 tangent to the manifold of unitary matrices $\UNITARY$ at $U$ under the tangent map
 $\DIFFERENTIAL{f}{}$ corresponding to $f$.
 The defect, equal to zero for
 a generic unitary matrix, gives an upper bound for the dimension
 of a smooth orbit (a manifold) of inequivalent unitary matrices
 mapped into the same doubly stochastic matrix. We demonstrate several properties of the defect and
 prove an explicit formula for the defect of the Fourier matrix $F_N$ of size $N$.
 In this way we obtain an upper bound for the dimension
 of a smooth orbit of inequivalent unitary complex Hadamard matrices stemming from $F_N$.
 It is equal to zero iff $N$ is prime and coincides
 with the dimension of the known orbits if $N$ is a power of a
 prime.
 Two constructions of these orbits are presented at the end of this work.
\end{abstract}


{\bf MSC-class:}\ \ 58K05, 51F25, 15A51, 65T50, 05B20, 15A90

%
%

\section{Introduction}
\label{sec_introduction}

Consider the set $\UNITARY$ of unitary matrices of finite size $N$.
For any unitary $U$ we define a matrix $B= f(U)$ with
non--negative entries
\begin{equation}
  \label{eq_f0}
  B_{i,j} = |U_{i,j}|^2 \ .
\end{equation}
From the unitarity condition, $UU^*={\mathbbm 1}$, it follows that
the resulting matrix $B$ is \BISTOCHASTIC\ (also called {sl doubly stochastic}), 
since the sum of elements in each of its columns or rows is equal to
unity. A \BISTOCHASTIC\ matrix $B$ for which there exists a
unitary (an orthogonal)
 $U$ satisfying (\ref{eq_f0}) is called {\sl unistochastic}
({\sl orthostochastic}). For $N=2$ all \BISTOCHASTIC\
matrices are unistochastic, even orthostochastic, but for $N\ge 3$
it is no longer the case \cite{MO79,ZKSS03}.

Our work is motivated by the following problem \cite{AMM91,BEKTZ05}.
\smallskip

\noindent {\bf (*)} {\sl For a given unitary $U\in {\cal U}_N$
find all other unitary matrices $V\in {\cal U}_N$ such that
$f(V)=B=f(U)$.}
\smallskip

This rather general question is closely related to several
problems in various branches of mathematics and theoretical
physics. For instance, taking the Fourier matrix $F_N$ as the
unitary $U$ in question we get the flat \BISTOCHASTIC\ matrix,
$B=J_N$ with $\ELEMENTof{J_N}{i}{j}=1/N$, so the above question
reduces to the problem of finding all unitary\footnote {We reserve
the term 'complex Hadamard matrix' for an $N \times N$ complex
matrix $H$ satisfying $H^*H = HH^* = N \cdot I_N$ and
  $\forall i,j \ \ |H_{i,j}|=1$, while $1/\sqrt{N} \cdot H$ is called
  a 'unitary complex Hadamard matrix'. Such matrices were called
   by Craigen 'unit Hadamard matrices' \cite{Cr91}.}
complex Hadamard matrices of size $N$.
This issue is related to construction of some  $*$-subalgebras in
finite von Neumann algebras \cite{Po83,MW92,Ha96}, analyzing
bi-unimodular sequences or finding cyclic $n$--roots
\cite{BF91,BS95} and equiangular lines \cite{GR05}. The search for
complex Hadamard matrices \cite{Ha96,Di04} is also motivated by
the theory of quantum information processing
\cite{We01,WGC03,TZ06}.

Furthermore, the general issue
of specifying all unitary matrices
such that their squared moduli give a fixed
\BISTOCHASTIC\ matrix was intensively studied
by high energy physicists investigating the parity violation
and analyzing the Cabibbo--Kobayashi--Maskawa matrices
\cite{Ja85,BD87,Di05}.
On the other hand, relation (\ref{eq_f0})
is relevant to investigation
of the semiclassical limit of quantum mechanics:
for a given \BISTOCHASTIC\ $B$,
representing the transition matrix of a Markov chain,
one looks for the set of unitary matrices $U$
which lead to the corresponding
quantum dynamics  \cite{KS97,Ta01,PZK01}.

To investigate problem (*) one uses the notion of {\sl equivalent}
unitary matrices  \cite{Ha96}, which differ by left and right
multiplication by diagonal unitary matrices followed by arbitrary
permutations of rows and columns. We suppose that for a generic
unitary $U$ all solutions of the problem (*) in a neighbourhood of
$U$ are equivalent, and we call such $U$ {\sl isolated}. However, for
some non--typical unitaries it is not the case. It is therefore
natural to ask for the dimension of a smooth orbit, i.e. a
manifold, stemming from $U$, if one exists, of non--equivalent
solutions $V$ of problem (*) posed for a given unitary $U$. The
upper bound for this dimension is obtained in this paper by
computing the difference between the dimension $(N-1)^2$ of
the minimal affine space $\BISTOCHSPACE$ containing all
\BISTOCHASTIC\ matrices, and the dimension of the image
 $\DIFFERENTIAL{f}{}(\TANGENTspace{\UNITARY}{U} )$
of the space
 $\TANGENTspace{\UNITARY}{U}$,
tangent to $\UNITARY$ at $U$, under the tangent map
$\DIFFERENTIAL{f}{}$.
 A non--negative integer number resulting from this
calculation for a  given  unitary matrix $U$
will be called its  {\sl defect}.
We  conjectured the defect to be equal to zero for a generic $U$,
as a non--zero
defect condition has the form of one additional equation imposed on
entries of $U$. As we have been recently informed, 
 the statement that the set of unitary matrices with a non--zero defect 
 has measure zero within the set of all unitaries 
 follows from an early work by Karabegov \cite{Ka89}.

Any non--zero value of the defect
may be considered as a kind of quantification
of the particular structure of $U$. 
    For instance,
the defect is positive if $U$ is an orthogonal matrix of size $N\ge 3$,
if $U$ has a certain {\sl pattern} i.e. some of its entries are equal to zero,
or if  $U$ has a tensor product structure \cite{TK07}.

After the definition of the defect was first proposed in our previous work
\cite{TZ06}, this concept was used in very recent papers
\cite{BN06,Sz06,MS07} to characterize  complex Hadamard matrices.
In this work we prove several properties of the defect,
demonstrating its invariance with respect to the equivalence
relation. We show that vanishing of the defect of $U$ implies that
$U$ is isolated and we find a relation to an analogous 'span
condition' by Nicoara \cite{Nicoara04}.

 The key result of this paper
we regard to be an explicit formula  
for the defect of the Fourier matrix $F_N$ of size $N$.
Equivalent, more transparent forms  
of this formula were obtained by  W. S{\l}omczy{\'n}ski and 
are proved in appendix B.
The defect vanishes iff $N$ is prime,
which implies the earlier statement by Petrescu
\cite{Petrescu97} that the Fourier matrix
is isolated if its dimension is a prime number.
This in turn implies that the flat \BISTOCHASTIC\ matrix $J_N$
belongs to the interior of the set of unistochastic matrices \cite{BEKTZ05},
if $N$ is prime.

For a composite $N$ the defect of $F_N$ is positive,
and it is usually greater than the dimension of
{\sl affine Hadamard families} stemming from $F_N$, introduced in
\cite{TZ06}.
     However, if the size of a matrix is a power of prime,
$N=p^k$, the defect and the dimension
coincide. So, in this very case, an explicit construction
of the defect--dimensional affine family of unitary complex 
Hadamard matrices stemming from the Fourier matrix is complete.
By 'complete' we mean that this solution cannot be embedded
inside any orbit of inequivalent complex Hadamard matrices
of a larger dimension.

This work is organized as follows. 
   In Section \ref{sec_defect}
the definition of the defect of a unitary matrix is provided.
   Several properties of the defect are investigated in 
Section \ref{sec_defect_properties}.
   In Section \ref{sec_defect_applications} we present some applications of the defect
analyzing the condition for a unitary matrix to be isolated
and discussing the unistochasticity problem.
   Section \ref{sec_Fourier_matrix_defect} contains derivation of the formula for the
defect of the Fourier matrix of an arbitrary size $N$
and a discussion of its special cases.
   In Section \ref{sec_defect_dimensional_Fourier_orbits} we provide two
constructions of the defect--dimensional orbit of unitary complex
Hadamard matrices stemming from
$F_N$ for $N$ being a power of a prime number.
   The paper is concluded in Section \ref{sec_conclusions}.

We use in the paper the notation $A \HADprod B$
for the Hadamard product of two matrices,
 $\ELEMENTof{A \HADprod B}{i}{j}=A_{i,j}B_{i,j}$,
while
  $\EXPentrywise(A)$
denotes entrywise exponentiation of a matrix,
  $\ELEMENTof{\EXPentrywise(A)}{i}{j}=\exp(A_{i,j})$.
Also, as functions of matrices are used, to avoid doubts about an
order of variables, for example when writing Jacobi matrices, we
introduce $\VEC,\VECR,\VECC$ notation for appropriate vector forms of
each matrix. Such notations make it possible for us to treat manifolds
of matrices and their tangent spaces as subsets of $\REALS^k$, identified with the
set of all real $k \times 1$ column matrices, and avoid more
abstract constructions. These and other symbols used are listed and explained in
Appendix \ref{sec_notation}.

%
%

\section{The defect of a unitary matrix}
\label{sec_defect}

%
%

\subsection{Definition of the defect}
\label{subsec_defect_definition}

Let $\VECR(\UNITARY)$ be a submanifold of $R^{2N^2}$ representing the
set $\UNITARY$ of all $N \times N$ unitary matrices (for the notation
consult Appendix \ref{sec_notation}).
Consider also the $(N-1)^2$ dimensional minimal hyperplane 
containing \BISTOCHASTIC\ matrices $\VEC(\BISTOCHSPACE) \subset \REALS^{N^2}$,
and a map $f:\ \REALS^{2N^2} \longrightarrow \REALS^{N^2}$,
effectively squaring the moduli of the entries of an $N\times N$ complex matrix
U:
\begin{equation}
  \label{eq_f_def}
  f( \VECR(U) ) = \VEC(B) \mbox{\ \ \ \ where\ \ \ } B_{i,j} = |U_{i,j}|^2 .
\end{equation}

Next consider the tangent map
$\DIFFERENTIAL{f}{\VECR(U)}:\ \REALS^{2N^2} \longrightarrow \REALS^{N^2}$,
realized by the appropriate Jacobi matrix:
\begin{eqnarray}
  \label{eq_Df}
  \lefteqn{
    \DIFFERENTIAL{f}{\VECR(U)}\left( \rule{0cm}{0.3cm} \VECR(V) \right) =
  }  &  &  \\
  &
  2 \cdot
  \left[
    \begin{array}{c|c}
      \diag\left( \rule{0cm}{0.4cm}  \VEC( \Re(U)) \right)
      \rule[-0.3cm]{0cm}{0.8cm}
      &
      \diag\left( \rule{0cm}{0.4cm}  \VEC( \Im(U)) \right)
    \end{array}
  \right]
  \ \ \cdot\ \  \VECR(V).  &  \nonumber \\
\end{eqnarray}

Consider also the tangent spaces, the space $\TANGENTspace{\UNITARY}{U}$
tangent to $\VECR(\UNITARY)$ at $\VECR(U)$ for some unitary $U
\in \UNITARY$, and the space $\TANGENTspace{\BISTOCHSPACE}{B}$  tangent to
$\VEC(\BISTOCHSPACE)$ at $\VEC(B) = f( \VECR(U) )$,
a \BISTOCHASTIC\ matrix.
$\TANGENTspace{\UNITARY}{U}$ here and further is understood as the
nullspace of the Jacobi matrix of the map
$\VECR(W) \longrightarrow \VECR(W^{*} W - I)$ calculated at $\VECR(U)$
(i.e. the kernel of the corresponding tangent map).
$\TANGENTspace{\BISTOCHSPACE}{B}$ is the space of all vectors
$\VEC(G)$ with $G$ being a real $N \times N$ matrix with 
sums of all entries in each row and collumn equal to zero, 
irrespectively of a \BISTOCHASTIC\ $B$.
It is clear that the image of
$\TANGENTspace{\UNITARY}{U}$ under $\DIFFERENTIAL{f}{\VECR(U)}$ must
be contained in $\TANGENTspace{\BISTOCHSPACE}{B}$, so its
dimension is not greater than $(N-1)^2$. It is reduced, with
respect to that value, by a number which will be called {\sl defect} of $U$:


\begin{definition}
  \label{def_defect}

  The defect of an $N \times N$ unitary matrix $U$, denoted
  $\DEFECT(U)$, is the following integer number:
  \begin{equation}
  \label{eq_def_defect}
    \DEFECT(U) =
    (N-1)^2 -
    \dim\left(
      \DIFFERENTIAL{f}{\VECR(U)}( \TANGENTspace{\UNITARY}{U} )
      \rule{0cm}{0.4cm}
    \right).
  \end{equation}
\end{definition}

It is obvious that $\DEFECT(U) = d$ is equivalent to the fact that the
dimension of the part of the nullspace of
$\DIFFERENTIAL{f}{\VECR(U)}$ contained in $\TANGENTspace{\UNITARY}{U}$
is equal to
\begin{eqnarray}
  \label{eq_Df_nulspace_part_within_TU_dim}
  \lefteqn{
    \dim\left(
      \rule{0cm}{0.4cm}
      \NULLR(\DIFFERENTIAL{f}{\VECR(U)}) \cap
      \TANGENTspace{\UNITARY}{U}
    \right) =
  } & & \\
  &
  \dim\left( \rule{0cm}{0.3cm} \TANGENTspace{\UNITARY}{U} \right) -
  \dim\left( \rule{0cm}{0.3cm} \TANGENTspace{\BISTOCHSPACE}{f(\VECR(U))} \right) + d\ \ =  &
  \nonumber \\
  &  N^2 - (N-1)^2 + d\ \ =\ \  2N-1 + d.  & \nonumber
\end{eqnarray}

If $\DEFECT(U)  > 0 $ then $\VECR(U)$ is also called a {\sl
  critical point} of map $f$ restricted to $\UNITARY$.

%
%

\subsection{Other characterizations of the defect}
\label{subsec_other_defect_characterizations}

The tangent space $\TANGENTspace{\UNITARY}{U}$ is equal
to the set:
\begin{eqnarray}
  \label{eq_T_U}
  \left\{ \rule{0cm}{0.4cm}  \VECR(EU):\ E \mbox{\ anti--hermitian}\right\} =  &
  \mbox{(or alternatively)}  &  \\
  \left\{ \rule{0cm}{0.4cm}  \VECR(UF):\ F \mbox{\ anti--hermitian}\right\}
  &   &  \nonumber
\end{eqnarray}
and is spanned by all the independent vectors from the set:
\begin{equation}
  \label{eq_T_U_spanning_vectors}
  \left\{ \rule{0cm}{0.4cm} \VECR(A^{(i,j)} U):\ 1 \leq i   <  j \leq N \right\} \cup
  \left\{ \rule{0cm}{0.4cm} \VECR(S^{(i,j)} U):\ 1 \leq i \leq j \leq N \right\},
\end{equation}
where
\begin{eqnarray}
  \label{eq_Aij_Sij_def}
  A^{(i,j)}_{k,l} =
  \left\{
    \begin{array}{cc}
       1 & \mbox{for\ \ } (k,l)=(i,j) \\
      -1 & \mbox{for\ \ } (k,l)=(j,i) \\
       0 & \mbox{otherwise}
    \end{array}
  \right.
  &       &
  S^{(i,j)}_{k,l} =
  \left\{
    \begin{array}{cc}
      \Ii & \mbox{for\ \ } (k,l)=(i,j) \\
      \Ii & \mbox{for\ \ } (k,l)=(j,i) \\
      0   & \mbox{otherwise}
    \end{array}
  \right. \\
   1 \leq i   <  j \leq N \ \ \ \
   &       &
   \ \ \ \  1 \leq i \leq j \leq N,  \nonumber
\end{eqnarray}
as matrices $A^{(i,j)}$ and $S^{(i,j)}$ span the real space of
anti--hermitian matrices.

Since $\DIFFERENTIAL{f}{\VECR(U)}( \VECR(S^{(i,i)} U) ) = \ZEROvect$ we
consider only ordered pairs $(i,j),\ i<j$, in construction of a matrix
$M$ containing vectors spanning 
$\DIFFERENTIAL{f}{\VECR(U)}( \TANGENTspace{\UNITARY}{U} )$ as its columns.

First, let us construct an $N^2 \times N(N-1)/2$ complex
matrix $M_{\COMPLEX}$ such that it's $\alpha(i,j)$-th column
(see Appendix \ref{sec_notation} for $\alpha(.,.)$) is defined by:
\begin{equation}
  \label{eq_Mc_ij_th_column}
  \ELEMENTof{M_{\COMPLEX}}{1:N^2}{\alpha(i,j)} =
  \VECC( U^{(i,j)} ),
\end{equation}
where $U^{(i,j)}$ is an $N \times N$ complex matrix, with the $i$-th and
$j$-th non--zero rows only, being negations of each other:
\begin{eqnarray}
  \label{eq_Uij_def}
  \ELEMENTof{U^{(i,j)}}{i}{1:N}  &  =  &  U_{i,1:N} \HADprod \CONJ{U}_{j,1:N}, \\
  \ELEMENTof{U^{(i,j)}}{j}{1:N}  &  =  &  - U_{i,1:N} \HADprod \CONJ{U}_{j,1:N}. \nonumber
\end{eqnarray}

Secondly, we form an $N^2 \times N(N-1)$ real matrix $M$,
\begin{equation}
  \label{eq_M_def}
  M =
  \left[
    \begin{array}{c|c}
      \Re(M_{\COMPLEX}) \rule[-0.1cm]{0cm}{0.5cm}  &  \Im(M_{\COMPLEX})
    \end{array}
  \right],
\end{equation}
which has that nice property:
\begin{eqnarray}
  \label{eq_DfTU_spanning_vectors_in_M}
  M_{1:N^2,\alpha(i,j)} =
  \ELEMENTof{\Re(M_{\COMPLEX})}{1:N^2}{\alpha(i,j)}
  =
  \DIFFERENTIAL{f}{\VECR(U)}
    \left( \rule{0cm}{0.4cm} \VECR( A^{(i,j)} U ) \right),
  &  &  \\
  M_{1:N^2, \frac{N(N-1)}{2} + \alpha(i,j)} =
  \ELEMENTof{\Im(M_{\COMPLEX})}{1:N^2}{\alpha(i,j)}
  =
  \DIFFERENTIAL{f}{\VECR(U)}
    \left( \rule{0cm}{0.4cm}  \VECR( S^{(i,j)} U ) \right) \  .
  &  &  \nonumber
\end{eqnarray}
Hence 
$\SPANR(M) = \DIFFERENTIAL{f}{\VECR(U)}( \TANGENTspace{\UNITARY}{U} )$,
and the defect of $U$ can be calculated as
\begin{equation}
  \label{eq_defect_def_with_M}
  \DEFECT(U) = (N-1)^2 - \RANK(M).
\end{equation}

Note also that
\begin{equation}
  \label{eq_null_span_relation}
  \dim( \NULLR(M^T) ) = N^2 - \dim( \SPANR( M^T ) ),
\end{equation}
where
\begin{equation}
  \label{eq_span_defect_relation}
  \dim( \SPANR( M^T ) ) = \dim( \SPANR( M ) ) =
  (N-1)^2 -\DEFECT(U),
\end{equation}
so
\begin{eqnarray}
  \label{eq_M_matrix_defect}
  \DEFECT(U) & = & N^2 - (2N-1) - \dim( \SPANR( M^T ) )  \nonumber \\
             & = & \dim( \NULLR( M^T ) ) - (2N-1). \label{eq_null_def_of_defect}
\end{eqnarray}

The nullspace of $M^T$ is the solution to the real system
\begin{equation}
  \label{eq_Mtransposed_system}
  M^T \cdot \VEC(R) = \ZEROvect
\end{equation}
with respect to a real $N \times N$ matrix variable $R$, which can be
rewritten with the matrix $M_{\COMPLEX}$:
\begin{equation}
  \label{eq_MCtransposed_system}
  M_{\COMPLEX}^T \cdot \VEC(R) = \ZEROvect,
\end{equation}
or explicitly
\begin{equation}
  \label{eq_U_Uconj_Rdiff_system}
  \forall\ \  1 \leq i < j \leq N \ \ \ \ \
  \sum_{k=1}^{N} U_{i,k} {\CONJ{U}}_{j,k} ( R_{i,k} - R_{j,k} )
  = 0.
\end{equation}
System (\ref{eq_U_Uconj_Rdiff_system}) is solved by the $(2N-1)$
dimensional real space spanned by matrices with only one row, or only
one column, filled with $1$'s, the other elements being zeros. If
the real solution space of (\ref{eq_U_Uconj_Rdiff_system}) is not
greater than that, then $\DEFECT(U) = 0$ according to the alternative
definition (\ref{eq_null_def_of_defect}) of the defect.

The solution space of system (\ref{eq_U_Uconj_Rdiff_system}) can also
be expressed as
\begin{eqnarray}
  \label{eq_iRU_equal_EU_nullpsace}
  \left\{
    \rule{0cm}{0.4cm}
    R:\ \Ii R \HADprod U  =  E U \mbox{  for some anti--Hermitian  } E
  \right\} =
  &  \mbox{(or alternatively)}  &  \\
  \left\{
    \rule{0cm}{0.4cm}
    R:\ (\Ii R \HADprod U)U^* \mbox{ is anti--Hermitian}
  \right\},
  &  &  \nonumber
\end{eqnarray}
that is the set of those $R$, for which the direction
$\VECR( \Ii R \HADprod U )$ of the zero first order change of moduli of
matrix $U$ sitting in $\VECR( U )$ belongs to the tangent space
$\TANGENTspace{\UNITARY}{U}$.

Those special $R$'s that solve (\ref{eq_U_Uconj_Rdiff_system}),
give rise, through $R \rightarrow \Ii R \HADprod U$, to matrices
$\Ii \cdot \diag(\STbasis{k}) \cdot U$,\ \  $U \cdot \Ii \cdot \diag(\STbasis{k})$,
which satisfy the equality in the  definition of the set in
(\ref{eq_iRU_equal_EU_nullpsace}). If a matrix  $U$ has no zero entries, like in
the case of unitary complex Hadamard matrices, 
then it spans a $(2N-1)$ dimensional real space, 
which can be represented in the vector form, 
\begin{equation}
  \label{eq_special_iRU_UiR_space_full_dimension}
  \SPANR\left(
    \rule{0cm}{0.5cm}
    \left\{
      \rule{0cm}{0.4cm}
      \VECR( \Ii \cdot \diag(\STbasis{k}) \cdot U ) :\ k=1..N
    \right\}
    \cup
    \left\{
      \rule{0cm}{0.4cm}
      \VECR( U \cdot \Ii \cdot \diag(\STbasis{k}) ) :\ k=1..N
    \right\}
  \right).
\end{equation}
This is due to the fact that for linear combinations we have this
equivalence with the special $R$'s,
$\STbasis{k} \ONESvect^T$, $\ONESvect \STbasis{k}^T$
(for $\ONESvect$ see Appendix \ref{sec_notation}):
\begin{eqnarray}
  \sum_{k=1}^{N} \alpha_k \VECR(\Ii \cdot \diag(\STbasis{k}) \cdot U)  +
  \sum_{l=1}^{N} \beta_l \VECR(U \cdot \Ii \cdot \diag(\STbasis{l}))
    &      =        &  \ZEROvect    \nonumber  \\
    &  \Updownarrow &   \label{eq_lin_comb_equivalence}      \\
  \sum_{k=1}^{N} \alpha_k \cdot \STbasis{k} \ONESvect^T  +
  \sum_{l=1}^{N} \beta_l  \cdot \ONESvect \STbasis{l}^T
    &      =        &  \ZEROvect,   \nonumber
\end{eqnarray}
if $|U_{i,j}| \neq 0$ for $i,j \in \{1..N\}$. Also in this case, the vectors
\begin{equation}
  \label{eq_phasing_tangent_vectors}
\VECR(\Ii \cdot \diag(\STbasis{k}) \cdot U),\ \ k=1..N
\ \ \ \  \mbox{and}\ \ \ \
\VECR(U \cdot \Ii \cdot \diag(\STbasis{l})),\ \ l=2..N
\end{equation}
span the space tangent
at $\VECR(U)$ to a $(2N-1)$ dimensional manifold:
\begin{equation}
  \label{eq_phasing_full_dim_manifold}
  \left\{
    \rule{0cm}{0.5cm}
    \VECR\left(
      \rule{0cm}{0.4cm}
      \diag( \PHASE{\alpha_1},\ldots,\PHASE{\alpha_N} )
      \cdot U \cdot
      \diag( 1,\PHASE{\beta_2},\ldots,\PHASE{\beta_N} )
    \right)
    :\ \alpha_k, \beta_k \in \REALS
  \right\},
\end{equation}
and if $\DEFECT(U)= d > 0$ then these vectors, together with additional
independent vectors $v_1,\ \ldots,\ v_d$, form a basis for the
space
\begin{equation}
  \label{eq_Df_nullspace_part_within_TU}
  \NULLR( \DIFFERENTIAL{f}{\VECR(U)} )  \cap
  \TANGENTspace{\UNITARY}{U}.
\end{equation}
In general, vectors (\ref{eq_phasing_tangent_vectors}) always belong
to the above space, but they may span a space of dimension
smaller than $(2N-1)$ (and not greater, through the $\Uparrow$
implication in (\ref{eq_lin_comb_equivalence}) for $\alpha_k=-\beta_l=1)$.
  Then also the manifold (\ref{eq_phasing_full_dim_manifold}),
obtained from $U$ by the left and the right multiplication of $U$ by unitary
diagonal matrices, will have its dimension reduced. This is the
subject of Lemma \ref{lem_dim_of_U_phasing_manifold} in
Section \ref{sec_defect_properties}.


In section \ref{subsec_Nicoara_result} we are going to apply another
characterization of the defect. New formulae and the ones already
introduced, all of which  will later be used when proving various
properties of the defect,  
are summarized by the following lemma.


\begin{lemma}
  \label{lem_alternative_defect_definitions}

  The defect of an $N \times N$ unitary matrix $U$ can be calculated as
  \begin{eqnarray}
    \DEFECT(U)  &  =  &  \dim( \NULLR( M^T ) ) - (2N-1)   \label{eq_defect_def_with_Mreal_null}  \\
                & = &   (N-1)^2 - \dim( \SPANR(M) )       \label{eq_defect_def_with_Mreal_span}    \\
              &  =  &  \dim( \NULLC( W^T ) ) - (2N-1)   \label{eq_defect_def_with_Wcomplex_null}  \\
                &  =  &  (N-1)^2 - \dim( \SPANC( W ) ),   \label{eq_defect_def_with_Wcomplex_span}
  \end{eqnarray}
  where (with $M_{\COMPLEX}$ of (\ref{eq_Mc_ij_th_column}))
  \begin{equation}
    \label{eq_W_def}
    W =
    \left[
      \begin{array}{c|c}
        M_{\COMPLEX}  \rule[-0.1cm]{0cm}{0.5cm}  &  - \CONJ{M_{\COMPLEX}}
      \end{array}
    \right].
  \end{equation}
\end{lemma}

\PROOFstart
Only the formulas 
(%
 \ref{eq_defect_def_with_Wcomplex_null},
 \ref{eq_defect_def_with_Wcomplex_span}%
) 
need explanation. Note that:
\begin{eqnarray}
  \label{eq_Wtransposed_nullspace_properties}
  v \in \NULLC( W^T ) & \Longrightarrow & \CONJ{v},\ \Re(v),\ \Im(v) \in \NULLC( W^T ), \\
  r \mbox{ real } \in \NULLC( W^T ) & \Longrightarrow & r \in \NULLR( M^T ), \\
  r \in \NULLR( M^T ) & \Longrightarrow & r \in \NULLC( W^T ).
\end{eqnarray}

Let $R_i,\ i=1..N$, denote matrices with the $i$-th row filled
with $1$'s, having $0$'s elsewhere, and let $C_j,\ j=2..N$,
denote matrices with the $j$-th column filled with $1$'s, having
$0$'s elsewhere. Obviously,
$\VEC(R_i),\ \VEC(C_j) \in \NULLR( M^T ),\ \NULLC( W^T )$,
and they are all independent.


Let vectors $v_1,\ \ldots,\ v_d \in \COMPLEX^N$ be such that the set
of complex vectors
\begin{equation}
  \label{eq_vectors_in_nullc_Wtrasposed}
  \left\{\rule{0cm}{0.4cm}  v_l:\ l=1..d\right\} \cup
  \left\{\rule{0cm}{0.4cm}  \VEC(R_i):\ i=1..N\right\} \cup
  \left\{\rule{0cm}{0.4cm}  \VEC(C_j):\ j=2..N\right\}
\end{equation}
is contained in $\NULLC( W^T )$ and consists of independent
vectors. In this reasoning the case when $d=0$, that is when
$\{v_1,\ldots,v_d\}$ is empty, is included.

Then one can choose real vectors
$r_1,\ \ldots,\ r_d \in   \{\Re(v_l),\ \Im(v_l):\ l=1..d \}$
such that the set of real vectors
\begin{equation}
  \label{eq_vectors_in_nullc_Mtrasposed}
  \left\{\rule{0cm}{0.4cm}  r_l:\ l=1..d\right\} \cup
  \left\{\rule{0cm}{0.4cm}  \VEC(R_i):\ i=1..N\right\} \cup
  \left\{\rule{0cm}{0.4cm}  \VEC(C_j):\ j=2..N\right\}
\end{equation}
is contained in $\NULLR( M^T )$ and consists of independent vectors.

  This choice is possible due to the following inclusion relation:
\begin{eqnarray}
  &
  \SPANC\left(
    \rule{0cm}{0.5cm}
    \left\{\rule{0cm}{0.4cm}  v_l:\ l=1..d\right\} \cup
    \left\{\rule{0cm}{0.4cm}  \VEC(R_i):\ i=1..N\right\} \cup
    \left\{\rule{0cm}{0.4cm}  \VEC(C_j):\ j=2..N\right\}
  \right)
  &  \nonumber \\
  &  \bigcap  &  \label{eq_span_of_v_R_C_within_span_of_REv_IMv_R_C}  \\
  &
  \SPANC\left(
    \rule{0cm}{0.5cm}
    \left\{\rule{0cm}{0.4cm}  \Re(v_l),\ \Im(v_l):\ l=1..d\right\} \cup
    \left\{\rule{0cm}{0.4cm}  \VEC(R_i):\ i=1..N\right\} \cup
    \left\{\rule{0cm}{0.4cm}  \VEC(C_j):\ j=2..N\right\}
  \right).
  &  \nonumber
\end{eqnarray}

On the other hand, if we asssume that all the vectors in the set 
(\ref{eq_vectors_in_nullc_Mtrasposed}) are independent and belong to
$\NULLR( M^T )$, then they form an independent set, as complex
vectors, in $\NULLC( W^T )$.


Thus we have come to that:
\begin{equation}
  \label{eq_equal_nullspace_dims}
  \dim\left( \rule{0cm}{0.4cm}  \NULLR( M^T ) \right)  =
  \dim\left( \rule{0cm}{0.4cm}  \NULLC( W^T ) \right).
\end{equation}
\PROOFend

To provide yet another characterization of the defect of $U$, used
later in 
    section \ref{subsec_isolated_unitary_matrices} in the proof of 
    Theorem \ref{theor_zero_defect_isolation},
let us define a function $g:\ \REALS^{N^2} \longrightarrow
\REALS^{N(N-1)}$ 
($g$ will also occur in 
    section \ref{sec_defect_dimensional_Fourier_orbits} in the proof
    of Theorem  \ref{theor_FpTOk_family_only_Hadamards}) :
\begin{eqnarray}
  {\left[
      \rule{0cm}{0.4cm}
      g(\VEC(R))
  \right]}_{\alpha(i,j)} =
  \Re
  \left(
    -\Ii \cdot
    \sum_{k=1}^{N}
    {
      U_{i,k} \CONJ{U}_{j,k} \PHASE{\left( R_{i,k} - R_{j,k} \right)}
    }
  \right)
  &  &  1 \leq i < j \leq N,    \nonumber    \\
  &  &  \label{eq_g_def}                    \\
  {\left[
      \rule{0cm}{0.4cm}
      g(\VEC(R))
  \right]}_{\frac{N(N-1)}{2} + \alpha(i,j)} =
  \Im
  \left(
    -\Ii \cdot
    \sum_{k=1}^{N}
    {
      U_{i,k} \CONJ{U}_{j,k} \PHASE{\left( R_{i,k} - R_{j,k} \right)}
    }
  \right)
  &  &  1 \leq i < j \leq N.   \nonumber
\end{eqnarray}
Note that $g(\VEC(R))=\ZEROvect$ precisely corresponds to the condition that
matrix $U \HADprod \EXPentrywise(\Ii R)$ is unitary. 
At this moment recall that  
           $f-f(U)$ can be interpreted as a function characterizing
           deviations of the moduli of an argument $V$ with respect to the moduli
           of $U$ while moving $V$ along $\UNITARY$. 
On the other hand,
           the function $g$ measures deviation of 
           $U \HADprod \EXPentrywise(\Ii R)$ from unitarity along the
           set of matrices with constant moduli.

The value of the linear map
$\DIFFERENTIAL{g}{\ZEROvect}:\ \REALS^{N^2} \longrightarrow \REALS^{N(N-1)}$,
being the differential of $g$ at $\ZEROvect$, at $\VEC(R)$, is the
vector
\begin{eqnarray}
  {\left[
      \rule{0cm}{0.4cm}
      \DIFFERENTIAL{g}{\ZEROvect}(\VEC(R))
  \right]}_{\alpha(i,j)} =
  \Re
  \left(
    \sum_{k=1}^{N}
    {
      U_{i,k} \CONJ{U}_{j,k} \left( R_{i,k} - R_{j,k} \right)
    }
  \right)
  &  &  1 \leq i < j \leq N,    \nonumber \\
  &  &  \label{eq_Dg_of_vecR}            \\
  {\left[
      \rule{0cm}{0.4cm}
      \DIFFERENTIAL{g}{\ZEROvect}(\VEC(R))
  \right]}_{\frac{N(N-1)}{2} + \alpha(i,j)} =
  \Im
  \left(
    \sum_{k=1}^{N}
    {
      U_{i,k} \CONJ{U}_{j,k} \left( R_{i,k} - R_{j,k} \right)
    }
  \right)
  &  &  1 \leq i < j \leq N.    \nonumber
\end{eqnarray}

The kernel of the differential $\DIFFERENTIAL{g}{\ZEROvect}$
corresponds to space of  solutions of system
(\ref{eq_U_Uconj_Rdiff_system}):
\begin{equation}
  \label{eq_g_differential_nullspace_defect}
  \DEFECT(U) =
  \dim\left(
    \rule{0cm}{0.4cm}
    \NULLR( \DIFFERENTIAL{g}{\ZEROvect} )
  \right) - (2N-1).
\end{equation}

%
%

\section{Properties of the defect}
\label{sec_defect_properties}


\begin{lemma}
  \label{lem_defect_permutation_invariant}

  For any $N \times N$ unitary matrix $U$ and permutation matrices
  $P_r,\ P_c$ :
  \begin{equation}
    \label{eq_defect_permutation invariant}
  \DEFECT( P_r \cdot U \cdot P_c )  =  \DEFECT( U ).
  \end{equation}
\end{lemma}

\PROOFstart
Consider $M_{\COMPLEX}^U$ of (\ref{eq_Mc_ij_th_column}) and
$M^U$ of (\ref{eq_M_def}), constructed
for $U$, and consider also $M_{\COMPLEX}^{UP}$ and $M^{UP}$
constructed for $UP$, where $P$ is a permutation
matrix. Then, for $I$ being the $N \times N$ identity matrix,
\begin{equation}
  \label{eq_Mmatrix_for_UtimesP}
  M_{\COMPLEX}^{UP} = (I \otimes P^T) M_{\COMPLEX}^{U}
  \Longrightarrow
  M^{UP} = (I \otimes P^T) M^{U}
\end{equation}
which, using (\ref{eq_defect_def_with_M}), results in\ \  $\DEFECT(UP)
= \DEFECT(U)$.

Now, let $P$ be given by $P_{i,:} = {\STbasis{\sigma(i)}}^T$, $\sigma$
being a permutation map. Then $M_{\COMPLEX}^{PU}$ of
(\ref{eq_Mc_ij_th_column}) for unitary $PU$ is obtained from
$M_{\COMPLEX}^{U}$ in the following steps:
\begin{itemize}
  \item
        negate and conjugate the $\alpha(i,j)$-th column of
        $M_{\COMPLEX}^{U}$ if $\sigma^{-1}(i) > \sigma^{-1}(j)$,
        for all $1 \leq i < j \leq N$
  \item
        if $\sigma^{-1}(i) < \sigma^{-1}(j)$ shift the
        $\alpha(i,j)$-th column of the result into \\ the
        $\alpha( \sigma^{-1}(i), \sigma^{-1}(j) )$-th position within
        a new result, otherwise shift it into the
        $\alpha( \sigma^{-1}(j), \sigma^{-1}(i) )$-th position,
        for all $1 \leq i < j \leq N$
  \item
        left multiply the result by $P \otimes I$
\end{itemize}
which amounts to permuting and negating columns of $M^{U}$ to get
the corresponding $M^{PU}$. Thus again $\DEFECT(PU) = \DEFECT(U)$.
\PROOFend 


\begin{lemma}
  \label{lem_defect_phasing_invariant}

  For any $N \times N$ unitary matrix $U$ and unitary diagonal matrices
  $D_r,\ D_c$ :
  \begin{equation}
    \label{eq_defect_phasing_invariant}
    \DEFECT( D_r \cdot U \cdot D_c )  =  \DEFECT( U ).
  \end{equation}
\end{lemma}

\PROOFstart
Right multiplication of $U$ by $D_c$ brings no change to $M_{\COMPLEX}$.
Left multiplication by $D_r$ stiffly rotates the chains of  coefficients
\begin{equation}
  \label{eq_coeff_chain}
  \left(\ \
    U_{i,1} \CONJ{U}_{j,1} ,\ \
    U_{i,2} \CONJ{U}_{j,2} ,\ \ \ldots,\ \
    U_{i,N} \CONJ{U}_{j,N}\ \
  \right)
\end{equation}
of system (\ref{eq_U_Uconj_Rdiff_system}), so it does not 
change the space of its solutions nor the value of the defect
(\ref{eq_null_def_of_defect}). (This rotation is equivalent to right
multiplication of each $N^2 \times 2$ sub--matrix of $M$ composed of the
real and imaginary part of some column of $M_{\COMPLEX}$, by a $2
\times 2$ real orthogonal matrix.)
\PROOFend


\begin{lemma}
  \label{lem_defect_THC_invariant}

  For any $N \times N$ unitary matrix $U$
  \begin{equation}
    \label{eq_defect_THC_invariant}
    \DEFECT(U) = \DEFECT(U^T) = \DEFECT(U^*) = \DEFECT(\CONJ{U}).
  \end{equation}
\end{lemma}

{\newcommand{\ANTIHERMITIAN}{\SET{A}}
\PROOFstart 
Since $M_{\COMPLEX}^{\CONJ{U}} = \CONJ{M_{\COMPLEX}^U}$,
for the $M_{\COMPLEX}$ matrices of (\ref{eq_Mc_ij_th_column}) constructed for
$\CONJ{U}$ and $U$,
we have that
$M^{\CONJ{U}}$ is obtained from $M^U$ by negating the right
$N^2 \times N(N-1)/2$\ \ sub--matrix of $M^U$, which leads to
$\DEFECT(\CONJ{U}) = \DEFECT(U)$ by (\ref{eq_defect_def_with_M}).

As for $U^T$, we will show that the set
(\ref{eq_iRU_equal_EU_nullpsace}),
used in characterization of the defect,
constructed either for $U$ or $U^T$, is a linear space of a fixed
dimension.
Let $\ANTIHERMITIAN$ denote the set of all $N \times N$
anti--hermitian matrices, $R$ denotes a real matrix. There holds:
\begin{eqnarray}
  \label{eq_solution_set_equalities}
  \{ R:\ \Ii R \HADprod U^T = E U^T \mbox{  for some  } E \in \ANTIHERMITIAN \}
  = & & \nonumber \\
  \{ R:\ \Ii R^T \HADprod U = (U E^T U^*) U \mbox{  for some  } (U E^T U^*) \in \ANTIHERMITIAN \}
  = & & \nonumber \\
  \{ R:\ \Ii R \HADprod U = E U \mbox{  for some  } E \in \ANTIHERMITIAN \}^T.
\end{eqnarray}
Thus the system (\ref{eq_U_Uconj_Rdiff_system}) 
(equivalently system (\ref{eq_Mtransposed_system}))
solved either for $U$ or for $U^T$ yields the solution space of the
same dimension in both cases.
By (\ref{eq_M_matrix_defect}) then  $\DEFECT(U) = \DEFECT(U^T)$.   
\PROOFend
}

Let us recall the definition of an equivalence class
in the set of  unitary matrices \cite{Ha96,TZ06}.


\begin{definition}
  \label{def_equivalence_relation}

  Two $N \times N$ unitary matrices $U$ and $V$ are $\RELofEQUI$
  \DEFINED{equivalent} if there exist permutation matrices $P_r,\ P_c$ and
  unitary diagonal matrices $D_r,\ D_c$ such that
  \begin{equation}
    \label{eq_equivalence_relation}
    V = P_r D_r \cdot U \cdot D_c P_c.
  \end{equation}
\end{definition}

Lemmas \ref{lem_defect_permutation_invariant} and
\ref{lem_defect_phasing_invariant} imply that
for any two  $\RELofEQUI$ equivalent unitary matrices, 
$V \RELofEQUI U$, their defect is the same,  $\DEFECT(V) = \DEFECT(U)$.
In particular, the defect is constant over the set of all unitary matrices
obtained from $U$ by left and right multiplying it by unitary diagonal
matrices. This set is the image under $\VECR^{-1}$ of what we shall call the phasing
manifold for $U$:


\begin{definition}
  \label{def_phasing_manifold}

  \DEFINED{The phasing manifold} for a unitary $N \times N$ matrix $U$ is the set
  \begin{equation}
    \label{eq_phasing_manifold}
    \left\{
      \rule{0cm}{0.5cm}
      \VECR\left(  \rule{0cm}{0.4cm}  D_r \cdot U \cdot D_c  \right)
      :\ D_r,D_c \mbox{ unitary diagonal}
    \right\}.
  \end{equation}
\end{definition}

The phasing manifold for $U$ is a differentiable manifold. 
Its dimension cannot be greater than $2N-1$, because
any element of (\ref{eq_phasing_manifold}) can be obtained
with $D_c$ having $\ELEMENTof{D_c}{1}{1} = 1$.
More formally,  we have


\begin{lemma}
  \label{lem_dim_of_U_phasing_manifold}

  Let
  \begin{equation}
    \label{eq_special_iRU_UiR_space}
    \SPANR\left(
      \rule{0cm}{0.6cm}
      \left\{
        \rule{0cm}{0.5cm}
        \VECR\left( \rule{0cm}{0.4cm}  \Ii \cdot \diag(\STbasis{k}) \cdot U \right)
        :\ k=1..N
      \right\}
      \cup
      \left\{
        \rule{0cm}{0.5cm}
        \VECR\left( \rule{0cm}{0.4cm}   U \cdot \Ii \cdot \diag(\STbasis{l}) \right)
        :\ l=1..N
      \right\}
    \right)
  \end{equation}
  be spanned by all the vectors from the set of independent vectors
  (where $p+r \leq 2N-1$):
  \begin{equation}
    \label{eq_special_iRU_UiR_space_spanning_vectors}
    \left\{
      \rule{0cm}{0.5cm}
      \VECR\left( \rule{0cm}{0.4cm}  \Ii \cdot \diag(\STbasis{i_k}) \cdot U \right)
      :\ k=1..p
    \right\}
    \cup
    \left\{
      \rule{0cm}{0.5cm}
      \VECR\left( \rule{0cm}{0.4cm}  U \cdot \Ii \cdot \diag(\STbasis{j_l}) \right)
      :\ l=1..r
    \right\}.
  \end{equation}

  Then the set
  \begin{equation}
    \label{eq_U_phasing_manifold_set}
    \left\{
      \rule{0cm}{0.5cm}
      \VECR\left( \rule{0cm}{0.4cm}  D_r \cdot U \cdot D_c \right)
      :\ D_r,D_c \mbox{  unitary diagonal}
    \right\}
  \end{equation}
  is equal to a $(p+r)$ dimensional differential manifold, given by
  the  parametrization:
  \begin{equation}
    \label{eq_U_phasing_manifold_parametrized}
    \left\{
        \rule{0cm}{1.0cm}
        \VECR\left(
           \rule{0cm}{0.7cm}
           \diag\left(
             \rule{0cm}{0.5cm}
             \EXPentrywise\left(
                   \rule{0cm}{0.4cm}
                   \Ii \textstyle{\sum_{k=1}^{p}} \phi_k \STbasis{i_k}
             \right)
           \right)
           \cdot \left[ \rule{0cm}{0.5cm}  U \right] \cdot
           \diag\left(
             \rule{0cm}{0.5cm}
             \EXPentrywise\left(
                   \rule{0cm}{0.4cm}
                   \Ii \textstyle{\sum_{l=1}^{r}} \psi_l \STbasis{j_l}
             \right)
           \right)
        \right):\ \phi_k,\psi_l \in \REALS
    \right\}.
  \end{equation}
\end{lemma} 

\PROOFstart 
First we show that sets (\ref{eq_U_phasing_manifold_set}) and
(\ref{eq_U_phasing_manifold_parametrized}) are equal, i.e. that each
matrix $D_r \ U \ D_c$ can be expressed with $D_r,\ D_c$  satisfying:
\begin{eqnarray}
  \label{eq_Dr_Dc_normalized}
  \ELEMENTof{D_r}{i}{i} = 1  &  \mbox{   for   }  &  i \notin \{i_1 .. i_p\}, \\
  \ELEMENTof{D_c}{j}{j} = 1  &  \mbox{   for   }  &  j \notin \{j_1 .. j_r\}.
\end{eqnarray}

Let
\begin{eqnarray}
  \label{eq_V_Vprim_definition}
  V  &  =  &  \diag( \PHASE{\phi_1},\ldots,\PHASE{\phi_N} )
              \cdot U \cdot
              \diag( \PHASE{\psi_1},\ldots,\PHASE{\psi_N} ),  \\
  V' &  =  &  \diag( \PHASE{\phi'_1},\ldots,\PHASE{\phi'_N} )
              \cdot U \cdot
              \diag( \PHASE{\psi'_1},\ldots,\PHASE{\psi'_N} ),  \nonumber
\end{eqnarray}
where
\begin{equation}
  \begin{array}{ccc}
    \phi'_i = \phi_i + \delta_i  &  \mbox{if}  &  i   \in  \{i_1..i_p\}, \\
    \phi'_i = 0                    &  \mbox{if}  &  i \notin \{i_1..i_p\},
  \end{array}
  \ \ \ \
  \begin{array}{ccc}
    \psi'_j = \psi_j + \epsilon_j  &  \mbox{if}  &  j  \in   \{j_1..j_r\}, \\
    \psi'_j = 0                    &  \mbox{if}  &  j \notin \{j_1..j_r\},
  \end{array}
\end{equation}
and where $\delta_i,\ \epsilon_j$ are uniquely defined by the equation:
\begin{eqnarray}
  \label{eq_delta_epsilon_spanning_definition}
  \sum_{i \in \{i_1..i_p\}}
      \delta_i \cdot \Ii \cdot \diag(\STbasis{i}) \cdot U
  +
  \sum_{j \in \{j_1..j_r\}}
      \epsilon_j \cdot U \cdot \Ii \cdot \diag(\STbasis{j})
  &  =   &   \nonumber \\
  &      &  \label{eq_phi_prim_psi_prim_definition}  \\
  \sum_{i \notin \{i_1..i_p\}}
      \phi_i \cdot \Ii \cdot \diag(\STbasis{i}) \cdot U
  +
  \sum_{j \notin \{j_1..j_r\}}
      \psi_j \cdot U \cdot \Ii \cdot \diag(\STbasis{j}).
  &      &   \nonumber
\end{eqnarray}

Let
\begin{equation}
  \label{eq_V_of_t_definition}
  V(t) =
  \diag\left(
    \rule{0cm}{0.4cm}
    \PHASE{(\phi_1 + t \delta_1)},\ldots,\PHASE{(\phi_N + t \delta_N)}
  \right)
  \cdot \left[ \rule{0cm}{0.4cm}  U \right] \cdot
  \diag\left(
    \rule{0cm}{0.4cm}
    \PHASE{(\psi_1 + t \epsilon_1)},\ldots,\PHASE{(\psi_N + t \epsilon_N)}
  \right),
\end{equation}
where $\delta_i,\ \epsilon_j$ are defined by
(\ref{eq_delta_epsilon_spanning_definition}) for
$i \in \{i_1..i_p\}$,\ \ $j \in \{j_1..j_r\}$, and for the remaining $i,\
j$ we define $\delta_i = -\phi_i$,\ \ $\epsilon_j = -\psi_j$.

Let also $L$ and $R$ denote the left and right hand side of equation
(\ref{eq_delta_epsilon_spanning_definition}), respectively.
Then
\begin{eqnarray}
  \label{eq_V_of_t_derivative}
  \lefteqn{ \frac{\delta}{\delta t} V(t) = }  &    &    \\
  &     &                                    \nonumber  \\
  &
  \diag\left(
    \rule{0cm}{0.4cm}
    \PHASE{(\phi_1 + t \delta_1)},\ldots,\PHASE{(\phi_N + t \delta_N)}
  \right)
  \cdot \left( \rule{0cm}{0.4cm}  L-R \right) \cdot
  \diag\left(
    \rule{0cm}{0.4cm}
    \PHASE{(\psi_1 + t \epsilon_1)},\ldots,\PHASE{(\psi_N + t \epsilon_N)}
  \right)
  &  = \ZEROvect,  \nonumber
\end{eqnarray}
so
\begin{equation}
  \label{eq_V_Vprim_equality}
  V = V(0) = V(1) = V',
\end{equation}
from which it generally follows that sets (\ref{eq_U_phasing_manifold_set}) and
(\ref{eq_U_phasing_manifold_parametrized}) are equal.

The equality $V=V'$  can be verified entry by entry,
using (\ref{eq_delta_epsilon_spanning_definition}), without
derivatives. Then one has to consider cases corresponding to possible
answers to the question whether $i$ and $j$
from the index pair $i,j$ of an entry 
belong or not to $\{i_1..i_p\}$, $\{j_1..j_r\}$, respectively.

To show that (\ref{eq_U_phasing_manifold_parametrized}) is a $(p+r)$
dimensional differential manifold, we need to check that the derivatives of the
vector function in
(\ref{eq_U_phasing_manifold_parametrized}) span a $(p+r)$ dimensional
tangent space. These derivatives are:
\begin{eqnarray}
  \label{eq_U_phasing_manifold_tangent_vectors}
  \left.
    \frac{\delta}{\delta \phi_k} \VECR(...)
  \right|_{\phi_k,\psi_l = 0}  &
  =  &
  \VECR\left(
    \rule{0cm}{0.4cm}
    \Ii \cdot \diag(\STbasis{i_k}) \cdot U
  \right),   \\
  \left.
    \frac{\delta}{\delta \psi_l} \VECR(...)
  \right|_{\phi_k,\psi_l = 0}  &
  =  &
  \VECR\left(
    \rule{0cm}{0.4cm}
    U \cdot \Ii \cdot \diag(\STbasis{j_l})
  \right),  \nonumber
\end{eqnarray}
and they form the $p+r$ element set
(\ref{eq_special_iRU_UiR_space_spanning_vectors}) which is assumed to
consist of independent vectors. 
Thus (\ref{eq_U_phasing_manifold_parametrized}) is a $p+r$ dimensional
manifold around $\VECR(U)$.

Since left and right multiplication of $U$ by unitary diagonal
matrices does not disturb the linear independence of vectors in
(\ref{eq_special_iRU_UiR_space_spanning_vectors}), and since the tangent
space for any other $\VECR(D_r \cdot U \cdot D_c)$ in
(\ref{eq_U_phasing_manifold_parametrized}) is spanned by vectors of the form
(\ref{eq_U_phasing_manifold_tangent_vectors}),
(\ref{eq_U_phasing_manifold_parametrized}) is globally a $(p+r)$
dimensional manifold.
\PROOFend 

We will also provide a lower bound for the defect of a real orthogonal
matrix, as well as a formula for the defect of a direct sum of unitary
matrices.


\begin{lemma}
  \label{lem_orth_and_direct_sum_matrix_defect}
  \begin{description}
    \item[a)]
          If $Q$ is a real $N \times N$ orthogonal matrix, then
          \begin{equation}
            \label{eq_orth_matrix_defect_lower_bound}
            \DEFECT(Q) \geq \frac{(N-1)(N-2)}{2}.
          \end{equation}

    \item[b)]
          If $U$ is an $N \times N$ block diagonal unitary matrix,
          \begin{equation}
            \label{eq_U_direct_sum}
            U\ =\ U_1 \oplus U_2 \oplus \ldots \oplus U_r,
          \end{equation}
          with\ \ $U_1$ of size $N_1$,\ \ ...,\ \ $U_r$ of size
          $N_r$,\ \ then
          \begin{eqnarray}
            \label{eq_direct_sum_matrix_defect}
            \lefteqn{\DEFECT(U) = }  &  &  \\
            &
            \displaystyle{
              \left(
                (N-1)^2 - \sum_{k=1}^r (N_k - 1)^2
              \right)\
              +\
              \sum_{k=1}^r \DEFECT(U_k)
            }  &
            >
            \ \ \sum_{k=1}^r \DEFECT(U_k).  \nonumber
          \end{eqnarray}
          where\ \  $\DEFECT(U_k)=0$\ \  if\ \  $N_k=1$,\ \ according to the
          definition of the defect.
  \end{description}
\end{lemma} 

\PROOFstart 
{\bf a)} The $M_{\COMPLEX}$ matrix of (\ref{eq_Mc_ij_th_column}) constructed
for $Q$,\ \ $M_{\COMPLEX}^{Q}$,\ \  is a real matrix. Then the corresponding
$M$ matrix of (\ref{eq_M_def}) is equal to:
\begin{equation}
  \label{eq_M_matrix_for_Q}
  M^Q\ =\
  \left[
    \begin{array}{c|c}
      M_{\COMPLEX} \rule[-0.1cm]{0cm}{0.5cm}  &  \ZEROvect
    \end{array}
  \right]
\end{equation}
and its rank is not greater then $(N-1)N/2$.
Then by (\ref{eq_defect_def_with_M}):
\begin{equation}
  \label{eq_orth_matrix_defect_lower_bound_calculation}
  \DEFECT(Q)\ \geq\
  (N-1)^2  - \frac{(N-1)N}{2}\ =\
  \frac{(N-1)(N-2)}{2}.
\end{equation}

{\bf b)} The $M_{\COMPLEX}$ matrix for
$U = U_1 \oplus \ldots \oplus U_r$,\ \ $M_{\COMPLEX}^U$\ , can be
permuted to take the form:
\begin{equation}
  \label{eq_Mc_for_U_direct_sum_form}
  M_{\COMPLEX}^{U_1}  \oplus
  M_{\COMPLEX}^{U_2}  \oplus
  \ldots              \oplus
  M_{\COMPLEX}^{U_r}  \oplus  \left[ \rule{0cm}{0.3cm} \ZEROvect \right]
\end{equation}
and consequently $M^U$ can be permuted to become:
\begin{equation}
  \label{eq_M_for_U_direct_sum_form}
  M^{U_1}  \oplus
  M^{U_2}  \oplus
  \ldots   \oplus
  M^{U_r}  \oplus  \left[ \rule{0cm}{0.3cm} \ZEROvect \right],
\end{equation}
where those summands $M^{U_k}$ for which $N_k = 1$ are 'empty'
matrices, that is they do not enter the direct sum.

Then by (\ref{eq_defect_def_with_M}):
\begin{eqnarray}
  \label{eq_direct_sum_matrix_defect_calculation}
  \lefteqn{ \DEFECT(U)\ =\ (N-1)^2 - \RANK( M^U )\ =}  &  &  \\
  &
  \displaystyle{
    (N-1)^2 - \sum_{k=1}^r (N_k - 1)^2 +
    \sum_{k=1}^r
      \left(
        \rule{0cm}{0.4cm}
        (N_k - 1)^2 - \RANK(M^{U_k})
      \right)  =
  }  &   \nonumber \\
  &
  \displaystyle{
    \left(
      (N-1)^2 - \sum_{k=1}^r (N_k - 1)^2
    \right)
    \ \ +\ \
    \sum_{k=1}^r \DEFECT(U_k),
  }  &   \nonumber
\end{eqnarray}
where in the second expression we define $\RANK(M^{U_k})=0$ if $M^{U_k}$
is 'empty'. As $\DEFECT(U_k)=0$ if $N_k=1$, the above formula is also
valid in the case of presence of $1 \times 1$ diagonal blocks in $U$.
\PROOFend 

%
%

\section{Exemplary applications}
\label{sec_defect_applications}

%
%

\subsection{Isolated unitary matrices and continuous families of unitary
  matrices with a fixed pattern of the moduli}
\label{subsec_isolated_unitary_matrices}


\begin{definition}
  \label{def_isolated}

  A $N \times N$ unitary  matrix is called \DEFINED{isolated} if there
  is a neighbourhood $\SET{W}$ around $\VECR(U)$ such that all
  unitaries $V$ with the properties:
  \begin{itemize}
    \item
          $V$ has the same pattern of moduli as $U$, i.e.
          $|V_{i,j}| = |U_{i,j}|$.
    \item
          $\VECR(V) \in \SET{W}$.
  \end{itemize}
  are those given by the intersection
  \begin{equation}
    \label{eq_phased_in_neighborhood}
    \VECR^{-1}(\SET{W})
    \ \ \cap\ \
    \left\{
      \rule{0cm}{0.4cm}
      D_r \cdot U \cdot D_c :\ D_r,\ D_c \mbox{ unitary diagonal }
    \right\}.
  \end{equation}
\end{definition}

A one way criterion for some $U$ being isolated, associated with
calculation of the defect of $U$, is stated as follows:


\begin{theorem}
  \label{theor_zero_defect_isolation}

  If the defect $\DEFECT(U) = 0$,  then matrix $U$ is isolated.
\end{theorem}

\PROOFstart
All the matrices with the same pattern of the moduli  as in $U$ are given by:
\begin{equation}
  \label{eq_U_pattern_matrices}
  U \HADprod \EXPentrywise (\Ii \cdot R),\ \ \ \  \VEC(R) \in \REALS^{N^2},
\end{equation}
and the unitarity condition for them can be expressed as
\begin{equation}
  \label{eq_U_pattern_matrices_unitarity}
  -\Ii \cdot
  \sum_{k=1}^{N}
  {
    U_{i,k} \CONJ{U}_{j,k} \PHASE{\left( R_{i,k} - R_{j,k} \right)}
  }
  = 0,
  \ \ \ \ \ \ \ \
  1 \leq i < j \leq N.
\end{equation}

We can rewrite (\ref{eq_U_pattern_matrices_unitarity}) with the use of
function  $g$ defined in (\ref{eq_g_def}) as:
\begin{equation}
  \label{eq_g_system}
  g(\VEC(R)) = \ZEROvect.
\end{equation}

From the characterization of the defect of $U$ with the kernel of the differential
of $g$ at $\ZEROvect$, see (\ref{eq_g_differential_nullspace_defect}),
we have that condition $\DEFECT(U) = 0$ implies
$\RANK(\DIFFERENTIAL{g}{\ZEROvect}) = N^2 - (2N-1)$.
Then one can choose a subsystem of system (\ref{eq_g_system}),
consisting of  $N^2 - (2N-1)$ equations
\begin{equation}
  \label{eq_g_subsystem}
  \tilde{g}(\VEC(R)) = \ZEROvect
\end{equation}
with the full rank
\begin{equation}
  \label{eq_g_subsystem_rank}
  \dim
  \left(
    \DIFFERENTIAL{\tilde{g}}{\ZEROvect}( \REALS^{N^2} )
  \right)
  = N^2 - (2N-1).
\end{equation}
System (\ref{eq_g_subsystem}) thus defines a $(2N-1)$ dimensional
manifold around $\ZEROvect$.

This must be a $(2N-1)$ dimensional space:
\begin{equation}
  \label{eq_phasing_plane}
  \left\{
    \VEC\left(
      \sum_{k=1}^{N} \alpha_k (\STbasis{k} \ONESvect^T) +
      \sum_{l=2}^{N} \beta_l  (\ONESvect \STbasis{l}^T)
    \right)
    :\
    \alpha_k,\beta_l \in \REALS
  \right\}.
\end{equation}

If $\VECR( U \HADprod \EXPentrywise(\Ii R) )$ is in $\SET{W}$, a small
neighbourhood of $\VECR(U)$, it can be expressed with $\VEC(R)$ in
a certain neighbourhood of $\ZEROvect$.
(The latter neighbourhood can be made sufficiently small by decreasing
the size of $\SET{W}$, for the purpose of the next argument.)

If $U \HADprod \EXPentrywise(\Ii R)$ is unitary, then
$\VEC(R)$ in this neighborhood of $\ZEROvect$ must satisfy system
(\ref{eq_g_system}), hence system (\ref{eq_g_subsystem}), so it must
belong to (\ref{eq_phasing_plane}).
Thus $U \HADprod \EXPentrywise(\Ii R)$ must be of the form:

\begin{equation}
  \label{eq_phased_U}
  \diag(\PHASE{\alpha_1},\ldots,\PHASE{\alpha_N} )
  \cdot U \cdot
  \diag( 1,\PHASE{\beta_2},\ldots,\PHASE{\beta_N} ),
\end{equation}
that is it belongs to the phasing manifold for $U$ (Definition
\ref{def_phasing_manifold}).

%
%
%
%
\PROOFend 
\bigskip

In general, the defect of $U$ allows us to calculate an upper bound
for the dimension of a differential manifold $\SET{F}'$, stemming
from $\VECR(U)$, generated by unitary matrices with the same 
pattern of the moduli as in $U$, if $\SET{F}'$ exists:
\begin{equation}
  \label{eq_Fprim_moduli_pattern_manifold}
  \SET{F}' \subset
  \left\{
    \rule{0cm}{0.5cm}
    \VECR(V) :\
    V \in \UNITARY \mbox{   and   }
    \left(
      \rule{0cm}{0.4cm}
      |V_{i,j}| = |U_{i,j}| \mbox{  for  } i,j \in \{1..N\}
    \right)
  \right\}.
\end{equation}
Such manifolds exist, the phasing manifold (Definition
\ref{def_phasing_manifold}) being a trivial example.

What is even more important for us, we will consider \emph{dephased} manifolds of
this kind. By a dephased manifold we mean a manifold $\SET{F}$ with
the property described in this definition:

\begin{definition}
  \label{def_dephased_manifold}
  A manifold (set) $\SET{F}$, consisting of vector forms $\VECR(V)$
  of unitary matrices $V$ with the moduli of their entries fixed at
  some nonnegative values, is called a \DEFINED{dephased manifold
    (set)} if the condition holds:
  \begin{equation}
    \label{eq_dephased_manifold}
    \left(\ 
      \rule{0cm}{0.4cm}
      \VECR(V) \in \SET{F}  
      \ \ \ \mbox{and}\ \ \ 
      D_r \cdot V \cdot D_c \neq V\ 
    \right)\ \ \ \
    \Longrightarrow\ \ \ \ 
    \VECR(D_r \cdot V \cdot D_c)\ \notin\ \SET{F},
  \end{equation}
  for any unitary diagonal matrices  $D_r,\ D_c$.
\end{definition}

The importance of this subclass of manifolds comes from our
interest in determining all $\RELofEQUI$-inequivalent (see Definition
\ref{def_equivalence_relation}) unitary matrices with the same 
pattern of the moduli, in particular unitary complex Hadamard matrices. This
question is connected to the unistochasticity problem of Section
\ref{subsec_unistochasticity_problem}. And the remark below explains
this importance more precisely.

\begin{remark}
  \label{rem_local_inequivalence_in_dephased_manifold}
  Let $\VECR(U)$ belong to such a dephased manifold $\SET{F}$, as described above. The
  number of different permuted versions of $U$: $P_r \cdot U \cdot P_c$
  is finite, so finite is the number of their images
  $D_r P_r \cdot U \cdot P_c D_c$ (obtained with the use of unitary
  diagonal matrices $D_r,\ D_c$) whose vector forms  sit in
  $\SET{F}$. This is because there can be at most one image for each
  $P_r U P_c$ in $\VECR^{-1}(\SET{F})$.  

  We conclude that there are finitely many vector forms of unitary
  matrices $\RELofEQUI$-equivalent to $U$ in $\SET{F}$,
  and  that $\VECR(U)$ has a neighbourhood in
  $\SET{F}$ in which there are no vector forms of matrices
  $\RELofEQUI$-equivalent to $U$.

  Moreover, there are infinitely many points in this neighbourhood,
  representing pairwise $\RELofEQUI$-inequivalent unitary matrices,
  forming a sequence converging to $\VECR(U)$.
\end{remark}

We need the notion of a dephased matrix:


\begin{definition}
  \label{def_dephased_matrix}

  Let $U$ be an $N \times N$ unitary matrix such that the linear space
  \begin{equation}
    \label{eq_spanned_space}
    \SPANR\left(
      \rule{0cm}{0.6cm}
      \left\{
        \rule{0cm}{0.5cm}
        \VECR\left(
          \rule{0cm}{0.4cm}
          \Ii \cdot \diag(\STbasis{k}) \cdot U
        \right)
        :\ k \in \{1..N\}
      \right\}
      \cup
      \left\{
        \rule{0cm}{0.5cm}
        \VECR\left(
          \rule{0cm}{0.4cm}
          U \cdot \Ii \cdot \diag(\STbasis{l})
        \right)
        :\ l \in \{1..N\}
      \right\}
    \right)
  \end{equation}
  is spanned by $(p+r) \leq 2N-1$ independent vectors from a
  \DEFINED{spanning set}:
  \begin{eqnarray}
    \label{eq_spanning_set}
    \SET{S}
      & = & \left\{
              \rule{0cm}{0.5cm}
              \VECR\left(
                \rule{0cm}{0.4cm}
                \Ii \cdot \diag(\STbasis{i_1}) \cdot U
              \right),
              \ \ldots,\
              \VECR\left(
                \rule{0cm}{0.4cm}
                \Ii \cdot \diag(\STbasis{i_p}) \cdot U
              \right),
            \right. \\
      &   & \left.
              \rule{0cm}{0.5cm}
              \VECR\left(
                \rule{0cm}{0.4cm}
                U \cdot \Ii \cdot \diag(\STbasis{j_1})
              \right),
              \ \ldots,\
              \VECR\left(
                \rule{0cm}{0.4cm}
                U \cdot \Ii \cdot \diag(\STbasis{j_r})
              \right)
            \right\}.  \nonumber
  \end{eqnarray}

  Let $\SET{I}$ be a $(p+r)$ element set of index pairs,
  $\SET{I} \subset \{1..N\} \times \{1..N\}$, called
  a \DEFINED{pattern set} onwards, associated with $U$ and the
  spanning set $\SET{S}$ in such a way that
  \begin{itemize}
    \item
          $(i,j) \in \SET{I}  \Longrightarrow  U_{i,j} \neq 0$
    \item
          $\VEC\left(\rule{0cm}{0.4cm} (\STbasis{i_1} \ONESvect^T) \HADprod F\right)$,\ ...,\
          $\VEC\left(\rule{0cm}{0.4cm} (\STbasis{i_p} \ONESvect^T) \HADprod F\right)$,\\
          $\VEC\left(\rule{0cm}{0.4cm} (\ONESvect \STbasis{j_1}^T) \HADprod F\right)$,\ ...,\
          $\VEC\left(\rule{0cm}{0.4cm} (\ONESvect \STbasis{j_r}^T) \HADprod F\right)$\\
          are independent vectors,
          where $F = \sum_{(i,j) \in \SET{I}} \STbasis{i}
          \STbasis{j}^T$ is a 'filtering matrix'.
  \end{itemize}

  Then an $N \times N$ unitary matrix $V$, with the same  pattern
  of the moduli as in $U$, is called \DEFINED{dephased} with respect to $U$,
  according to the pattern set $\SET{I}$ associated with the spanning set
  $\SET{S}$, if
  \begin{equation}
    \label{eq_dephasing_condition}
    V_{i,j} = U_{i,j} \ \ \ \ \mbox{for any}\ \ \ \  (i,j) \in \SET{I}.
  \end{equation}
\end{definition}

We use the notions introduced above in this lemma:


\begin{lemma}
  \label{lem_phased_is_not_dephased}

  Let $V$ be dephased with respect to $U$, according to a pattern set
  $\SET{I}$ associated with a spanning set $\SET{S}$.

  If $D_r V D_c \neq V$, then $D_r V D_c$ is not dephased (in the same
  manner), for any unitary diagonal matrices $D_r,\ D_c$.
\end{lemma}

\PROOFstart 
Using the procedure applied in the proof of Lemma
\ref{lem_dim_of_U_phasing_manifold}, one can find $D'_r,\ D'_c$ such
that $D'_r V D'_c = D_r V D_c$, and (where $i_k,\ j_l$ characterize
$\SET{S}$, as in Definition \ref{def_dephased_matrix}):
\begin{eqnarray}
  \label{eq_D_prim_property}
  \ELEMENTof{D'_r}{i}{i} = 1 & \mbox{for} & i \notin \{i_1,\ldots,i_p\}, \\
  \ELEMENTof{D'_c}{j}{j} = 1 & \mbox{for} & j \notin \{j_1,\ldots,j_r\}.
  \nonumber
\end{eqnarray}
If $D'_r V D'_c \neq V$, some of the remaining diagonal entries of
$D'_r$: $\ELEMENTof{D'_r}{i_k}{i_k} = \PHASE{\phi_{i_k}}$ and of
$D'_c$: $\ELEMENTof{D'_c}{j_l}{j_l} = \PHASE{\psi_{j_l}}$ must differ
from $1$. Now assume that $D'_r V D'_c$ is also dephased with respect
to $U$, according to the pattern set $\SET{I}$ associated with the
spanning set $\SET{S}$. Then,
for $F = \sum_{(i,j) \in \SET{I}} \STbasis{i} \STbasis{j}^T$:
\begin{eqnarray}
  F \HADprod (D'_r V D'_c) =\
  F \ \HADprod\
  \EXPentrywise\left(
    \rule{0cm}{0.8cm}
    \Ii
    \left(
      \sum_{k=1}^{p} \phi_{i_k} \STbasis{i_k} \ONESvect^T +
      \sum_{l=1}^{r} \psi_{j_l} \ONESvect \STbasis{j_l}^T
    \right)
  \right)
  \ \HADprod\  V  &  &  \nonumber \\
  =
  \EXPentrywise\left(
    \rule{0cm}{0.8cm}
    \Ii
    \left(
      \sum_{k=1}^{p}
        \phi_{i_k}
        \left(\rule{0cm}{0.4cm} (\STbasis{i_k} \ONESvect^T) \HADprod F\right) +
      \sum_{l=1}^{r}
        \psi_{j_l}
        \left(\rule{0cm}{0.4cm} (\ONESvect \STbasis{j_l}^T) \HADprod F\right)
    \right)
  \right)
  \ \HADprod\  (U \HADprod F)
  =\
  U \HADprod F,  &  &  \nonumber \\
    &  &  \label{eq_DVD_dephased_condition}
\end{eqnarray}
where the last equality is the consequence of the assumption that
$D'_r V D'_c$ is also dephased.

Since $U_{i,j} \neq 0$ for $(i,j) \in \SET{I}$, the respective phases
must be equal to zero:
\begin{equation}
  \label{eq_dephased_condition_on_DVD_phases}
  \sum_{k=1}^{p}
    \phi_{i_k}
    \left(\rule{0cm}{0.4cm} (\STbasis{i_k} \ONESvect^T) \HADprod F\right) +
  \sum_{l=1}^{r}
    \psi_{j_l}
    \left(\rule{0cm}{0.4cm} (\ONESvect \STbasis{j_l}^T) \HADprod F\right)
  = \ZEROvect
\end{equation}
and, as the vector forms of matrices standing in combination
(\ref{eq_dephased_condition_on_DVD_phases}) are independent
(the property of the pattern set $\SET{I}$,
 according to which $D'_r V D'_c$ is dephased),
all $\phi_{i_k},\ \psi_{j_l}$ are equal to zero, which
contradicts that $D'_r V D'_c \neq V$.
\PROOFend 

\bigskip
We will further consider manifolds (stemming from $\VECR(U)$) of
vector forms of matrices $V$ dephased with respect to $U$ in a chosen
way. The above lemma implies that such manifolds are dephased
in the sense of Definition \ref{def_dephased_manifold}. 


\begin{theorem}
  \label{theor_manifold_dim_defect_bound}

  Let $U$ be an $N \times N$ unitary matrix and let $\SET{F}$ be a
  differential manifold in $\REALS^{2N^2}$ stemming
  from $\VECR(U)$, generated, through $V \rightarrow \VECR(V)$, purely by unitary
  matrices $V$ with the same  pattern of the moduli as in $U$, and dephased
  with respect to $U$ according to a pattern set $\SET{I}$ associated
  with a spanning set $\SET{S}$.

  Then
  \begin{equation}
    \label{eq_manifold_dim_defect_bound}
    \dim \SET{F} =
    \dim\left(
        \rule{0cm}{0.4cm}
        \TANGENTspace{\SET{F}}{U}
    \right)
    \ \ \ \  
    \leq
    \ \ \ \ 
    \DMPbound(U)
  \end{equation}
  
  where
  \begin{description}
    \item[$\TANGENTspace{\SET{F}}{U} \stackrel{def}{=}
           \TANGENTspace{\SET{F}}{\VECR(U)}$]
           the space tangent to $\SET{F}$ at $\VECR(U)$,

    \item[$\DMPbound(U)=$]
          $\DEFECT(U) + (2N-1) - \NUMBERofIN{U}{0} -
          \NUMBERofIN{\SET{S}}{}$,\ \ \  where
          \begin{description}
              \item[$\NUMBERofIN{U}{0}$]
                    the number of zero entries in $U$,
              \item[$\NUMBERofIN{\SET{S}}{}$]
                    the number of elements of the spanning set $\SET{S}$,
                    equal to $p+r \leq 2N-1$,
                    where $p,\ r$ bear the same meaning as in Definition
                    \ref{def_dephased_matrix}.
          \end{description}
  \end{description}
\end{theorem}

{
  \newcommand{\SETofZEROS}{\tilde{\SET{I}}} %
\PROOFstart 
Let $v \in \TANGENTspace{\SET{F}}{U}$, that is $v = \gamma'(0)$ for
some smooth curve $\gamma(t) \subset \SET{F}$ such that $\gamma(0) =
\VECR(U)$. Since $\VECR^{-1}(\gamma) \subset \UNITARY$,
$v \in \TANGENTspace{\UNITARY}{U}$ i.e. it satisfies
\begin{equation}
  \label{eq_v_tangent_to_vecrUNITARY}
  v = \VECR(E \cdot U)\ \ \ \ \mbox{for some anti--hermitian}\ \  E.
\end{equation}
Since the moduli of the entries of a matrix do not change over
$\VECR^{-1}(\gamma(t))$, $v \in \NULLR(\DIFFERENTIAL{f}{\VECR(U)})$ of
(\ref{eq_Df}), and in particular zero entries stay intact, so $v$
satisfies also:
\begin{equation}
  \label{eq_v_zero_mouli_change_direction}
  v = \VECR(\Ii R \HADprod U)\ \ \ \   \mbox{for some real matrix}\ \  R.
\end{equation}
Thus $v$ belongs to a space parametrized by the solution space of
(\ref{eq_U_Uconj_Rdiff_system})
(or equivalently (\ref{eq_Mtransposed_system})), namely:
\begin{equation}
  \label{eq_v_in_B_space}
  v\ \  \in\ \
  \SET{D}\  =\
  \left\{
    \rule{0cm}{0.4cm}
    \VECR(\Ii R \HADprod U):\ \VEC(R) \in \SET{R}
  \right\},
\end{equation}
where
\begin{equation}
  \label{eq_B_parametrizing_space}
  \SET{R}\ =\
  \left\{
    \rule{0cm}{0.4cm}
    \VEC(R):\ \Ii R \HADprod U = E U \mbox{  for some anti--hermitian } E
  \right\}.
\end{equation}

Let $\SETofZEROS \subset \{1..N\} \times \{1..N\}$ be such that
$(\imath,\jmath) \in \SETofZEROS \Leftrightarrow U_{\imath,\jmath}=0$.
Because of potential zeros in $U$, we can reduce the parametrizing
space $\SET{R}$ of (\ref{eq_B_parametrizing_space}):
\begin{equation}
  \label{eq_B_newly_parametrized}
  \SET{D} =
  \left\{
    \rule{0cm}{0.4cm}
    \VECR(\Ii R \HADprod U):\ \VEC(R) \in \SET{R}'
  \right\},
\end{equation}
where
\begin{equation}
  \label{eq_B_new_parametrizing_space}
  \SET{R}' =
  \SET{R} \cap
  \left\{
      \VEC(R):\ \forall (\imath,\jmath) \in \SETofZEROS\ \
      R_{\imath,\jmath}=0
  \right\},
\end{equation}
and since
\begin{equation}
  \label{eq_Rparam_ofB_direct_sum}
  \SET{R}\ \ \ \ =\ \ \ \
  \bigoplus_{(\imath,\jmath) \in \SETofZEROS}
    \left\{
      \rule{0cm}{0.5cm}
      \VEC\left(
        \rule{0cm}{0.4cm}
        \alpha \cdot \STbasis{\imath} \STbasis{\jmath}^T
      \right)
      :\ \alpha \in \REALS
    \right\}
  \ \ \ \oplus\ \ \ \SET{R}',
\end{equation}
we obtain a bound for the dimension of $\SET{D}$ of
(\ref{eq_v_in_B_space}) or (\ref{eq_B_newly_parametrized}), using the
characterization (\ref{eq_M_matrix_defect}) of the defect of $U$,
stated also by Lemma
\ref{lem_alternative_defect_definitions}\ :
\begin{eqnarray}
  \label{eq_B_dim_bound}
  \dim(\SET{D})  &  \leq  &
  \dim(\SET{R}') = \dim(\SET{R}) - \NUMBERofIN{\SETofZEROS}{}  \\
  &  &  = \DEFECT(U) + (2N-1) - \NUMBERofIN{U}{0},  \nonumber
\end{eqnarray}
where $\SET{R} = \NULLR(M^T)$ with $M$ of (\ref{eq_M_def}).

Further, for 'filtering matrices'
\begin{eqnarray}
  F  &  =  &  \sum_{(i,j) \in \SET{I}}
                \STbasis{i} \STbasis{j}^T,
                \label{eq_F_filtering_matrix} \\
  G  &  =  &  \ONESvect \ONESvect^T - \sum_{(\imath,\jmath) \in \SETofZEROS}
                                 \STbasis{\imath} \STbasis{\jmath}^T,
                \label{eq_G_filtering_matrix}
\end{eqnarray}
and for all the matrices (see the description of a spanning set in
Definition \ref{def_dephased_matrix})
\begin{eqnarray}
  \Ii (\STbasis{i_k} \ONESvect^T) \HADprod U,  &  k\ =\ 1..p,  &  \\
  \label{eq_spanning_set_row_phasing_tangent_matrices}
  \Ii (\ONESvect \STbasis{j_l}^T) \HADprod U,  &  l\ =\ 1..r,  &
  \label{eq_spanning_set_column_phasing_tangent_matrices}
\end{eqnarray}
in
$\VECR^{-1}( \SET{S} )$, vectors
\begin{eqnarray}
  \label{eq_G_filtered_spanning_vectors}
  \VECR\left(
    \rule{0cm}{0.4cm}
    \Ii (\STbasis{i_k} \ONESvect^T) \HADprod G \HADprod U
  \right),
  &  k\ =\ 1..p,  &  \\
  \VECR\left(
    \rule{0cm}{0.4cm}
    \Ii (\ONESvect \STbasis{j_l}^T) \HADprod G \HADprod U
  \right),
  &  l\ =\ 1..r,  &    \nonumber
\end{eqnarray}
are still independent.

As the considered manifold $\SET{F}$
is composed of $\VECR(V)$ with $V$ dephased with respect
to $U$, that is with non--zero entries $V_{i,j}$, for $(i,j) \in \SET{I}$,
fixed, a non--zero $v \in \TANGENTspace{\SET{F}}{U}$ must not belong to a
$(p+r)$-dimensional subspace of $\SET{D}$ defined with the use of
basis vectors (\ref{eq_G_filtered_spanning_vectors}) by
(we parametrize
$\SET{D}$ with $\SET{R}'$ as in (\ref{eq_B_newly_parametrized}),
hence we use $G$ in the formula below):
\begin{equation}
  \label{eq_Bprim_subspace_of_B}
  \SET{D}' =
  \left\{
    \rule{0cm}{0.8cm}
    \VECR\left(
      \sum_{k=1}^{p}
        \alpha_k \cdot \Ii (\STbasis{i_k} \ONESvect^T) \HADprod G
        \HADprod U
      +
      \sum_{l=1}^{r}
        \beta_l \cdot \Ii (\ONESvect \STbasis{j_l}^T) \HADprod G
        \HADprod U
    \right)
    :\ \alpha_k,\ \beta_l \in \REALS
  \right\}.
\end{equation}

If $v$ were in $\SET{D}'$, that is if it were a combination like that
in (\ref{eq_Bprim_subspace_of_B}), then the dephasing condition would
force for this tangent vector that
\begin{equation}
  \label{eq_v_of_TF_dephasing_consequence}
  \ZEROvect = F \HADprod \VECR^{-1}(v) =
  \Ii \left(
    \sum_{k=1}^{p}
        \alpha_k \cdot (\STbasis{i_k} \ONESvect^T) \HADprod F
    +
    \sum_{l=1}^{r}
        \beta_l \cdot  (\ONESvect \STbasis{j_l}^T) \HADprod F
   \right)
   \HADprod ( U \HADprod F )
\end{equation}
implying $\alpha_k = 0,\  \beta_l = 0$, because vector forms of matrices
standing in the last combination are independent, being a requirement for
the proper choice of a pattern set $\SET{I}$ in Definition
\ref{def_dephased_matrix}.


$\TANGENTspace{\SET{F}}{U}$ is thus bound to be contained in some
space $\SET{D}''$ such that $\SET{D} = \SET{D}' \oplus \SET{D}''$.
The dimension of $\SET{D}''$ reads, using (\ref{eq_B_dim_bound}):
\begin{eqnarray}
  \label{eq_Bbis_dim_bound}
  \dim(\SET{D}'')  &  =   & \dim(\SET{D}) - \dim(\SET{D}') \\
                   &  =   & \dim(\SET{D}) - (p+r)          \\
                   & \leq & \DEFECT(U) + (2N-1) - \NUMBERofIN{U}{0} -
                            \NUMBERofIN{S}{},  \nonumber
\end{eqnarray}
which completes the proof.

Note that if $\SET{F}(...)$ is a parametrization of $\SET{F}$ around
$\VECR(U)$, $\SET{F}$ having the properties stated in
Theorem \ref{theor_manifold_dim_defect_bound}, then
\begin{equation}
  \label{eq_Fprim_parametrization}
  \SET{F}'(...) =
  \VECR\left(
    \rule{0cm}{1cm}
    \diag\left(
      \EXPentrywise\left(
        \Ii \sum_{k=1}^{p}
            \phi_k \cdot \STbasis{i_k}
      \right)
    \right)
    \cdot
    \VECR^{-1}\left( \rule{0cm}{0.4cm} \SET{F}(...) \right)
    \cdot
    \diag\left(
      \EXPentrywise\left(
        \Ii \sum_{l=1}^{r}
            \psi_l \cdot \STbasis{j_l}
      \right)
    \right)
  \right)
\end{equation}
parametrizes a $(p+r) + \dim \SET{F}$ dimensional manifold $\SET{F}'$ around
$\VECR(U)$. The additional independent vectors
(\ref{eq_G_filtered_spanning_vectors})  spanning 
(together with a basis of $\TANGENTspace{\SET{F}}{U}$) the space 
$\TANGENTspace{\SET{F}'}{U}$ can  obtained by differentiating
(\ref{eq_Fprim_parametrization}) with respect to $\phi_k,\ \psi_l$.
 That is to say,
$\TANGENTspace{\SET{F}'}{U} =
    \SET{D}' \oplus \TANGENTspace{\SET{F}}{U}$.

\PROOFend  
} 



Note that the bound $\DMPbound(U)$ defined in 
Theorem \ref{theor_manifold_dim_defect_bound}
is independent of the choice of a
spanning set for $U$. Also, it is natural to suppose that for $U$ having a block diagonal
structure $U = U_1 \oplus \ldots \oplus U_r$ this bound could be the
sum of the bounds calculated for its diagonal components. This rule of
a total bound, not necessarily our $\DMPbound(U)$, being the sum of some
bounds for $U_p$, applies to a very special
construction of a manifold $\SET{F}$ in
Theorem \ref{theor_manifold_dim_defect_bound}, in which the direct sum
of matrix forms of respective parametrized dephased manifolds
constructed for $U_1$, ..., $U_r$ is taken to get the matrix form of
$\SET{F}$:
\begin{equation}
  \label{eq_Fmanifold_direct_sum_construction}
  \VECR^{-1}(\SET{F}) =
     \VECR^{-1}(\SET{F}_1) \oplus \ldots \oplus \VECR^{-1}(\SET{F}_r),
\end{equation}
and it is because the dimensions of the component manifolds add
up to the dimension of $\SET{F}$.

In fact, this rule holds for the quantity  $\DMPbound(U)$
defined in theorem   \ref{theor_manifold_dim_defect_bound}.


\begin{lemma}
  \label{lem_U_direct_sum_DMPbound}
  Let $U = U_1 \oplus U_2 \oplus \ldots \oplus U_r$ be a block
  diagonal unitary matrix of size $N$, where $N_p$ denotes the size of
  $U_p$. Then
  \begin{equation}
    \label{eq_U_direct_sum_DMPbound}
    \DMPbound(U) =
       \DMPbound(U_1) + \DMPbound(U_2) + \ldots + \DMPbound(U_r)
  \end{equation}
\end{lemma}

\PROOFstart 
Let $\SET{S}_p$ be a spanning set for $U_p$. Let us construct a
spanning set $\SET{S}$ for $U$ using the rules
\begin{eqnarray}
  \VECR\left(
    \rule{0cm}{0.4cm}
    \Ii \cdot \diag(\STbasis{k}) \cdot U_p
  \right)
  \in \SET{S}_p
  &  \Longrightarrow  &   
  \VECR\left(
    \rule{0cm}{0.5cm}
    \Ii \cdot \diag\left(\STbasis{\sum_{m=1}^{p-1} N_m + k}\right) \cdot U
  \right)
  \in \SET{S},  \nonumber  \\
&  &  \label{eq_setS_building_rules}  \\
  \VECR\left(
    \rule{0cm}{0.4cm}
    U_p \cdot \Ii \cdot \diag(\STbasis{l})
  \right)
  \in \SET{S}_p
  &  \Longrightarrow  &    
  \VECR\left(
    \rule{0cm}{0.5cm}
    U \cdot \Ii \cdot \diag\left(\STbasis{\sum_{m=1}^{p-1} N_m + l}\right)
  \right)
  \in \SET{S},    \nonumber
\end{eqnarray}
and let every element of $\SET{S}$ be put into it in this way. Thus
$\SET{S}$ is properly constructed and it is clear that every spanning
set for $U$ must be created in this manner.

Though it is not a part of the proof, let us mention that a pattern
set $\SET{I}$ associated with the set $\SET{S}$ must have all its
elements put into it using the rule:
\begin{equation}
  \label{eq_setI_building_rule}
  (i,j) \in \SET{I}_p
  \ \ \ \Longrightarrow \ \ \ \ 
  \left(
    \sum_{m=1}^{p-1} N_m \  +\ i,\ \
    \sum_{m=1}^{p-1} N_m \  +\ j
  \right)
  \ \ \in\ \ \SET{I},
\end{equation}
where $\SET{I}_p$'s are some pattern sets associated with the
sets $\SET{S}_p$  used in the construction of $\SET{S}$.

Therefore,
\begin{equation}
  \label{eq_S_number_of_elements}
  \NUMBERofIN{\SET{S}}{} =
                           \NUMBERofIN{\SET{S}_1}{} +
                           \NUMBERofIN{\SET{S}_2}{} +
                           \ldots                 +
                           \NUMBERofIN{\SET{S}_r}{},
\end{equation}
where again $\NUMBERofIN{\SET{S}_p}{}$ stands for the number of
elements in $\SET{S}_p$.

Using the above equality as well as formula (\ref{eq_direct_sum_matrix_defect})
in Lemma \ref{lem_orth_and_direct_sum_matrix_defect} {\bf b)} we find
that:
\begin{eqnarray}
  \label{eq_U_direct_sum_DMPbound_calculation}
  \lefteqn{\DMPbound(U) =} & & \\
  &
  \DEFECT(U)\ \  +\ \  (2N-1)\ \ -\ \ \NUMBERofIN{U}{0}\ \ -\ \
  \NUMBERofIN{\SET{S}}{}
  &  =  \nonumber  \\    
  &  &  \nonumber \\
  &
  \displaystyle{
    \left(
      (N-1)^2 - \sum_{p=1}^r (N_p - 1)^2
    \right)\
    +\
    \sum_{p=1}^r \DEFECT(U_p)
    \ \ +\ \
    \left( 2\sum_{p=1}^r N_p\ \ -\ \ 1\right)
  }
  &     \nonumber  \\
  &
  \displaystyle{
  -\ \
  \left(
    N^2\ \  -\ \ \sum_{p=1}^r N_p^2\ \ +\ \
    \sum_{p=1}^r \NUMBERofIN{U_p}{0}
  \right)
  \ \ -\ \
  \sum_{p=1}^r \NUMBERofIN{\SET{S}_p}{}
  }
  &  =  \nonumber    \\   
  &  &  \nonumber \\
  &
  \displaystyle{
  \sum_{p=1}^r
  \left(
    \rule{0cm}{0.6cm}
    \DEFECT(U_p)\ \  +\ \  (2N_p-1)\ \
    -\ \ \NUMBERofIN{U_p}{0}\ \
    -\ \ \NUMBERofIN{\SET{S}_p}{}
  \right)
  }
  &  \nonumber   \\
  &
  \displaystyle{
  \ \ \ +\ \ \ (N-1)^2 - \sum_{p=1}^r (N_p - 1)^2\ \ \
  +\ \ \ r - 1\ \ \
  -\ \ \  N^2 + \sum_{p=1}^r N_p^2
  }
  & =\ \ \ \  \sum_{p=1}^r \DMPbound(U_p).  \nonumber
\end{eqnarray}

\PROOFend    

%
%

\subsection{Relation to results of  Nicoara}
\label{subsec_Nicoara_result}

In a  paper \cite{Nicoara04} on commuting squares of 
von~Neumann algebras Nicoara  introduced the 'span condition' for such a square to
be isolated.
Consider the simple case of a commuting square of orthogonal
maximal abelian *--subalgebras of the algebra $\SET{M}_N (\COMPLEX)$ of
complex $N \times N$ matrices:
\begin{eqnarray}
  \DIAGONAL           &  \subset  &  \SET{M}_N (\COMPLEX)                       \nonumber  \\
  \cup              &           &  \cup
                                             \ \ \ \ \ \ \ \ ,\label{eq_commuting_square}  \\
  \COMPLEX \cdot I  &  \subset  & U^* \DIAGONAL U                               \nonumber
\end{eqnarray}
where $\COMPLEX \cdot I$ is the algebra of all $N \times N$ scalar matrices, $\DIAGONAL$
the algebra of all $N \times N$ diagonal matrices, $U$ a unitary complex Hadamard
matrix, i.e. $|U_{i,j}|^2 = 1/N$. Any abelian *--subalgebra, as closed
with respect to the hermitian transposition ${(...)}^*$, is unitarily
diagonalizable, and if it is maximal, then it is diagonalizable into
$\DIAGONAL$.
The property that a commuting square (\ref{eq_commuting_square}) is isolated
is equivalent to $U$ being isolated in accordance with Definition \ref{def_isolated}. The span
condition in this case reads:


\begin{lemma}
  \label{lem_span_condition}

  A unitary complex Hadamard matrix $U$ is isolated if
  \begin{equation}
    \label{eq_span_condition}
    \dim\left(
      \rule{0cm}{0.3cm}
      [ \DIAGONAL,\  U^* \DIAGONAL U ]
    \right)
    \ \  =\ \  (N-1)^2,
  \end{equation}
  where
  \begin{equation}
    \label{eq_set_of_commutators}
    [ \DIAGONAL,\  U^* \DIAGONAL U ]
    \ \stackrel{def}{=}\
    \SPANC\left(
      \rule{0cm}{0.5cm}
      \left\{
        \rule{0cm}{0.4cm}
        \VECC\left(
          \rule{0cm}{0.3cm}
          D_1 \cdot U^* D_2 U \  - \   U^* D_2 U \cdot D_1
        \right)
        :\ \  D_1,\ D_2 \in \DIAGONAL
      \right\}
    \right).
  \end{equation}
\end{lemma}

Condition (\ref{eq_span_condition}) is equivalent to
\begin{equation}
  \label{eq_span_condition_comparable}
  \dim\left(
    \rule{0cm}{0.7cm}
    \SPANC\left(
      \rule{0cm}{0.6cm}
      \left\{
        \rule{0cm}{0.5cm}
        \VECC\left(
          \rule{0cm}{0.4cm}
          B^{(i,j)}
        \right)
        :\ i,j \in \{1..N\}
      \right\}
    \right)
  \right)
  \ \  =\ \  (N-1)^2, \ \ \ \ \ \ \ \ \ \ \ \
\end{equation}
where $B^{(i,j)}$ is an $N \times N$ matrix filled all with $0$'s
except for the $i$-th row and the $i$-th column:
\begin{eqnarray}
  B^{(i,j)}_{i,1:N}  &  =  &
  \left[\ \
    \rule{0cm}{0.4cm}
    \left( U_{j,1} \CONJ{U}_{j,i} \right),\ \ \ldots,\ \
    \left( U_{j,i-1} \CONJ{U}_{j,i} \right),
    \ \ 0,\ \
    \left( U_{j,i+1} \CONJ{U}_{j,i} \right),\ \ \ldots,\ \
    \left( U_{j,N} \CONJ{U}_{j,i} \right)\ \
  \right],                                  \nonumber   \\
  B^{(i,j)}_{1:N,i}  &  =  &
  -\left( B^{(i,j)}_{i,1:N} \right)^*,      \label{eq_Bij_ith_row_and_column}
\end{eqnarray}
as for\ \ $D_1 = \diag(\alpha_1,\ldots,\alpha_N)$,\ \ $D_2 =
\diag(\beta_1,\ldots,\beta_N)$,
\begin{equation}
  \label{eq_commutator_from_the_set}
  D_1 \cdot U^* D_2 U  \ - \  U^* D_2 U \cdot D_1
  \ \ \ = \ \
  \sum_{i,j \in \{1..N\}}  \alpha_i \beta_j \cdot B^{(i,j)}.
\end{equation}

We will show that the sufficient condition (\ref{eq_span_condition})
is equivalent to our condition $\DEFECT(U)=0$ for $U$ being
isolated. That is, using also Lemma \ref{lem_defect_THC_invariant}, that
the equivalence holds:
\begin{equation}
  \label{eq_zero_defect_and_span_equivalence}
  \DEFECT(U) =0      \ \ \ \Longleftrightarrow\ \ \
  \DEFECT(U^*) = 0   \ \ \ \Longleftrightarrow\ \ \
  (\ref{eq_span_condition_comparable}).
\end{equation}

To show this, take $U^*$ and form matrix $W$ of
Lemma \ref{lem_alternative_defect_definitions}
for $U^*$:
\begin{equation}
  \label{eq_Ustar_Wmatrix}
  W^{U^*} =
  \left[
    \begin{array}{c|c}
      M^{U^*}_{\COMPLEX}  &  -\CONJ{M}^{U^*}_{\COMPLEX}
    \end{array}
  \right],
\end{equation}
then concatenate it horizontally with an $N^2\times N$ matrix filled only with $0$'s, and
reorder the columns of the resulting matrix to obtain a square matrix $B$ having the
property that it's $k$-th\ \  $N^2 \times N$ sub--matrix, $k=1..N$, is
equal to:
\begin{eqnarray}
  \label{eq_B_submatrix}
  \lefteqn{B_{1:N^2,(k-1)N+1:kN}\ \ \  =} & & \\
  &
  \left[
    \rule{0cm}{0.5cm}
    \begin{array}{c|c|c|c|c|c|c}
      \VECC\left(
        \rule{0cm}{0.4cm}
        U^{*(k,1)}
      \right)                            \rule{0cm}{0.5cm} &
      \ldots                             \rule{0cm}{0.5cm} &
      \VECC\left(
        \rule{0cm}{0.4cm}
        U^{*(k,k-1)}
      \right)                            \rule{0cm}{0.5cm} &
      \ZEROvect                          \rule{0cm}{0.5cm} &
      \VECC\left(
        \rule{0cm}{0.4cm}
        U^{*(k,k+1)}
      \right)                            \rule{0cm}{0.5cm} &
      \ldots                             \rule{0cm}{0.5cm} &
      \VECC\left(
        \rule{0cm}{0.4cm}
        U^{*(k,N)}
      \right)                            \rule{0cm}{0.5cm}
    \end{array}
  \right], &       \nonumber
\end{eqnarray}
where $U^{(i,j)}$ is defined by (\ref{eq_Uij_def}) for $i<j$, and we
additionally define:
\begin{equation}
  \label{eq_Uji_def}
  U^{(j,i)} \ \ =\ \ - \CONJ{U^{(i,j)}}  \ \ \ \ \ \ \mbox{for}\ \ \   j > i.
\end{equation}

Then the rows of $B$ correspond to matrices $B^{(i,j)}$\ :
\begin{equation}
  \label{eq_Bij_B_relation}
  B^{(i,j)}
  =
  \VECC^{-1}\left(
    \rule{0cm}{0.5cm}
    \left(
      \rule{0cm}{0.4cm}
      B_{(i-1)N+j, 1:N^2}
    \right)^T
  \right),
\end{equation}
and thus, also by Lemma \ref{lem_alternative_defect_definitions}:
\begin{eqnarray}
  \lefteqn{
    \dim\left(
      \rule{0cm}{0.7cm}
      \SPANC\left(
        \rule{0cm}{0.5cm}
        \left\{
          \rule{0cm}{0.4cm}
          \VECC\left(
            \rule{0cm}{0.3cm}
            B^{(i,j)}
          \right)
          :\ i,j \in \{1..N\}
        \right\}
      \right)
    \right)\ \ =
  }
  &  &                                                                   \nonumber \\
  &
  \dim\left(
    \rule{0cm}{0.4cm}
    \SPANC\left(
      \rule{0cm}{0.3cm}
      B^T
    \right)
  \right)
  \ \ =\ \
  \dim\left(
    \rule{0cm}{0.4cm}
    \SPANC\left(
      \rule{0cm}{0.3cm}
      B
    \right)
  \right) & = \nonumber\\
  &
  \dim\left(
    \rule{0cm}{0.4cm}
    \SPANC\left(
      \rule{0cm}{0.3cm}
      W^{U^*}
    \right)
  \right)
  \ \ =\ \
  (N-1)^2 - \DEFECT(U^*) & =  \nonumber \\
  & (N-1)^2 - \DEFECT(U). &    \label{eq_dim_equalities}
\end{eqnarray}

From the above (\ref{eq_zero_defect_and_span_equivalence}) immediately
follows. We can formulate yet another characterization of the defect:
\begin{equation}
  \label{eq_defect_def_with_Bij}
  \DEFECT(U)\ \  =\ \
  (N-1)^2
  \ \  -\ \
  \dim\left(
    \rule{0cm}{0.7cm}
    \SPANC\left(
      \rule{0cm}{0.5cm}
      \left\{
        \rule{0cm}{0.4cm}
        \VECC\left(
          \rule{0cm}{0.3cm}
          B^{(i,j)}
        \right)
        :\ i,j \in \{1..N\}
      \right\}
    \right)
  \right),
\end{equation}
with $B^{(i,j)}$ described by (\ref{eq_Bij_ith_row_and_column}).

%
%

\subsection{The unistochasticity problem}
\label{subsec_unistochasticity_problem}

Related to some applications in physics is the unistochasticity
problem, that is the problem of extracting full information about a
unitary matrix from the moduli of its entries only.


\begin{definition}
  \label{def_orthostochastic}

  An $N \times N$ \BISTOCHASTIC\ matrix $B$ is called
  \DEFINED{unistochastic}  (\DEFINED{orthostochastic}) if there exists
  an $N \times N$ unitary (real orthogonal)
  matrix $U$ such that\ \ \  $\forall i,j\in \{1..N\} \ \ B_{i,j} = |U_{i,j}|^2$.
\end{definition}

In other words, $B$ is unistochastic if $\VEC(B) = f(\VECR(U))$ for
some unitary $U$, having $f$ defined by (\ref{eq_f_def}). In physical
applications, $i,j$-th entries of $B$ correspond to probabilities of
obtaining the $i$-th possible result of an experiment, being one of
some chosen $N$ 'orthogonal' states of a measured quantum system, given
the $j$-th initial state was prepared. In this framework $U$, a
unitary preimage of $B$, describes possible evolution of the state of
the measured system between the moments of preparation and
measurement of the state.

A more detailed question concerning the unistochasticity issue
is the following: does there exist a unistochastic ball around the
flat matrix $J_N$, $\ELEMENTof{J_N}{i}{j} = 1/N$, within the Birkhoff's
polytope, the set of all \BISTOCHASTIC\ matrices? Note that $J_N$
is unistochastic for every $N$, since the Fourier matrix $F_N$ of
(\ref{eq_Fourier_matrix}) is its unitary preimage. A partial answer to
the posed question, which uses the notion of the defect, is provided by the lemma:


\begin{theorem}
  \label{theor_unistochastic_ball}

  Let $U$ be a unitary complex  Hadamard matrix, i.e
  $|U_{i,j}|=1/\sqrt{N}$, thus being a unitary preimage of the flat \BISTOCHASTIC\
  matrix,
  $\VECR(U) \in f^{-1}\left( \rule{0cm}{0.3cm} \VEC(J_N) \right)$,
  where $\ELEMENTof{J_N}{i}{j} = 1/N$.

  If $\DEFECT(U) = 0$ then there exists a unistochastic ball around
  $J_N$ in the set of all \BISTOCHASTIC\ matrices.
\end{theorem}

\PROOFstart   
Consider the maps (differential symbols also denote their matrix
representation
$\ELEMENTof{\DIFFERENTIAL{p}{v}}{i}{j} = \frac{\delta p_i}{\delta x_j}(v)$\ ):
\begin{itemize}
  \item
        $u:\ \SET{W} \subset \REALS^{N^2} \longrightarrow
        \REALS^{2N^2}$ parametrizing the unitary manifold
        $\VECR(\UNITARY)$ around $\VECR(U)$,
        \ \  $\SET{W}$ open ,
        such that $u(\ZEROvect \in \SET{W}) = \VECR(U)$.
        Then $\DIFFERENTIAL{u}{\ZEROvect}(\REALS^{N^2}) =
        \TANGENTspace{\UNITARY}{U}$ and columns of
        $\DIFFERENTIAL{u}{\ZEROvect}$ form a basis of
        $\TANGENTspace{\UNITARY}{U}$
  \item
        $f:\ \REALS^{2N^2} \longrightarrow \REALS^{N^2}$ defined
        earlier in (\ref{eq_f_def}),
        $\ELEMENTof{
          \VEC^{-1}\left( \rule{0cm}{0.4cm} f\left( \rule{0cm}{0.3cm} \VECR(U) \right) \right)
        }{i}{j}
        = |U_{i,j}|^2$, a
        map squaring the moduli of the entries of a complex matrix
  \item
        $m:\ \REALS^{N^2} \longrightarrow \REALS^{(N-1)^2}$,\ \
        $m\left( \rule{0cm}{0.3cm} \VEC(B) \right) =
        \VEC\left( \rule{0cm}{0.3cm} B_{1:N-1,1:N-1} \right)$,
        where the second $\VEC$
        is over $(N-1) \times (N-1)$ matrices.\ \  $m$ provides a one
        to one map between \BISTOCHASTIC\ $N \times N$ matrices and
        $(N-1) \times (N-1)$ matrices with non--negative entries,
        $m\left( \rule{0cm}{0.3cm} \VEC(J_N) \right) =
        \VEC\left( \rule{0cm}{0.4cm} \frac{N-1}{N} J_{N-1} \right)$,
        $\DIFFERENTIAL{m}{\VEC(J_N)}
        \left( \TANGENTspace{\BISTOCHSPACE}{\VEC(J_N)} \right)
 = \REALS^{(N-1)^2}$
\end{itemize}

Since $\DEFECT(U)=0$, from the definition of the defect
$\DIFFERENTIAL{f}{\VECR(U)}\left( \rule{0cm}{0.3cm} \TANGENTspace{\UNITARY}{U} \right) =
\TANGENTspace{\BISTOCHSPACE}{\VEC(J_N)}$, that is an $N^2 \times N^2$
matrix $\DIFFERENTIAL{f}{\VECR(U)} \cdot \DIFFERENTIAL{u}{\ZEROvect}$
contains a basis for $\TANGENTspace{\BISTOCHSPACE}{\VEC(J_N)}$, say at
$j_1,j_2,\ldots,j_{(N-1)^2}$-th column positions.

Consider the map $\tilde{u}:\ \REALS^{(N-1)^2} \longrightarrow
\REALS^{2N^2}$, being map $u$ restricted to its $j_1$-th, ...,
$j_{(N-1)^2}$-th variables, the other variables being set to $0$. Then
of course columns of $\DIFFERENTIAL{f}{\VECR(U)} \cdot
\DIFFERENTIAL{\tilde{u}}{\ZEROvect}$ form the above basis of
$\TANGENTspace{\BISTOCHSPACE}{\VEC(J_N)}$, and columns of the
$(N-1) \times (N-1)$ matrix
$\DIFFERENTIAL{m}{\VEC(J_N)} \cdot
\DIFFERENTIAL{f}{\VECR(U)} \cdot
\DIFFERENTIAL{\tilde{u}}{\ZEROvect}$
must form a basis for $\REALS^{(N-1)^2}$, i.e. this matrix is
non--singular.

Thus the differentiable map
$m(f(\tilde{u})):\ \REALS^{(N-1)^2} \longrightarrow \REALS^{(N-1)^2}$
satisfies the Inverse Function Theorem.
This means that each point $v$ in an open set
$\SET{V} \subset \REALS^{(N-1)^2}$ around $\VEC(\frac{N-1}{N}
J_{N-1})$ has its preimage $\tilde{w}$ in an open set $\tilde{\SET{W}}
\subset \REALS^{(N-1)^2}$ containing $\ZEROvect$.

Consider a 'pseudo--inverse' of $m$: $n$, such that:
\begin{equation}
  \label{eq_n_def}
  n\left( v \in \REALS^{(N-1)^2} \right) =
  \VEC\left(
    \left[
      \begin{array}{c|c}
        \VEC^{-1}(v)  &  \vdots \\
        \hline 
        \cdots        &  \cdot
      \end{array}
    \right]
    \in  \BISTOCHSPACE
  \right),
\end{equation}
where the second $\VEC$ is over $(N-1) \times (N-1)$ matrices, and the
$N$-th row and $N$-th column of the matrix in brackets on the right hand side
is completed to form a \BISTOCHASTIC\ matrix.
Note that $n(m) = \IDENTITY$ over $\VEC(\BISTOCHSPACE)$, and of course
$m(n) = \IDENTITY$.

The stated above property of $m(f(\tilde{u}))$ can now be rephrased
as:

Any point $n(v)$ in the set $n(\SET{V})$ open in
$\VEC(\BISTOCHSPACE)$ around $\VEC(J_N)$,
corresponding to a \BISTOCHASTIC\ matrix $\VEC^{-1}(n(v))$,
 has its preimage
$\tilde{u}(\tilde{w}) \in \tilde{u}(\tilde{\SET{W}}) \subset \VECR(\UNITARY)$,
corresponding to a unitary
matrix $\VECR^{-1}(\tilde{u}(\tilde{w}))$:
\begin{eqnarray}
  v     &  =             &  m(f(\tilde{u}(\tilde{w}))) \nonumber \\
        &  \Updownarrow  &          \label{eq_unitary_preimages} \\
  n(v)  &  =             &  f(\tilde{u}(\tilde{w})).    \nonumber
\end{eqnarray}
\PROOFend    

In section \ref{sec_Fourier_matrix_defect} we show that the defect of
the Fourier matrix $F_N$ is equal to zero only for $N$
prime. Therefore it is tempting to suppose that a unistochastic ball
around $J_N$ may not exist for composite $N$. This is indeed true for
$N=4$, as there exists a ray, stemming from $J_4$, of \BISTOCHASTIC\
matrices with the property that they have not unitary preimages \cite{BEKTZ05}:
\begin{equation}
  \label{eq_Ingemar_ray}
  \left\{\ \ \
    J_4 \ +\
    t \cdot
    \left[
      \begin{array}{rrrr}
        \rule{0cm}{0.4cm} \frac{9}{4}  &  -\frac{3}{4}  &  -\frac{3}{4}  &  -\frac{3}{4}  \\
        \rule{0cm}{0.4cm} -\frac{3}{4}  &   \frac{1}{4}  &   \frac{1}{4}  &   \frac{1}{4}  \\
        \rule{0cm}{0.4cm} -\frac{3}{4}  &   \frac{1}{4}  &   \frac{1}{4}  &   \frac{1}{4}  \\
        \rule{0cm}{0.4cm} -\frac{3}{4}  &   \frac{1}{4}  &   \frac{1}{4}  &   \frac{1}{4}
      \end{array}
    \right]\ \ \ :\ \ \
    t \in \left\langle \ -\frac{1}{9},\ 0\ \right)
  \ \ \
  \right\}.
\end{equation}
However, for $N=6$ there exists a unistochastic ball around $J_6$.
This is because the so called 'spectral matrix' $S_6$\ :
\begin{equation}
  \label{eq_S6_matrix}
  S_6\ \ \ =\ \ \
  \frac{1}{\sqrt{6}} \cdot
  \left[
    \begin{array}{cccccc}
      1 & 1 & 1 &  1 & 1 &  1\\
      1 & 1 & \omega &  \omega & \omega^2 &  \omega^2\\
      1 & \omega & 1 &  \omega^2 & \omega^2 &  \omega\\
      1 & \omega & \omega^2 &  1 & \omega &  \omega^2\\
      1 & \omega^2 & \omega^2 & \omega  & 1 &  \omega\\
      1 & \omega^2 & \omega &  \omega^2 & \omega &  1
    \end{array}
  \right],
  \ \ \ \ \ \mbox{where}\ \ \
  \omega\  =\ \exp\left( \Ii \cdot \frac{2\pi}{3} \right),
\end{equation}
found independently by Tao \cite{Tao04} and by Moorhouse
\cite{Moor01} (denoted by $S_6^{(0)}$ in our catalogue \cite{TZ06}), 
 has the defect equal to zero, so  
Theorem \ref{theor_unistochastic_ball} can be applied.
Similar examples for $N=9$ and $N=10$ can be found in \cite{BN06}, see
matrices $H_{9}$ and $BN_{10}$ there (in \cite{TZ06} they are denoted
by $N_9^{(0)}$ and $N_{10}^{(0)}$). Thus we also have unistochastic
balls around $J_9$ and $J_{10}$.

%
%

\section{The defect of a Fourier matrix}
\label{sec_Fourier_matrix_defect}

In this section we will use system (\ref{eq_U_Uconj_Rdiff_system}) to
obtain the value of defect of the $N \times N$ unitary Fourier matrix
$F_N$,
\begin{equation}
  \label{eq_Fourier_matrix}
  \ELEMENTof{F_N}{i}{j} =
  \frac{1}{\sqrt{N}} \PHASE{\frac{2\pi}{N} (i-1)(j-1)}
  \ \ \ \ \
  i,j \in \{1,2,\ldots,N\}.
\end{equation}
This value, as well as the defect of any unitary complex  Hadamard
matrix $H$ ($|H_{i,j}| = 1/\sqrt{N}$), is interesting from the point
of view of the unistochasticity issue, discussed in the previous
subsection.
If the defect of such a matrix equals zero then there exists
a unistochastic ball around the flat matrix $J_N$
($\ELEMENTof{J_N}{i}{j} = 1/N$) in the set of all \BISTOCHASTIC\
matrices, see Theorem \ref{theor_unistochastic_ball}.
Furthermore,  the result $d(F_N)=0$ implies that the Fourier matrix
is isolated, see Theorem \ref{theor_zero_defect_isolation}.
On the other hand, any positive defect $d$ of a given unitary complex  Hadamard matrix
$H$ (of which $F_N$ is an example)
gives the upper bound for the dimension of a smooth orbit of
complex Hadamard matrices, dephased with respect to H,  stemming from
$H$. This is stated by  Theorem \ref{theor_manifold_dim_defect_bound}.

A similar approach, to the one presented below in calculation of the
defect, can be used to calculate the defect of any Kronecker product
of unitary Fourier matrices. In fact, one needs to calculate the defect
only for representatives of permutation equivalence classes of such
products, see \cite{Ta06}. For instance, $F_6$ is permutation
equivalent to $F_2 \otimes F_3$, so Lemma \ref{lem_defect_permutation_invariant}
implies that their defects are equal. On the other hand,
$F_4 \otimes F_2 \otimes F_2$ and $F_4 \otimes F_4$ are permutation
inequivalent, even if we pre--multiply both products by unitary
diagonal matrices \cite{Ta06}. Thus
Lemmas \ref{lem_defect_permutation_invariant} and
       \ref{lem_defect_phasing_invariant}
cannot be used and these defects need not to be equal.

%
%

\subsection{Statement of the main result}
\label{subsec_Fourier_defect_main_result}

Before we prove a formula for the defect of the Fourier matrix of size
$N$, we need the definition of a {\sl parameter cycle matrix}, in 
which the notion of least common multiple ({\bf lcm}) is used.


\begin{definition}
  \label{def_parameter_cycle_matrix}

  A \DEFINED{parameter cycle matrix} (\DEFINED{PCM}) of size $N$ is
  any complex $N \times N$ matrix $P$ built using these rules (where
  $\PCM{P}{x}{y}$ designates parameters in matrix $P$  in a way
  different from ordinary indexing of rows and columns; we call $x$ the step
  index and $y$ the cycle index):
  \begin{itemize}
    \item
        The first column of $P$ is filled with $N$ arbitrary \emph{real}
        numbers,
        $P^{0,0},\dots,P^{0,N-1}$,
        running  from the top to the bottom.
    \item
        For the step index $j \in \{2,\ldots,\frac{N+1}{2}\}$ if $N$ is odd,
        or for $j \in \{2,\ldots,\frac{N}{2}\}$ if $N$ is even,
        the $j$-th and $(N-j+2)$-th column of $P$ are filled in such a
        way that
        \begin{eqnarray}
          \label{eq_P_jth_column_pattern}
          \lefteqn{\PCM{P}{j-1}{k-1} =}  &  &  \\
          &
          P_{k,\ j}                   \ \  =\ \
          P_{(k+(j-1))\!\!\!\!\!\mod N,\ j}     \ \  =\ \
          \ldots                      \ \  =\ \
          P_{\left(
               k +
               \left( \frac{\lcm(N,j-1)}{j-1}-1\right ) (j-1)
             \right)\!\!\!\!\!
             \mod N,\ j}\ \ ,    &                          \nonumber \\
          \lefteqn{\CONJ{\PCM{P}{j-1}{k-1}} =}  &  &          \nonumber \\
          &
          P_{k,\ N-j+2}                   \ \ =\ \
          P_{(k+(j-1))\!\!\!\!\!\mod N,\ N-j+2}     \ \ =\ \
          \ldots                      \ \  =\ \
          P_{\left(k +
               \left( \frac{\lcm(N,j-1)}{j-1}-1 \right) (j-1)\!\!
             \right)\!\!\!\!\!
             \mod N,\ N-j+2}\ \ ,  &                        \nonumber
        \end{eqnarray}
        where
        \begin{equation}
          \label{eq_cycle_indices}
          k = 1,\ 2,\ \ldots,\ \frac{N}{\frac{\lcm(N,j-1)}{j-1}}
        \end{equation}
        designates the $(k-1)$-th\ \
        $\left( \rule{0cm}{0.3cm} \lcm(N,j-1) \right)/(j-1)$ element 'cycle',
        and $\PCM{P}{j-1}{k-1}$ are arbitrary \emph{complex}
        parameters of $P$.
    \item
        If $N$ is even, then the $(N/2 + 1)$-th column is filled
        according to the pattern:
        \begin{equation}
          \label{eq_P_central_column_pattern}
          \PCM{P}{\frac{N}{2}}{k-1}      \ \ \  = \ \ \
          P_{k,\frac{N}{2}+1}              \ \  = \ \
          P_{k+\frac{N}{2},\frac{N}{2}+1}
          \ \ \ \ \ \mbox{for}\ \ \  k = 1,\ 2,\ \ldots,\ \frac{N}{2}
          \ \ ,\ \ \ \ \ \ \ \ \ \ \ \ \ \ \ \
        \end{equation}
        where $\PCM{P}{\frac{N}{2}}{k-1}$ are arbitrary \emph{real} parameters.
  \end{itemize}
\end{definition}

As an example we provide a parameter cycle matrix of order $6$, 
\begin{equation}
  \label{eq_P_6x6_matrix}
  \left[
    \begin{array}{c|c|c|c|c|c}
      \rule{0cm}{0.5cm}
      \PCM{P}{0}{0} &
      \PCM{P}{1}{0} & \PCM{P}{2}{0} &
      \PCM{P}{3}{0} &
      \CONJ{\PCM{P}{2}{0}} & \CONJ{\PCM{P}{1}{0}} \\
      \hline 
      \rule{0cm}{0.5cm}
      \PCM{P}{0}{1} &
      \PCM{P}{1}{0} & \PCM{P}{2}{1} &
      \PCM{P}{3}{1} &
      \CONJ{\PCM{P}{2}{1}} & \CONJ{\PCM{P}{1}{0}} \\
      \hline 
      \rule{0cm}{0.5cm}
      \PCM{P}{0}{2} &
      \PCM{P}{1}{0} & \PCM{P}{2}{0} &
      \PCM{P}{3}{2} &
      \CONJ{\PCM{P}{2}{0}} & \CONJ{\PCM{P}{1}{0}} \\
      \hline 
      \rule{0cm}{0.5cm}
      \PCM{P}{0}{3} &
      \PCM{P}{1}{0} & \PCM{P}{2}{1} &
      \PCM{P}{3}{0} &
      \CONJ{\PCM{P}{2}{1}} & \CONJ{\PCM{P}{1}{0}} \\
      \hline 
      \rule{0cm}{0.5cm}
      \PCM{P}{0}{4} &
      \PCM{P}{1}{0} & \PCM{P}{2}{0} &
      \PCM{P}{3}{1} &
      \CONJ{\PCM{P}{2}{0}} & \CONJ{\PCM{P}{1}{0}} \\
      \hline 
      \rule{0cm}{0.5cm}
      \PCM{P}{0}{5} &
      \PCM{P}{1}{0} & \PCM{P}{2}{1} &
      \PCM{P}{3}{2} &
      \CONJ{\PCM{P}{2}{1}} & \CONJ{\PCM{P}{1}{0}} \\
    \end{array}
  \right].
\end{equation}

%
%
%

The notion of the parameter cycle matrices allows us to obtain 
concrete results on the defect of the Fourier matrix $F_N$ of size $N$.
It can be expressed by a sum of greatest common divisors (gcd).


\begin{theorem}
  \label{theor_Fourier_matrix_defect}

  For $N$ being a natural number
  \begin{equation}
    \label{eq_Fourier_matrix_defect}
    \DEFECT(F_N) =
    \left\{
      \begin{array}{cc}
        \rule[-0.3cm]{0cm}{0.8cm}
        1 - N + 2\sum_{l=1}^{\frac{N-1}{2}} \gcd(N,l) &
        \mbox{ for $N$ odd},  \\
        \rule[-0.3cm]{0cm}{0.8cm}
        1 -\frac{N}{2} + 2\sum_{l=1}^{\frac{N}{2}-1} \gcd(N,l) &
        \mbox{ for $N$ even}. \nonumber \\
      \end{array}
    \right.
  \end{equation}
\end{theorem}

\PROOFstart 
We rewrite system (\ref{eq_U_Uconj_Rdiff_system}) for $F_N$, denoted
further as $F$:
\begin{equation}
  \label{eq_F_Fconj_Rdiff_system}
  \forall 1 \leq i < j \leq N\ \ \ \
  \sum_{k=1}^{N} F_{i,k} \CONJ{F}_{j,k} ( R_{i,k} - R_{j,k} )
  = 0
\end{equation}
as
\begin{equation}
  \label{eq_RFtransposed_property_1}
  \sum_{k=1}^{N} R_{i,k} F_{j-i+1,k}  =
  \sum_{k=1}^{N} R_{j,k} F_{j-i+1,k}\ \ ,
\end{equation}
and there also generally holds that
\begin{equation}
  \label{eq_RFtransposed_property_2}
  \sum_{k=1}^{N} R_{i,k} F_{j-i+1,k}  =
  \sum_{k=1}^{N} R_{i,k} \CONJ{F}_{N-(j-i)+1,k}  =
  \CONJ{ \sum_{k=1}^{N} R_{i,k} F_{N-(j-i)+1,k} }\ \ .
\end{equation}

We introduce a complex $N \times N$ matrix $P$ such that:
\begin{equation}
  \label{eq_P_definition}
  \CONJ{P}\ \stackrel{\mbox{def}}{=}\  R \cdot F^T = R \cdot F
  \ \ \Longleftrightarrow\ \
  P = RF^*
  \ \ \Longleftrightarrow\ \
  R = PF
\end{equation}
Then statements (\ref{eq_RFtransposed_property_1}) and
(\ref{eq_RFtransposed_property_2}) can be expressed in terms the
elements of matrix $P$ as:
\begin{equation}
  \label{eq_P_matrix_properties}
  \begin{array}{ccc}
    \CONJ{P}_{i,j-i+1} & = & P_{i,N-(j-i)+1}  \\
           \|          &   &       \|         \\
    \CONJ{P}_{j,j-i+1} & = & P_{j,N-(j-i)+1}
  \end{array}
\end{equation}

Rules (\ref{eq_P_matrix_properties}) as well as the requirement of
matrix $R$ being real force matrix $P$ to be a parameter cycle
matrix of Definition \ref{def_parameter_cycle_matrix}, and the
solution space of (\ref{eq_F_Fconj_Rdiff_system}) is
fully parametrized by the formula $R = PF$, where $P$ is any such PCM matrix.

The total number of \emph{real} parameters in $P$, parametrizing the solution
space of (\ref{eq_F_Fconj_Rdiff_system}), reduced by $(2N-1)$ to
become the defect of $F_N$, reads:
\begin{itemize}
  \item
        $N$ odd:
        \begin{eqnarray}
          \label{eq_F_N_odd_defect}
          \DEFECT(F_N) \ \ =\ \
          N\  +\
          2 \left(
              \sum_{l=1}^{\frac{N-1}{2}}
                \frac{N}{\frac{\lcm(N,l)}{l}}
            \right)
          \ -\  (2N-1) \ \ = \nonumber
          & & \\ 
          2 \sum_{l=1}^{\frac{N-1}{2}}
                    \left( \frac{N}{\frac{\lcm(N,l)}{l}} - 1 \right)
          \ \ =\ \
          2 \sum_{l=1}^{\frac{N-1}{2}}
                    \left( \rule{0cm}{0.4cm} \gcd(N,l) - 1 \right), & &  \nonumber
        \end{eqnarray}
\item
        $N$ even:
        \begin{eqnarray}
          \label{eq_F_N_even_defect}
          \DEFECT(F_N) \ \ =\ \
          N \ +\
          2 \left(
              \sum_{l=1}^{\frac{N}{2}-1}
                \frac{N}{\frac{\lcm(N,l)}{l}}
            \right)
          \ +\  \frac{N}{2} \ -\  (2N-1) =  \nonumber
          & & \\  
          2\sum_{l=1}^{\frac{N}{2}-1}
                    \left( \frac{N}{\frac{\lcm(N,l)}{l}} - \frac{1}{2}
                    \right)
          \ \ =\ \
          2\sum_{l=1}^{\frac{N}{2}-1}
                    \left( \gcd(N,l) - \frac{1}{2} \right). & & \nonumber
        \end{eqnarray}
\end{itemize}

That is:
\begin{eqnarray}
  \label{eq_F_N_defect}
  \mbox{$N$ odd:} &
  \DEFECT(F_N)\ \  =\ \   1  -  N  +  2\sum_{l=1}^{\frac{N-1}{2}} \gcd(N,l)  &  \\
  \rule{0cm}{0.8cm}
 \label{eq_F_N_defect1}
  \mbox{$N$ even:} &
  \DEFECT(F_N)\ \  =\ \   1  -  \frac{N}{2}  +  
                2\sum_{l=1}^{\frac{N}{2}-1} \gcd(N,l) &  
\end{eqnarray}
\PROOFend 

Alternative  formulas,
for the defect of $F_N$ can be useful.


\begin{theorem}
  \label{theor_WS_Fourier_matrix_defect_formula}

  For any natural $N \geq 2$ with the factorization into prime numbers:
  \begin{equation}
    \label{eq_N_factorization}
    N = \prod_{j=1}^{n} {p_j}^{k_j}
    \ \ \ \ \
    p_1 > p_2 > \ldots > p_n
  \end{equation}
  there holds
  \begin{description}
    \item{a)}
          \begin{equation}
            \label{eq_WS_aux_formula}
            \DEFECT(F_N) =
            \sum_{l=1}^{N-1} \left( \gcd(N,l) - 1 \right)
            \text{ ,}
          \end{equation}
    \item{b)}
          \begin{equation}
            \label{eq_WS_main_formula}
            \DEFECT(F_N) =
            N \cdot
            \left(
              \prod_{j=1}^{n} (1 + k_j - \frac{k_j}{p_j})   -   2
            \right)  +  1
            \text{ .}
          \end{equation}
   \end{description}
\end{theorem}

The proof is provided in Appendix \ref{sec_proof}.

%
%

\subsection{Some special cases}
\label{subsec_Fourier_defect_special_cases}

Since the explicit formula
(\ref{eq_WS_main_formula}) is not very transparent
the defects of Fourier matrices for small dimensionalities
are collected in table 1.
\bigskip
\begin{table} [ht]
\caption{Defect for
the Fourier matrix $F_N$ of size $N$.}
\smallskip
{\renewcommand{\arraystretch}{1.45}
\begin{tabular}
[c]{||c||c|c|c|c|c|c|c|c|c|c|c|c|c|c|c|c||}%
\hline \hline  $N$   &
  $1$ & $2$ & $3$ & $4$ & $5$ &$6$ &$7$ &$8$ &$9$ &$10$ &$11$ &$12$ &$13$ &$14$ &$15$ &$16$ \\
\hline $\DEFECT(F_N)$ &
  $0$ & $0$ & $0$ & $1$ & $0$ & $4$ & $0$ &$5$ & $4$ &$8$ &$0$ &$17$ &$0$ &$12$ & $16$ & $17$ \\
\hline \hline  $N$   &
 $17$ &$18$ &$19$ &$20$ &$21$ &$22$ &$23$ &$24$ &$25$ &$26$ &$27$ &$28$ &$29$ &$30$ &$31$ &$32$ \\
\hline $\DEFECT(F_N)$ &
 $0$ & $28$ & $0$ & $33$ & $24$ & $20$ & $0$ &$53$ & $16$ &$24$ &$28$ &$49$ &$0$ &$76$ & $0$ & $49$ \\
\hline \hline
\end{tabular} }
\label{tab1}  \end{table}

Let us now discuss some special cases
of the formula
for the defect of the Fourier matrix $F_N$.
\smallskip

\noindent
{\bf i).}  $N$ is prime.

If $N=p$ then $n=1$, $k_1=1$, so
the right hand side of equation (\ref{eq_WS_main_formula})
reads $p(2-1/p-2)+1$ and provides the result
$\DEFECT(F_{p}) = 0$ as advertised.
Hence the Fourier matrix of a prime dimension is isolated.
\medskip

\noindent
{\bf ii).} $N$ is a product of two distinct primes.

If $N=pq$ then $n=2$, $k_1=k_2=1$
so (\ref{eq_WS_main_formula}) reads\\
$pq \left[ \rule{0cm}{0.3cm} (2-1/p)(2-1/q)-2 \right]  +  1$, which gives:
\begin{equation}
  \label{eq_F_product_of_2primes_defect}
            \DEFECT(F_{pq}) \ =\ \  2(p-1)(q-1)
\end{equation}
It is worth to emphasize that the upper bound $\DMPbound(F_{pq})$
(see equation (\ref{eq_manifold_dim_defect_bound})) 
for the dimension of an orbit of dephased
unitary complex Hadamard matrices stemming from $F_{pq}$
implied by this formula is exactly {\sl twice}
the dimension $D$ of the orbits actually constructed  in  \cite{Ha96,Di04,TZ06}
for a product of primes, $D=(p-1)(q-1)$. The problem of
describing the entire (possibly existing) manifold of dephased unitary complex Hadamard
matrices stemming from the Fourier matrix $F_{pq}$
is open even in the simplest case $N=2\cdot 3=6$ \cite{TZ06},
but a recent discovery of a new 'non--affine' (according to the proper
definition in our catalogue \cite{TZ06}) $N=6$ orbit
of unitary complex Hadamard matrices \cite{BN06},
 and some further results
seem to suggest that in this case a full $4$ dimensional orbit does exist \cite{Bxxx07}.
\medskip

\noindent
 {\bf iii).} $N$ is a product of three distinct primes.

If $N=pqr$ then $n=3$, $k_1=k_2=k_3=1$
and eq. (\ref{eq_WS_main_formula}) amounts to:
\begin{equation}
  \label{eq_F_product_of_3primes_defect}
  \DEFECT(F_{pqr}) \ = \
  2 \left[ \rule{0cm}{0.3cm} 3pqr\ -\ 2(pq+pr+qr)\ +\ (p+q+r) \right] \ .
\end{equation}
\medskip

\noindent
{\bf iv).} $N$ is a power of two.
If $N=2^k$ then $p=2,\ n=1$ and  $k_1=k$,
so  (\ref{eq_WS_main_formula}) leads to:
\begin{equation}
   \label{eq_F_power_of_2_defect}
   \DEFECT(F_{2^k}) \ = \ 2^{k-1}(k-2)+1  \ .
\end{equation}
\medskip

\noindent
{\bf v).} $N$ is a power of a prime.

If $N=p^k$ then $n=1$ and  $k_1=k$,
so  (\ref{eq_WS_main_formula}) takes the form of:
\begin{equation}
  \label{eq_F_power_of_prime_defect}
  \DEFECT(F_{p^k}) \ = \
  p^{k-1} \left[ \rule{0cm}{0.3cm} (p-1)k\ -\ p \right]\ \ +\ \ 1  \ .
\end{equation}
Interestingly, in this very case the defect is equal to
the dimension of the known
smooth orbits of dephased unitary complex Hadamard matrices stemming from $F_{p^k}$,
featured in Section \ref{sec_defect_dimensional_Fourier_orbits}.
This shows that these solutions are complete in the sense that they
are not contained in smooth orbits (of the
respective type) of a higher dimension.

%
%

\section{Orbits of the maximal dimension stemming from Fourier
         matrices of a prime power size}
\label{sec_defect_dimensional_Fourier_orbits}

In this section we present examples of $N \times N$ unitary matrices $U$ with no
zero entries, for which there exist $\DEFECT(U)$--dimensional smooth
families (manifolds) generated, through $V \rightarrow \VECR(V)$, by unitary
matrices $V$ with the same  pattern of moduli as in $U$, and dephased with
respect to $U$.

As $U$ has no zero entries, a spanning set $\SET{S}$ for $U$ (see
Definition \ref{def_dephased_matrix}) will always have $p+r = 2N-1$
independent vectors as its elements. So, according to
Theorem \ref{theor_manifold_dim_defect_bound}, a manifold 
of the type described above will have its dimension bounded just by
$\DEFECT(U)$.

We will consider Fourier matrices $F_{p^k}$, of the size being the $k$-th
natural power of a prime number $p$, as examples for which this bound
is saturated. To make the notion of being dephased with respect to
$F_{p^k}$ precise, as Definition \ref{def_dephased_matrix} requires,
and also for practical reasons, we arbitrarily choose the spannig set
for $F_{p^k}$ to be:
\begin{equation}
  \label{eq_FpTOk_spanning_set}
  \SET{S}_{F_{p^k}} =
  \left\{
    \rule{0cm}{0.4cm}
    \VECR\left(
      \rule{0cm}{0.3cm}
      \Ii \cdot \diag(\STbasis{r}) \cdot F_{p^k}
    \right)
    :\ r \in \{1..p^k\}
  \right\}
  \
  \cup
  \
  \left\{
    \rule{0cm}{0.4cm}
    \VECR\left(
      \rule{0cm}{0.3cm}
      F_{p^k} \cdot \Ii \cdot \diag(\STbasis{c})
    \right)
    :\ c \in \{2..p^k\}
  \right\},
\end{equation}
and the pattern set to be:
\begin{equation}
  \label{eq_FpTOk_pattern_set}
  \SET{I}_{F_{p^k}} =
  \left\{
    (1,1),\ (2,1),\ \ldots,\ (p^k,1)
  \right\}
  \
  \cup
  \
  \left\{
    (1,2),\ (1,3),\ \ldots,\ (1,p^k)
  \right\}.
\end{equation}
In other words, $V$ is dephased with respect to $F_{p^k}$, according to
$\SET{I}_{F_{p^k}}$ associated with $\SET{S}_{F_{p^k}}$, if the
entries in the first row and column of $V$ are equal
to the corresponding entries in $F_{p^k}$, 
i.e. they are all equal to $1/\sqrt{p^k}$. 
Of course
$V$ is assumed to be a unitary complex Hadamard matrix, that is $VV^*
= I$,\ \ $|V_{i,j}| = 1/\sqrt{p^k}$\ \  for\ \ $i,j=1..p^k$.

To construct a $\DEFECT(F_{p^k})$-dimensional manifold, generated by
dephased unitary complex Hadamard matrices, stemming from $\VECR(F_{p^k})$, we
have to take a subspace of the space of all parameter cycle matrices
(PCM matrices) $P$ of size $p^k$, introduced in
Definition \ref{def_parameter_cycle_matrix}.
Because of the dephasing condition, and this will be made clear in the
proof of the theorem below, we have to impose on $P$ additional
constraints.
We have to set all $p^k$ real parameters of $P$ sitting in the
first row to zero, and each of the remaining $(p^k-1)$ (out of the
total of $2p^k-1$ to be fixed) real parameters
sitting in the first column to minus the sum of the remaining complex
parameters sitting in the same row as the  parameter (in the $1$-st column)
being set. This leaves $\DEFECT(F_{p^k})$ real parameters free. Then
the second column of $P$ as well as each $j$-th column with $(j-1)$ not
divided by $p$ are filled all with zeros.
Using these constraints as well as the alternative indexing of
parameters in $P$ 
(see Definition \ref{def_parameter_cycle_matrix}), we state that:


\begin{theorem}
  \label{theor_Fp_family_PCM_construction}

  The $\DEFECT(F_{p^k})$-parameter family
  \begin{equation}
    \label{eq_FpTOk_family_PCM}
    \left\{
      \rule{0cm}{1.5cm}
      \VECR\left(
        \rule{0cm}{0.3cm}
        F_{p^k} \HADprod \EXPentrywise(\Ii P F_{p^k})
      \right)
      :\
      P\ \mbox{is PCM, and}
      \left(
          \begin{array}{c}
              \PCM{P}{j}{0} = 0 \ \mbox{for}\
              \left\{
                  \begin{array}{c}
                    j \in \left\{ \rule{0cm}{0.4cm} 1..(\frac{p^k-1}{2})\right\}\ \mbox{if}\ \ p \neq 2
                    \\
                    j \in \left\{ \rule{0cm}{0.4cm} 1..(\frac{p^k}{2}) \right\}\ \mbox{if}\ \ p = 2
                  \end{array}
              \right.
              \\
              \PCM{P}{0}{i} = \displaystyle{-\sum_{j=2}^{p^k}
                P_{i+1,j}}\ \ \mbox{for}\ \
              i \in \left\{ \rule{0cm}{0.3cm} 0..(p^k-1)\right\}
          \end{array}
          \!\!\!\!\!
      \right)
    \right\}
  \end{equation}
  is a differentiable manifold stemming from $\VECR(F_{p^k})$, and
  is generated, through $V \rightarrow \VECR(V)$, by unitary complex Hadamard matrices $V$
  dephased with respect to $F_{p^k}$ according to $\SET{I}_{F_{p^k}}$
  associated with $\SET{S}_{F_{p^k}}$.
\end{theorem}

{
\newcommand{\PreBASIS}[2]{P^{(#1,#2)}_{\Re}}
\newcommand{\PimBASIS}[2]{P^{(#1,#2)}_{\Im}}
\PROOFstart
Let $P$ satisfy the constraints formulated in
(\ref{eq_FpTOk_family_PCM}).

The first row of $PF_{p^k}$ is filled with $0$'s, as $P_{1,:} =
\ZEROvect$. In the first column the entries satisfy, due to
the constraints imposed on $P$:
\begin{equation}
  \label{eq_PF_first_column_entries}
  \ELEMENTof{PF_{p^k}}{i+1}{1} \ =\
  \frac{1}{\sqrt{p^k}} \cdot  \sum_{j=1}^{p^k} P_{i+1,j} \ =\
  \frac{1}{\sqrt{p^k}} \cdot
  \left(
    \PCM{P}{0}{i} + \sum_{j=2}^{p^k} P_{i+1,j}
  \right)
  = 0\ .
\end{equation}
Thus $PF_{p^k}$, which is real thanks to the PCM structure of $P$, has its first
row and column filled with $0$'s, so
$F_{p^k} \HADprod \EXPentrywise(\Ii PF_{p^k})$, if unitary, is indeed
dephased with respect to $F_{p^k}$.

Next we will show that $P F_{p^k} \in \SET{R}_{F_{p^k}}$, a linear subspace of
real $p^k \times p^k$ matrices, defined in the following
Theorem \ref{theor_Fp_family_direct_construction}, presenting another
construction of the considered family. From this it will follow that
$F_{p^k} \HADprod \EXPentrywise(\Ii PF_{p^k})$ is unitary, as it is shown in
the proof of Theorem \ref{theor_Fp_family_direct_construction}.
Note that $P F_{p^k}$ already satisfies the 
dephasing constraints parltly defining $\SET{R}_{F_{p^k}}$.

To have $P F_{p^k}  \in \SET{R}_{F_{p^k}}$ it is enough to show that for $p^k
\times p^k$ PCM matrices
\begin{description}
  \item[$\PreBASIS{j}{i}, \PimBASIS{j}{i}$]
        with the step index $j$ and the cycle index $i$ in the ranges:
        \begin{equation}
          \label{eq_side_P_Re_Im_ji_j_i_range}
          j \in
          \left\{
            \begin{array}{ccc}
              \left\{ \rule{0cm}{0.5cm} 1..\left( \frac{p^k-1}{2} \right) \right\} &
              \!\!\mbox{if} & \!\!\! p \neq 2 \\
              \left\{ \rule{0cm}{0.5cm} 1..\left( \frac{p^k}{2}-1 \right) \right\} &
              \!\!\mbox{if} & \!\!\! p=2
            \end{array}
          \right.
          \ \ \ \ \ \mbox{and}\ \ \ \ \
          i \in
          \left\{
            \rule{0cm}{0.7cm}
            1,2,\ldots,\left( \frac{p^k}{\frac{\lcm(p^k,j)}{j}} - 1 \right)
          \right\}\ ,
        \end{equation}
  \item[$\PreBASIS{\frac{p^k}{2}}{i}$ additionally for $p=2$]
        with the cycle index $i$ in the range:
        \begin{equation}
          \label{eq_central_P_Re_i_range}
          i \in
          \left\{
            1,2,\ldots,\left(\frac{p^k}{2}-1 \right)
          \right\}\ ,
        \end{equation}
\end{description}
and such that they satisfy the additional constraints imposed on $P$ in
Theorem \ref{theor_Fp_family_PCM_construction},
and (PCM indexing introduced in
Definition \ref{def_parameter_cycle_matrix} is used):
\begin{itemize}
  \item
        $\pcmELEMENTof{\PreBASIS{j}{i}}{j}{i} = 1$,\ \ \ \
        $\pcmELEMENTof{\PreBASIS{j}{i}}{t}{s} = 0$\ \  for other allowed
        $(t,s) \neq (j,i),\ t \geq 1$,
  \item
        $\pcmELEMENTof{\PimBASIS{j}{i}}{j}{i} = \Ii$,\ \ \ \
        $\pcmELEMENTof{\PimBASIS{j}{i}}{t}{s} = 0$\ \  for other allowed
        $(t,s) \neq (j,i),\ t \geq 1$,
  \item
        if $p=2$ then
        $\pcmELEMENTof{\PreBASIS{\frac{p^k}{2}}{i}}{\frac{p^k}{2}}{i}
        = 1$,\ \ \ \
        $\pcmELEMENTof{\PreBASIS{\frac{p^k}{2}}{i}}{t}{s} = 0$\ \  for other allowed
        $(t,s) \neq (\frac{p^k}{2},i),\ t \geq 1$,
\end{itemize}
$\left( \rule{0cm}{0.3cm} \PreBASIS{j}{i} F_{p^k} \right)$,\ \
$\left( \rule{0cm}{0.3cm} \PimBASIS{j}{i} F_{p^k} \right)$\ \
(and\ \
$\left( \rule{0cm}{0.3cm} \PreBASIS{\frac{p^k}{2}}{i}F_{p^k} \right)$\ \
if $p=2$)
belong to $\SET{R}_{F_{p^k}}$.
It is because\ \  $\PreBASIS{j}{i}$,\ \  $\PimBASIS{j}{i}$ form
a basis of the space of PCM matrices $P$ satisfying the constraints of
Theorem \ref{theor_Fp_family_PCM_construction}, and there are
$\DEFECT(F_{p^k})$ of them (see the calculation of the number of real
parameters in a PCM matrix leading to the calculation of the defect of
$F_N$ in the proof of Theorem \ref{theor_Fourier_matrix_defect}).

The\ \
$\left( \rule{0cm}{0.3cm} \PreBASIS{j}{i} F_{p^k} \right)$,\ \
$\left( \rule{0cm}{0.3cm} \PimBASIS{j}{i} F_{p^k} \right)$\ \
matrices have the properties:
\begin{itemize}
   \item
        The only non--zero rows of
        $\left( \rule{0cm}{0.3cm} \PreBASIS{j}{i} F_{p^k} \right)$ are
        identical and equal to:
        \begin{equation}
          \label{eq_side_PreF_nonzero_row}
          -2 \ELEMENTof{F_{p^k}}{1}{:} +
             \ELEMENTof{F_{p^k}}{j+1}{:} +
             \ELEMENTof{F_{p^k}}{p^k-j+1}{:}
          =
          -2 \frac{1}{\sqrt{p^k}} \ONESvect^T +
           2 \Re\left( \ELEMENTof{F_{p^k}}{j+1}{:} \right)\ ,
         \end{equation}
         and they are spaced at row index distance $\gcd(j,p^k)$ one from the
         next one below.
   \item
         The only non--zero rows of
         $\left( \rule{0cm}{0.3cm} \PimBASIS{j}{i} F_{p^k} \right)$ are
        identical and equal to:
        \begin{equation}
          \label{eq_side_PimF_nonzero_row}
          0     \ELEMENTof{F_{p^k}}{1}{:} +
          \Ii   \ELEMENTof{F_{p^k}}{j+1}{:}
         -\Ii   \ELEMENTof{F_{p^k}}{p^k-j+1}{:}
          =
          -2 \Im\left( \ELEMENTof{F_{p^k}}{j+1}{:} \right)\ ,
         \end{equation}
         and they are spaced at row index distance $\gcd(j,p^k)$ one from the
         next one below.
   \item
         If $p=2$, the only non--zero rows of
         $\left( \rule{0cm}{0.3cm} \PreBASIS{\frac{p^k}{2}}{i}F_{p^k} \right)$ are
        identical and equal to:
        \begin{equation}
          \label{eq_central_PreF_nonzero_row}
          -1 \ELEMENTof{F_{p^k}}{1}{:} +
             \ELEMENTof{F_{p^k}}{\frac{p^k}{2}+1}{:}
          =
          -\frac{1}{\sqrt{p^k}} \ONESvect^T +
           \ELEMENTof{F_{p^k}}{\frac{p^k}{2}+1}{:}
          \ \ \ \ \ \ \mbox{(a real one!)},
         \end{equation}
         and they are spaced at row index distance $p^k/2$ one from the
         next one below.
\end{itemize}

Now take any allowed step index $j$ and the corresponding matrix
$\left( \rule{0cm}{0.3cm} \PreBASIS{j}{i} F_{p^k} \right)$
or matrix 
$\left( \rule{0cm}{0.3cm} \PimBASIS{j}{i} F_{p^k} \right)$,
for any allowed cycle index $i$.
Let $\gcd(j,p^k)=p^m,\ m<k$, i.e. $j=ap^m$ with $a,p$
relatively prime.

Since any two non--zero rows of the considered matrix are spaced at a
row index distance being a multiplicity of $p^m$, the matrix
automatically satisfies
the constraints of order $p^m$, order $p^{m+1}$, ..., order
  $p^{k-1}$
defining $\SET{R}_{F_{p^k}}$ in
Theorem \ref{theor_Fp_family_direct_construction}.
For\ \  $p=2$\ \  and\ \  $j=p^k/2$\ \
$\left( \rule{0cm}{0.3cm} \PreBASIS{\frac{p^k}{2}}{i}F_{p^k} \right)$
surely satisfies  the constraints of order $p^{k-1}$.
The above holds because the \emph{respective} differences of rows in
the considered matrix are zero rows in all cases.

For the other constraints to hold, it is obviously sufficient that 
the constraints 
of order $1$, order $p$, ..., order $p^{m-1}$ 
hold for the universal rows $\DELTA{}{}$:
\begin{equation}
  \label{eq_universal_DELTA}
  \DELTA{}{} = \Re\left( \ELEMENTof{F_{p^k}}{j+1}{:} \right)
  \ \ \ \
  \mbox{or}
  \ \ \ \
  \DELTA{}{} = \Im\left( \ELEMENTof{F_{p^k}}{j+1}{:} \right)\ .
\end{equation}

The constraints of order $p^{\tilde{m}}$, with $0 \leq \tilde{m} < m$,
require that
\begin{equation}
  \label{eq_universal_DELTA_constraints}
  \DELTA{}{}_l =
  \DELTA{}{}_{l+p^{k-(\tilde{m}+1)}} =
  \DELTA{}{}_{l+2 p^{k-(\tilde{m}+1)}} =
  \ldots =
  \DELTA{}{}_{l+(p-1) p^{k-(\tilde{m}+1)}}\ \ ,
\end{equation}
and note that (\ref{eq_universal_DELTA_constraints}) is true for either definition
of $\DELTA{}{}$ in (\ref{eq_universal_DELTA}), because for the allowed
(in the definition of $\SET{R}_{F_{p^k}}$)
natural $l \geq 1$ and $s \geq 0$:
\begin{eqnarray}
  \label{eq_Fourier_row_phase_rotation}
  \lefteqn{
    \ELEMENTof{F_{p^k}}{j+1}{l+sp^{k-(\tilde{m}+1)}}
    =
    \ELEMENTof{F_{p^k}}{j+1}{l} \cdot
    \PHASE{\left( \frac{2\pi}{p^k} j \right) \cdot
      sp^{k-(\tilde{m}+1)}}
    =
  } & & \\ 
  &
  \ELEMENTof{F_{p^k}}{j+1}{l} \cdot
    \PHASE{\frac{2\pi}{p^{k-m}} a  \cdot sp^{k-(\tilde{m}+1)}}
  =
  \ELEMENTof{F_{p^k}}{j+1}{l} \cdot
       \PHASE{ 2\pi \cdot a  \cdot s \cdot p^{(m - \tilde{m} - 1)}}
  =
  \ELEMENTof{F_{p^k}}{j+1}{l}\ \ .
  & \nonumber
\end{eqnarray}
This also applies to\ \
$\left( \rule{0cm}{0.3cm} \PreBASIS{\frac{p^k}{2}}{i}F_{p^k} \right)$
\ \  if\ \  $p=2$.

From the above we conclude that the space of $P F_{p^k}$, with $P$
satisfying the constraints of
Theorem \ref{theor_Fp_family_PCM_construction},
is contained within $\SET{R}_{F_{p^k}}$.
For the final argument that (\ref{eq_FpTOk_family_PCM}) is indeed a
manifold, we refer the reader to the similar one in the ending of the
proof of Theorem \ref{theor_Fp_family_direct_construction}, that
follows. See also the remark on the possible dimension of
$\SET{R}_{F_{p^k}}$ there.

\PROOFend 

Another construction of the discussed $\DEFECT(F_{p^k})$
dimensional family stemming from $F_{p^k}$ is presented in the theorem
below. However, the way in which the free parameters (phases)  are
scattered around a member matrix of the family seems to be  more sophisticated in
comparison with the pattern of parameters in a PCM matrix of the
previous theorem.


\begin{theorem}
  \label{theor_Fp_family_direct_construction}

  Let, for $p$ prime and $k\in \NATURAL$ such that $k>1$, $\SET{R}_{F_{p^k}}$ be
  the set of all real $p^k \times p^k$ matrices $R$
  satisfying the independent constraints (where $\DELTA{i}{j}_l$
  denotes the difference $R_{i,l} - R_{j,l}$):

  \begin{description}
    \item[the constraints of order $p^0 = 1$ :]
          \begin{equation}
            \label{eq_order_1_constraints}
            \DELTA{i}{i+1}_l = \DELTA{i}{i+1}_{l+p^{k-1}} =
            \DELTA{i}{i+1}_{l+2p^{k-1}} = \ldots =
            \DELTA{i}{i+1}_{l+(p-1)p^{k-1}}
          \end{equation}
          for $l = 1,\ 2,\ \ldots,\ p^{k-1}$,\\
          and for $i \in \{1,2,\ldots,(p^{k}-1)\}$\\
          (one $p^k$--element cycle of differences of rows, counting
          $\DELTA{p^k}{1}_{...}$).

    \item[the constraints of order $p$ :]
          \begin{equation}
            \label{eq_order_p_constraints}
            \DELTA{i}{i+p}_l = \DELTA{i}{i+p}_{l+p^{k-2}} =
            \DELTA{i}{i+p}_{l+2p^{k-2}} = \ldots =
            \DELTA{i}{i+p}_{l+(p-1)p^{k-2}}
          \end{equation}
          for $l = 1,\ 2,\ \ldots,\ p^{k-2}$,\\
          and for
          $i \in \bigcup_{r \in \{1,\ldots,p\}}
                   \left\{
                       r+sp:\ s \in \{0,1,\ldots,(p^{k-1}-2)\}
                   \right\}$\\
          ($p$\ \ $p^{k-1}$--element cycles of differences, counting final
          wrappings).

    \item[the constraints of order $p^2$ :]
          \begin{equation}
            \label{eq_order_p_to_2_constraints}
            \DELTA{i}{i+p^2}_l = \DELTA{i}{i+p^2}_{l+p^{k-3}} =
            \DELTA{i}{i+p^2}_{l+2p^{k-3}} = \ldots =
            \DELTA{i}{i+p^2}_{l+(p-1)p^{k-3}}
          \end{equation}
          for $l = 1,\ 2,\ \ldots,\ p^{k-3}$,\\
          and for
          $i \in \bigcup_{r \in \{1,\ldots,p^2\}}
                   \left\{
                       r+sp^2:\ s \in \{0,1,\ldots,(p^{k-2}-2)\}
                   \right\}$\\
          ($p^2$\ \  $p^{k-2}$--element  cycles of differences).

    \item[$\ldots$]

    \item[the constraints of order $p^m$ :]
          \begin{equation}
            \label{eq_order_p_to_m_constraints}
            \DELTA{i}{i+p^m}_l = \DELTA{i}{i+p^m}_{l+p^{k-(m+1)}} =
            \DELTA{i}{i+p^m}_{l+2p^{k-(m+1)}} = \ldots =
            \DELTA{i}{i+p^m}_{l+(p-1)p^{k-(m+1)}}
          \end{equation}
          for $l = 1,\ 2,\ \ldots,\ p^{k-(m+1)}$,\\
          and for
          $i \in \bigcup_{r \in \{1,\ldots,p^m\}}
                   \left\{
                       r+sp^m:\ s \in \{0,1,\ldots,(p^{k-m}-2)\}
                   \right\}$\\
          ($p^m$\ \ $p^{k-m}$--element  cycles of differences).

    \item[$\ldots$]

    \item[the constraints of order $p^{k-1}$ :]
          \begin{equation}
            \label{eq_order_p_to_k_minus_1_constraints}
            \DELTA{i}{i+p^{k-1}}_1 = \DELTA{i}{i+p^{k-1}}_2 =
            \DELTA{i}{i+p^{k-1}}_3 = \ldots =
            \DELTA{i}{i+p^{k-1}}_p
          \end{equation}
          for
          $i \in \bigcup_{r \in \{1,\ldots,p^{k-1}\}}
                   \left\{
                       r+sp^{k-1}:\ s \in \{0,1,\ldots,(p-2)\}
                   \right\}$\\
          ($p^{k-1}$\ \ $p$--element  cycles of differences).

    \item[the dephasing constraints:]
          \begin{equation}
            \label{eq_dephasing_constraints}
            R_{1,1}=R_{2,1}=\ldots=R_{p^k,1}=R_{1,2}=R_{1,3}=\ldots=R_{1,p^k}
            \ \ =\ \ 0
          \end{equation}
  \end{description}

  Then $\SET{R}_{F_{p^k}}$ is a
  $\DEFECT(F_{p^k}) = p^{k-1}\left( \rule{0cm}{0.3cm} (k-1)p - k \right) + 1$
  dimensional subspace of real $p^k \times p^k$ matrices, and if
  $R^{(1)}$, $R^{(2)}$,..., $R^{(\DEFECT(F_{p^k}))}$ form a basis of
  $\SET{R}_{F_{p^k}}$, then
  \begin{equation}
    \label{eq_Fp_family_manifold}
    \left\{
      \rule{0cm}{0.5cm}
       \SET{F}\left(
         \rule{0cm}{0.3cm}
         \phi_1,\ \ldots,\ \phi_{\DEFECT(F_{p^k})}
       \right)
       :\  \phi_i \in \REALS
    \right\}\ ,
  \end{equation}
  where
  \begin{equation}
    \label{eq_Fp_family_manifold_parametrization}
    \SET{F}\left(
      \rule{0cm}{0.3cm}
      \phi_1,\ \ldots,\ \phi_{\DEFECT(F_{p^k})}
    \right)
    \ \ =\ \
    \VECR\left(
      \rule{0cm}{1cm}
      F_{p^k} \HADprod
      \EXPentrywise\left(
        \Ii \sum_{i=1}^{\DEFECT(F_{p^k})}
        \phi_i \cdot R^{(i)}
      \right)
    \right)\ ,
  \end{equation}
  is a $\DEFECT(F_{p^k})$ dimensional manifold (around
  $\VECR(F_{p^k})$) generated by dephased, with respect to $F_{p^k}$,
  unitary complex Hadamard
  matrices, stemming from $\VECR(F_{p^k})$, parametrized by function
  $\SET{F}$ given by (\ref{eq_Fp_family_manifold_parametrization}).
\end{theorem}

{
\newcommand{\levelINDEX}[2]{g^{(#1)}_{#2}}
\PROOFstart
First we show that\ \  $F_{p^k} \HADprod \EXPentrywise(\Ii R)$\ \  is unitary
for any $R \in \SET{R}_{F_{p^k}}$. That is, that for any
$i<j,\ i,j\in \{1,\ldots,p^k\}$, with $n$ further denoting the
difference $j-i$,
the entries of the vector of the summands in the inner product:
\begin{equation}
  \label{eq_FiFj_inner_product_vector}
  \ELEMENTof{
    F_{p^k} \HADprod \EXPentrywise(\Ii R)
  }{i}{:}
  \HADprod
  \ELEMENTof{
    \rule{0cm}{0.4cm}
    \CONJ{ F_{p^k} \HADprod \EXPentrywise(\Ii R) }
  }{j}{:}
  \ =\
  \ELEMENTof{
    \rule{0cm}{0.4cm}
    \CONJ{F_{p^k}}
  }{j-i+1}{:}
  \HADprod
  \EXPentrywise\left(
    \Ii
    \left[
      \rule{0cm}{0.3cm}
      \DELTA{i}{j}_1\
      \DELTA{i}{j}_2\
      \ldots\
      \DELTA{i}{j}_{p^k}
    \right]
  \right),
\end{equation}
where again $\DELTA{i}{j}_l = R_{i,l} - R_{j,l}$,\ \ all add up to zero.

Let $\gcd(p^k,n) = p^m$, i.e. $j-i = ap^m$ with $a,\ p$ relatively
prime. Then, for the 'initial index' 
                        $l=1,2,\ldots,p^{k-(m+1)}$ 
and for the  'factor of rotation by $2\pi$' 
                        \mbox{$r = 0,1,\ldots,(p^m-1)$}, 
we have, within the
$(j-i+1)$-th row of $F_{p^k}$, the groups:
\begin{equation}
  \label{eq_FpTOk_row_zero_groups}
  \ELEMENTof{F_{p^k}}{n+1}{g_1} + \ELEMENTof{F_{p^k}}{n+1}{g_2} +
  \ldots +
  \ELEMENTof{F_{p^k}}{n+1}{g_p}  =  0\ ,
\end{equation}
with
\begin{eqnarray}
  \label{eq_g_indices_definition}
  g_1  &  =  &  \left( \rule{0cm}{0.3cm} l\ +\ 0 \cdot p^{k-(m+1)} \right)\ \  +\ \  rp^{k-m},                   \\
  g_2  &  =  &  \left( \rule{0cm}{0.3cm} l\ +\ 1 \cdot p^{k-(m+1)} \right)\ \  +\ \  rp^{k-m},        \nonumber  \\
       &  \ldots\ ,  &                                                                    \nonumber\\
  g_p  &  =  &  \left( \rule{0cm}{0.3cm} l\ +\ (p-1) \cdot p^{k-(m+1)} \right)\ \  +\ \  rp^{k-m}.        \nonumber
\end{eqnarray}
We are aiming to show that the corresponding groups of $\DELTA{}{}$'s
in (\ref{eq_FiFj_inner_product_vector}) are groups of equal numbers.

Note that $\DELTA{i}{j}_l$'s, for the chosen pair of $i,j$, are
subject to the constraints of order $p^0=1$, order $p^2$, ...,
order $p^m$, \emph{extended} to row index distance $j-i = ap^m$.

Let the indices $g$ in (\ref{eq_g_indices_definition}) be denoted by:
\begin{equation}
  \label{eq_level_1_g_indices}
  \levelINDEX{1}{1} = g_1,\ \ \ \
  \levelINDEX{1}{2} = g_2,\ \ \ \
  \ldots,\ \ \ \
  \levelINDEX{1}{p} = g_p\ .
\end{equation}
Then (where we define\ \  $uv\!\!\! \mod v \stackrel{def}{=} v$\ \ for $u,v \in
\NATURAL$):
\begin{itemize}
  \item
        on account of the constraints of order $1$ 
        $\DELTA{i}{j}_{\levelINDEX{2}{t}} =
        \DELTA{i}{j}_{\levelINDEX{1}{t}}$
        for $t \in \{1,2,\ldots,p\}$, where
        \begin{equation}
          \label{eq_level_2_indices}
          \levelINDEX{2}{t}\ \  =\ \
          \levelINDEX{1}{t}\!\!\! \mod p^{k-1}\ ,
        \end{equation}
        and $\levelINDEX{2}{t}$'s fall into the range covered by 
        the constraints of order $p$,

  \item
        on account of the constraints of order $p$ 
        $\DELTA{i}{j}_{\levelINDEX{3}{t}} =
        \DELTA{i}{j}_{\levelINDEX{2}{t}}$
        for $t \in \{1,2,\ldots,p\}$, where
        \begin{equation}
          \label{eq_level_3_indices}
          \levelINDEX{3}{t}\ \ =\ \
          \levelINDEX{2}{t}\!\!\! \mod p^{k-2}\ ,
        \end{equation}
        and $\levelINDEX{3}{t}$'s fall into the range covered by the
        constraints of order $p^2$,

  \item
        ...

  \item
        on account of the constraints of order $p^{m-1}$ 
        $\DELTA{i}{j}_{\levelINDEX{m+1}{t}} =
        \DELTA{i}{j}_{\levelINDEX{m}{t}}$
        for $t \in \{1,2,\ldots,p\}$, where
        \begin{equation}
          \label{eq_level_mPLUS1_indices}
          \levelINDEX{m+1}{t}\ \  =\ \
          \levelINDEX{m}{t}\!\!\! \mod p^{k-m}\ ,
        \end{equation}
        and $\levelINDEX{m+1}{t}$'s fall into the range covered by the
        constraints of order $p^m$
        .
\end{itemize}

It is easy to notice that numbers $\levelINDEX{1}{1}$,
$\levelINDEX{1}{2}$, ..., $\levelINDEX{1}{p-1}$ are all not divided by
$p^{k-m}$, $p^{k-m+1}$, ..., $p^{k-1}$.
Thus at each step of the above procedure the sequence
$\levelINDEX{u}{1},\ \ldots,\ \levelINDEX{u}{p-1},\ \levelINDEX{u}{p}$
remains a sequence of the appropriate remainders equally spaced by
$p^{k-(m+1)}$, so the final sequence has this form:
\begin{eqnarray}
  \label{eq_final_g_sequence}
  \levelINDEX{m+1}{1}  &  =  &  \tilde{l}\ +\ 0 \cdot p^{k-(m+1)}\ ,           \\
  \levelINDEX{m+1}{2}  &  =  &  \tilde{l}\ +\ 1 \cdot p^{k-(m+1)}\ , \nonumber \\
  &  \ldots\ \ ,  &                                         \nonumber \\
  \levelINDEX{m+1}{p}  &  =  &  \tilde{l}\ +\ (p-1) \cdot p^{k-(m+1)}\ ,  \nonumber
\end{eqnarray}
for some $\tilde{l} \in \{1,2,\ldots,p^{k-(m+1)}\}$.

On account of the constraints of order $p^m$  all
$\DELTA{i}{j}_{\levelINDEX{m+1}{t}}$\ \ for $t\in\{1,...,p\}$ are
equal.

Hence we have obtained that for any $l \in \{1,2,\ldots,p^{k-(m+1)}\}$, and
for any $r \in \{0,1,\ldots,p^m-1\}$, with $g_t$'s defined in
(\ref{eq_g_indices_definition}):
\begin{equation}
  \label{eq_required_DELTA_equality}
  \DELTA{i}{j}_{g_1}\ \  =\ \  \DELTA{i}{j}_{g_2}\ \  =\ \
  \ldots\ \  =\ \
  \DELTA{i}{j}_{g_p}\ .
\end{equation}
Taking (\ref{eq_FpTOk_row_zero_groups}) into account, the entries of
the vector of the summands in the inner
product (\ref{eq_FiFj_inner_product_vector}) add up
to zero, which confirms the unitarity of\ \
$F_{p^k} \HADprod \EXPentrywise(\Ii R)$\ \
for any $R \in \SET{R}_{F_{p^k}}$.

The number of independent equations, imposing the constraints of order
$1$, order $p$, ..., order $p^{k-1}$  
on $R$ for it to belong to
$\SET{R}_{F_{p^k}}$ reads:
\begin{eqnarray}
  \label{eq_no_of_difference_constraints}
  p^{k-1}     (p-1) \left( 1       (    p^k - 1) \right)   &   +   &  \\
  p^{k-2}     (p-1) \left( p       (p^{k-1} - 1) \right)   &  +    &              \nonumber \\
  p^{k-3}     (p-1) \left( p^2     (p^{k-2} - 1) \right)   &  +\ \  \ldots\ \ +  &  \nonumber \\
  p^{k-(m+1)} (p-1) \left( p^m     (p^{k-m} - 1) \right)   &   +\ \ \ldots\ \ +  & \nonumber  \\
            1 (p-1) \left( p^{k-1} (      p - 1) \right)\ ,&  &                     \nonumber
\end{eqnarray}
which is equal to:
\begin{eqnarray}
  {
    (p-1)
    \left(
      \rule{0cm}{0.5cm}
      p^k
      \left(
        \rule{0cm}{0.4cm}
        p^{k-1} + p^{k-2} + p^{k-3} + \ldots + p^{k-(m+1)} +         \ldots + 1
      \right)
      \ \ -\ \  p^{k-1} \cdot k
    \right)
    \ \ \ =
  }         &  &  \nonumber  \\
  (p^k)^2\  -\  (k+1)p^k\  +\  k \cdot p^{k-1}\ .\ \ \ \ \ \ \ \ \
  &  &    \label{eq_no_of_difference_constraints_final}
\end{eqnarray}
Taking also the dephasing constraints into consideration, the number of
independent parameters in $R \in \SET{R}_{F_{p^k}}$ is the difference:
\begin{eqnarray}
  \lefteqn{
    (p^k)^2\ \  -\ \  (2p^k-1)\ \  -\ \
    \left(
      \rule{0cm}{0.4cm}
      (p^k)^2\  -\  (k+1)p^k\  +\  k \cdot p^{k-1}
    \right)\ \
  =}  \nonumber  &  &  \\
  &  p^{k-1}
     \left(
       \rule{0cm}{0.4cm}
       (k-1)p \ -\  k
     \right)
     \ \  +\ \ 1
  & =\ \  \DEFECT(F_{p^k})\ ,         \label{eq_FpTOk_family_phase_matrix_space_dim}
\end{eqnarray}
in accordance with formula (\ref{eq_WS_main_formula}). This gives us
the dimension of $\SET{R}_{F_{p^k}}$.

Function $\SET{F}$ of (\ref{eq_Fp_family_manifold_parametrization})
parametrizes a $\DEFECT(F_{p^k})$ dimensional manifold around
$\VECR(F_{p^k})$, since
\begin{equation}
  \label{eq_Fp_family_parametrization_derivative}
  \left.
    \frac{\delta}{\delta\phi_i}
    \SET{F}\left(
      \rule{0cm}{0.3cm}
      \phi_1,\ \ldots,\ \phi_{\DEFECT(F_{p^k})}
    \right)
  \right|_{\mathbf{\phi} = \ZEROvect}
  \ \ =\ \
  \VECR( \Ii R^{(i)} \HADprod F_{p^k} )
\end{equation}
are independent vectors in $\REALS^{2(p^k)^2}$. In fact, a similar
argument leads to the conclusion that $\SET{F}$ of
(\ref{eq_Fp_family_manifold_parametrization})
parametrizes a manifold around any point
$\SET{F}\left( \rule{0cm}{0.3cm} \phi_1,\ldots,\phi_{\DEFECT(F_{p^k})}
\right)$.

Note that although we do not show here explicitly that constraints
defining $\SET{R}_{F_{p^k}}$ are independent, the dimension of
$\SET{R}_{F_{p^k}}$ cannot exceed $\DEFECT(F_{p^k})$, as the
dimension of the manifold generated with  $\SET{R}_{F_{p^k}}$,  
(\ref{eq_Fp_family_manifold}), cannot be greater than that, according
to Theorem \ref{theor_manifold_dim_defect_bound}.

\PROOFend
}

As examples, let us examine the $\DEFECT(F_N)$--dimensional families
steming from $\VECR(F_{2^3})$ and $\VECR(F_{3^2})$.
We present both forms,
featured in Theorems
\ref{theor_Fp_family_PCM_construction}
and
\ref{theor_Fp_family_direct_construction}, 
of these.

\begin{eqnarray}
  \label{eq_F8_family}
  \lefteqn{
    \SET{F}_{F_8}(a,b,c,d,e)\ \ \  =\ \ \
  }  &  &  \\
  &
  \left\{
    \rule{0cm}{0.5cm}
    \VECR\left(
      \rule{0cm}{0.4cm}
      F_{8} \HADprod \EXPentrywise(\Ii \cdot P_8(a,b,c,d,e) \cdot F_{8})
    \right)
    \ :  \ \
    a,b,c,d,e \in \REALS
  \right\}
  &  =      \nonumber  \\
  &
  \left\{
    \rule{0cm}{0.5cm}
    \VECR\left(
      \rule{0cm}{0.4cm}
      F_{8} \HADprod \EXPentrywise(\Ii \cdot R_8(a,b,c,d,e) )
    \right)
    \ :\ \
    a,b,c,d,e \in \REALS
  \right\}\ ,
  &        \nonumber
\end{eqnarray}
where
\begin{equation}
  \label{eq_P8_phases_factor_for_F8}
  P_8(a,b,c,d,e)\ \ \  =\ \ \
  \left[
    \begin{array}{c|ccc|c|ccc}
      0       &   0 & 0           & 0 &       0 &       0 & 0           & 0 \\
      -(2a+c) &   0 & (a + \Ii b) & 0 &       c &       0 & (a - \Ii b) & 0 \\
      -d      &   0 & 0           & 0 &       d &       0 & 0           & 0 \\
      -(2a+e) &   0 & (a + \Ii b) & 0 &       e &       0 & (a - \Ii b) & 0 \\
      0       &   0 & 0           & 0 &       0 &       0 & 0           & 0 \\
      -(2a+c) &   0 & (a + \Ii b) & 0 &       c &       0 & (a - \Ii b) & 0 \\
      -d      &   0 & 0           & 0 &       d &       0 & 0           & 0 \\
      -(2a+e) &   0 & (a + \Ii b) & 0 &       e &       0 & (a - \Ii b) & 0
    \end{array}
  \right]
\end{equation}
and
\begin{equation}
  \label{eq_R8_phases_for_F8}
   R_8(a,b,c,d,e)\ \ \  =\ \ \
   \left[
     \begin{array}{cccc|cccc}
       \bO & \bO & \bO & \bO & \bO & \bO & \bO & \bO\\
       \bO & a & b & c & \bO & a & b & c\\
       \bO & d & \bO & d & \bO & d & \bO & d\\
       \bO & e & b & c - a + e & \bO & e & b & c - a + e\\
       \hline
       \bO & \bO & \bO & \bO & \bO & \bO & \bO & \bO\\
       \bO & a & b & c & \bO & a & b & c\\
       \bO & d & \bO & d & \bO & d & \bO & d\\
       \bO & e & b & c - a + e & \bO & e & b & c - a + e
     \end{array}
   \right] \ .
\end{equation}

\begin{eqnarray}
  \label{eq_F9_family}
  \lefteqn{
    \SET{F}_{F_9}(a,b,c,d)\ \ \  =\ \ \
  }  &  &  \\
  &
  \left\{
    \rule{0cm}{0.5cm}
    \VECR\left(
      \rule{0cm}{0.4cm}
      F_{9} \HADprod \EXPentrywise(\Ii \cdot P_9(a,b,c,d) \cdot F_{9})
    \right)
    \ :  \ \
    a,b,c,d \in \REALS
  \right\}
  &  =      \nonumber  \\
  &
  \left\{
    \rule{0cm}{0.5cm}
    \VECR\left(
      \rule{0cm}{0.4cm}
      F_{9} \HADprod \EXPentrywise(\Ii \cdot R_9(a,b,c,d) )
    \right)
    \ :\ \
    a,b,c,d \in \REALS
  \right\}\ ,
  &        \nonumber
\end{eqnarray}
where
\begin{equation}
  \label{eq_P9_phases_factor_for_F9}
  P_9(a,b,c,d)\ \ \  =\ \ \
  \left[
    \begin{array}{c|cccc|cccc}
        0 &    0 & 0 &           0 & 0 &    0 &           0 & 0 & 0 \\
      -2a &    0 & 0 & (a + \Ii b) & 0 &    0 & (a - \Ii b) & 0 & 0 \\
      -2c &    0 & 0 & (c + \Ii d) & 0 &    0 & (c - \Ii d) & 0 & 0 \\
        0 &    0 & 0 &           0 & 0 &    0 &           0 & 0 & 0 \\
      -2a &    0 & 0 & (a + \Ii b) & 0 &    0 & (a - \Ii b) & 0 & 0 \\
      -2c &    0 & 0 & (c + \Ii d) & 0 &    0 & (c - \Ii d) & 0 & 0 \\
        0 &    0 & 0 &           0 & 0 &    0 &           0 & 0 & 0 \\
      -2a &    0 & 0 & (a + \Ii b) & 0 &    0 & (a - \Ii b) & 0 & 0 \\
      -2c &    0 & 0 & (c + \Ii d) & 0 &    0 & (c - \Ii d) & 0 & 0 \\
    \end{array}
  \right]
\end{equation}
and
\begin{equation}
  \label{eq_R9_phases_for_F9}
   R_9(a,b,c,d)\ \ \  =\ \ \
   \left[
    \begin{array}{ccc|ccc|ccc}
      \bO & \bO & \bO & \bO & \bO & \bO & \bO & \bO & \bO\\
      \bO & a & b & \bO & a & b & \bO & a & b\\
      \bO & c & d & \bO & c & d & \bO & c & d\\
      \hline
      \bO & \bO & \bO & \bO & \bO & \bO & \bO & \bO & \bO\\
      \bO & a & b & \bO & a & b & \bO & a & b\\
      \bO & c & d & \bO & c & d & \bO & c & d\\
      \hline
      \bO & \bO & \bO & \bO & \bO & \bO & \bO & \bO & \bO\\
      \bO & a & b & \bO & a & b & \bO & a & b\\
      \bO & c & d & \bO & c & d & \bO & c & d
    \end{array}
  \right] \ .
\end{equation}

Reasoning very much like in the proof of Theorem
\ref{theor_zero_defect_isolation}, one can prove the fact stated below 
about the discussed continuous families stemming from $F_{p^k}$. By 'dephased'
matrices we mean dephased with respect to $F_{p^k}$ in the manner
described in the introductory part of this section.


\begin{theorem}
  \label{theor_FpTOk_family_only_Hadamards}

  There exists a neighbourhood $\SET{W}$ of $\VECR(F_{p^k})$ in $\REALS^{2(p^k)^2}$
  such that the only vectors
  $v \in \SET{W} \setminus \{\VECR(F_{p^k})\}$
  generated by dephased unitary complex Hadamard matrices:
  \begin{equation}
    \label{eq_v_potential_Hadamards_around_FpTOk}
    v = \VECR\left( F_{p^k} \HADprod \EXPentrywise( \Ii R ) \right),
    \ \ \ \ \ \ \mbox{where}\ \ \ \
    R_{1,j} = R_{i,1}= 0,\ \ \ i,j \in \{1..p^k\},
  \end{equation}
  are those generated by members of the continuous
  $\DEFECT(F_{p^k})$-dimensional family of
  Theorem \ref{theor_Fp_family_direct_construction}:
  \begin{equation}
    \label{eq_FpTOk_family}
    \left\{
        F_{p^k} \HADprod \EXPentrywise( \Ii R ):\ R \in
        \SET{R}_{F_{p^k}}
    \right\}\ ,
  \end{equation}
  where $\SET{R}_{F_{p^k}}$ is the $\DEFECT(F_{p^k})$-dimensional
  linear space defined in Theorem \ref{theor_Fp_family_direct_construction}.
\end{theorem}

\PROOFstart 
Unitarity of
$\VECR^{-1}(v) = F_{p^k} \HADprod \EXPentrywise( \Ii R )$
(see (\ref{eq_v_potential_Hadamards_around_FpTOk}))
and the dephasing condition can be expressed, for $\VEC(R) \in
\REALS^{(p^k)^2}$, as the system of equations:
\begin{equation}
  \label{eq_v_potential_Hadamards_conditions}
  \left\{
    \begin{array}{cccc}
      R_{1,j} & = & 0  &\ \ \ \   j \in \{2,\ldots,p^k\} \\
      R_{i,1} & = & 0  &\ \ \ \   i \in \{1,\ldots,p^k\} \\
      g\left( \rule{0cm}{0.3cm} \VEC(R) \right) & = & \ZEROvect &
    \end{array}\ \ ,
  \right.
\end{equation}
where $g$ is defined at the end of Section
\ref{subsec_other_defect_characterizations}
in (\ref{eq_g_def}), and where $U$ is taken to be $F_{p^k}$. The
collective system will be denoted by
\begin{equation}
  \label{eq_h_def}
  h\left( \rule{0cm}{0.3cm} \VEC(R) \right)\ \  =\ \  \ZEROvect\ .
\end{equation}

Looking at the form  of the differential of $g$ at $\ZEROvect$
(see (\ref{eq_Dg_of_vecR}) and the description there)
we notice that the differential of $h$ at
$\ZEROvect$ satisfies:
\begin{equation}
  \label{eq_Dh_rank}
  \dim\left(
    \rule{0cm}{0.5cm}
    \DIFFERENTIAL{h}{\ZEROvect}
    \left(
      \rule{0cm}{0.3cm}
      \REALS^{(p^k)^2}
    \right)
  \right)
  \ \ =\ \
  (p^k)^2 - \DEFECT(F_{p^k}),
\end{equation}
as $\VEC(\STbasis{k}\ONESvect^T)$, $k=1..p^k$ and
   $\VEC(\ONESvect\STbasis{l}^T)$, $l=2..p^k$,
spanning a part of the kernel of $\DIFFERENTIAL{g}{\ZEROvect}$,
are no longer in the kernel of $\DIFFERENTIAL{h}{\ZEROvect}$.

Thus one can choose a
$\left( \rule{0cm}{0.4cm} (p^k)^2 - \DEFECT(F_{p^k})
\right)$--equation subsystem
of (\ref{eq_h_def}) with the full rank:
\begin{equation}
  \label{eq_h_subsystem_def_and_property}
  \tilde{h}\left(
    \VEC(R)
  \right)\ \  =\ \ \ZEROvect,
  \ \ \ \ \ \ \ \ \mbox{where}\ \ \ \
  \dim\left(
    \rule{0cm}{0.5cm}
    \DIFFERENTIAL{\tilde{h}}{\ZEROvect}
    \left(
      \rule{0cm}{0.3cm}
      \REALS^{(p^k)^2}
    \right)
  \right)
  \ \ =\ \ (p^k)^2 - \DEFECT(F_{p^k})\ .
\end{equation}

Therefore system (\ref{eq_h_subsystem_def_and_property}) defines a
$\DEFECT(F_{p^k})$-dimensional manifold around $\ZEROvect$, and this
must be the $\DEFECT(F_{p^k})$-dimensional linear space
$\VEC(\SET{R}_{F_{p^k}})$ with $\SET{R}_{F_{p^k}}$ defined in Theorem
\ref{theor_Fp_family_direct_construction}, for it satisfies
(\ref{eq_h_subsystem_def_and_property}).

If $v \in \SET{W}$ and $\VECR^{-1}(v)$ is a dephased unitary complex
Hadamard matrix, then $\VECR^{-1}(v) = F_{p^k} \HADprod
\EXPentrywise(\Ii R)$ for some $\VEC(R)$ close to $\ZEROvect$, and $\VECR(R)$
must satisfy (\ref{eq_h_subsystem_def_and_property}), i.e.
$R \in \SET{R}_{F_{p^k}}$.
\PROOFend

%
%

\section{Conclusions}
\label{sec_conclusions}

In this work we proposed a definition of the defect of
a unitary matrix of size $N$. This notion is shown to be useful
while investigating certain properties of unitary matrices.
Demonstrating that the defect of any Fourier matrix of a prime size
is equal to zero we infer that in this case $F_N$ is an isolated
unitary complex Hadamard matrix. This result also allows us to
prove that for prime dimensions there exists a unistochastic
ball around the flat \BISTOCHASTIC\ matrix $J_N$.

A positive value of the defect of $F_N$ for a composite $N$
provides a direct upper bound for the dimension
of an orbit of dephased (and thus locally $\RELofEQUI$-inequivalent)
unitary complex Hadamard matrices.
Already for $N=6$ this bound, equal to $4$, is larger than the
dimension $D=2$ of the largest orbit known, which may suggest
that the list of known Hadamard matrices is incomplete.

The defect of any $U$ may be expressed using the rank of certain
matrix associated with $U$ and computed numerically. Such computations
were performed for several unitaries
of size $N=6$ belonging to the known families of inequivalent
unitary complex Hadamard matrices. In all cases studied
the defect was equal to $\DEFECT(F_{6}) = 4$,
which provides a hint \cite{Bxxx07} that these families
may be embedded inside an unknown orbit of dimension $4$.
This reasoning allows us to believe that the notion of the defect
will be useful in further search for new families
of (unitary) complex Hadamard matrices.

In this paper we presented two constructions of $\DEFECT(F_N)$-dimensional smooth
families of inequivalent complex Hadamard matrices which stem from the
Fourier matrix  $F_N$. These constructions work
 for $N$ being a power of a prime number.
 One of them involves the 'parameter cycle matrices', which proved to be useful by 
computing the defect $\DEFECT(F_N)$. 
The family of complex Hadamard matrices obtained in this way 
has a particularly nice form which is 
due to the symmetric structure of $F_N$,
 and is closely related to the fact that $F_N$ diagonalizes circulant matrices of size $N$. 
Analogous properties of orbits of inequivalent matrices stemming
from tensor products of Fourier matrices need further investigations
for other composite $N$ which are not a power of prime.

The defect of a unitary matrix is related to the map
(\ref{eq_f0}) projecting the $N^2$ dimensional set of unitary matrices
into the $(N-1)^2$ dimensional set of unistochastic matrices.
The actual value of the defect provides a kind of characterization of
the space of unitary matrices and allows one to classify
its elements. For a generic unitary matrix $\DEFECT(U) = 0$,
while any deviation from this value for a given
$U$ confirms certain special properties of the analyzed matrix.
For instance, we find that the defect of a generic
real orthogonal matrix $O_N$ of size $N>2$ is positive
%
%
and satisfies $\DEFECT(O_N) \geq (N-1)(N-2)/2$.
Although we have some knowledge
on the defect of unitary matrices with a tensor product
structure \cite{TK07}, the general problem of characterizing
a class of unitary matrices of size $N$ with a fixed defect
remains open.

%
%

\section*{Acknowledgments}

\medskip
It is a pleasure to thank Wojciech S{\l}omczy{\'n}ski
for proving an alternative  formula  (\ref{eq_WS_main_formula})
for the defect of the Fourier matrix and for writing it down in
appendix B.
We enjoyed numerous fruitful discussions with
I. Bengtsson, W. Bruzda, {\AA}. Ericsson and  M. Ku{\'s}.
We  are also grateful to P.~Di{\c t}{\v a},
M. Matolcsi, R. Nicoara, and  F. Sz\"{o}ll\H{o}si
for helpful correspondence and for
providing us with their results prior to publication.
We acknowledge financial support by
Polish Ministry of Science and Information Technology
under the grant 1\, P03B\, 042\, 26 
and by the European Research Project SCALA.
\appendix

%
%

\section{Notation}
\label{sec_notation}

We shall adopt the following conventions:
\begin{description}
  \item[$A \HADprod B$]
        denotes the Hadamard product of matrices $A$ and $B$

  \item[$\diag(v)$]
        denotes an $N \times N$ diagonal matrix, for an $N$ element
        vector or sequence $v$,
        such that
        $[{\rm {\bf diag}}(v)]_{i,i}=v_i$

  \item[$\EXPentrywise(A)$]
        denotes the entrywise operation $\exp$ on matrix $A$

  \item[$\Re,\ \Im$]
        denotes also the entrywise  operations $\Re,\ \Im$ on matrices

  \item[$A_{i_1:i_2,j_1:j_2}$]
        denotes a sub--matrix of matrix $A$
        \begin{displaymath}
          A_{i_1:i_2,j_1:j_2} =
          \left[
            \begin{array}{ccc}
              A_{i_1,j_1} & \ldots & A_{i_1,j_2} \\
              \ldots      & \ldots & \ldots      \\
              A_{i_2,j_1} & \ldots & A_{i_2,j_2}
            \end{array}
          \right]
        \end{displaymath}
        If $i_1 = i_2$ or $j_1 = j_2$ we write $i_1,\ j_1$ instead of
        $i_1:i_2,\ j_1:j_2$, respectively.

  \item[$\STbasis{k}$]
        denotes the $k$-th standard basis column vector

  \item[$\ONESvect$]
        denotes a vertical vector $[1,1,\ldots,1]^T$ filled all with ones

 \item[$\VECC(A)\    \left( \rule{0cm}{0.4cm} \VEC(A) \right)$]
       denotes the 'row by row' vertical complex (real) vector form of a
        complex (real) $N \times N$ matrix $A$:
        \begin{displaymath}
          \VECC(A) = [ A_{1,1}, \ldots, A_{1,N},
                       A_{2,1}, \ldots, A_{2,N}, \ldots
                       A_{N,1}, \ldots, A_{N,N}        ]^T
        \end{displaymath}

        We identify $\REALS^k$ ($\COMPLEX^k$) with the set of all real
        (complex) vertical $k \times 1$ vectors (matrices).

  \item[$\VECR(A)$]
        denotes the 'row by row' vertical real vector form of a
        complex $N \times N$ matrix $A$:
        \begin{displaymath}
          \VECR(A) =
          \left[
            \begin{array}{c}
              \Re( \VECC(A) ) \\
              \Im( \VECC(A) )
            \end{array}
          \right]
        \end{displaymath}

 \item[$\SPANC(S)\    \left( \rule{0cm}{0.4cm} \SPANR(S) \right)$]
      denotes a complex (real) linear space spanned by vectors from a set or
      columns of a matrix $S$

 \item[$\NULLC(D)\    \left( \rule{0cm}{0.4cm} \NULLR(D) \right)$]
        denotes, for an operator or complex matrix $D$,
       the complex space  $\{ v \in \COMPLEX^N:\ D(v)=\ZEROvect \}$
        $\left(\rule{0cm}{0.4cm}\right.$the real space $\{ r \in
        \REALS^N:\ D(r)=\ZEROvect \}$$\left.\rule{0cm}{0.4cm}\right)$,
         for a given $N$

  \item[$\BISTOCHSPACE$]
        denotes the set of all real matrices with all row and column sums
        equal to $1$, for a given size $N$; this includes \BISTOCHASTIC\ matrices
        which contain non--negative entries only

  \item[$\UNITARY$]
        denotes the set of all unitary matrices, for a given size $N$

  \item[$\alpha$]
        denotes the following function  generating indices into a
        matrix
        \begin{displaymath}
          \alpha\ :\ \left\{ \rule{0cm}{0.5cm} (i,j):\ 1 \leq i < j \leq N \right\}
                     \longrightarrow
                     \left\{ 1,\ 2,\ \ldots,\ \frac{(N-1)N}{2} \right\}
        \end{displaymath}
        such that
        \begin{displaymath}
          \begin{array}{c||c|c|c|c|c|c|c|c}
            (i,j) & (1,2) & \ldots & (1,N) & (2,3) & \ldots & (2,N) & \ldots & (N-1,N) \\
            \hline
    \alpha(i,j) &   1   & \ldots &   N-1   &  N  & \ldots & 2N-3  & \ldots & \frac{(N-1)N}{2}
          \end{array}
        \end{displaymath}
\end{description}

(\ref{eq_WS_aux_formula}).

%
%

\section{Proof of Theorem
  \ref{theor_WS_Fourier_matrix_defect_formula} on alternate 
  formulae for the defect of the Fourier matrix $F_N$ \\
  (by Wojciech S\l omczy\'nski)}

\label{sec_proof}

a. From Theorem \ref{theor_Fourier_matrix_defect} and from the symmetry relations%
\[
\sum_{l=1}^{\frac{N-1}{2}}\gcd\left(  N,l\right)  =\frac{1}{2}\left(
\sum_{l=1}^{N-1}\gcd\left(  N,l\right)  \right)  \text{ ,}%
\]
for odd $N$, and%
\[
\sum_{l=1}^{\frac{N}{2}-1}\gcd\left(  N,l\right)  =\frac{1}{2}\left(
\sum_{l=1}^{N-1}\gcd\left(  N,l\right)  -\frac{N}{2}\right)  \text{ ,}%
\]
for even $N$, we deduce (\ref{eq_WS_aux_formula}).

b. From (\ref{eq_WS_aux_formula}) we get%
\begin{align}
\mathbf{d}(F_{N}) &  =1-2N+\sum_{l=1}^{N}\gcd\left(  N,l\right)  \nonumber\\
&  =1-2N+\sum_{d|N}d\cdot\left|  \left\{  l:1\leq l\leq N\text{ and }%
\gcd\left(  l,N\right)  =d\right\}  \right|  \nonumber\\
&  =1-2N+\sum_{d|N}d\cdot\left|  \left\{  k:1\leq k\leq N/d\text{ and }%
\gcd\left(  k,N/d\right)  =1\right\}  \right|  \nonumber\\
&  =1-2N+\sum_{d|N}d\cdot\varphi\left(  N/d\right)  \nonumber\\
&  =N\left(  \sum_{d|N}\frac{d}{N}\varphi\left(  N/d\right)  -2\right)
+1\nonumber\\
&  =N\left(  \sum_{d|N}\psi\left(  N/d\right)  -2\right)  +1\nonumber\\
&  =N\left(  \psi_{\leq}\left(  N\right)  -2\right)  +1\text{ ,}%
\label{defect1}%
\end{align}
where $\varphi$ is the Euler function \cite[p. 158]{Aig79} given by:%
\[
\varphi\left(  M\right)  :=\left|  \left\{  l:1\leq k\leq M\text{ and }%
\gcd\left(  k,M\right)  =1\right\}  \right|  \text{ ,}%
\]
$\psi$ is an arithmetic function defined by:%
\[
\psi\left(  M\right)  :=\varphi\left(  M\right)  /M\text{ ,}%
\]
and the M\"{o}bius inverse function $\psi_{\leq}$ is given by:%
\[
\psi_{\leq}\left(  M\right)  :=\sum_{d|M}\psi\left(  \frac{M}{d}\right)
\text{ .}%
\]
\qquad\qquad\qquad\qquad\qquad\qquad

We shall show that%
\begin{equation}
\psi_{\leq}\left(  N\right)  =\prod_{j=1}^{n}\left(  1+k_{j}-\frac{k_{j}%
}{p_{j}}\right)  =:R\left(  N\right)  \text{ .}\label{PsiMobinv}%
\end{equation}
To prove (\ref{PsiMobinv}) it suffices to apply the M\"{o}bius inversion
formula \cite[p. 154]{Aig79}%
\[
\psi\left(  N\right)  :=\sum_{d|N}\psi_{\leq}\left(  \frac{N}{d}\right)
\mu\left(  d\right)
\]
and to show that%
\begin{equation}
\psi\left(  N\right)  :=\sum_{d|N}R\left(  \frac{N}{d}\right)  \mu\left(
d\right)  \text{ ,}\label{Psi}%
\end{equation}
where $\mu$ is the M\"{o}bius function defined as%
\[
\mu\left(  d\right)  :=\left\{
\begin{tabular}
[c]{ll}%
$\left(  -1\right)  ^{s}$ & $d=\prod_{j=1}^{s}p_{j}$ , where $p_{j}$ $\left(
j=1,\ldots,s\right)  $ are different primes\\
$0$ & otherwise
\end{tabular}
\right.  \text{ .}%
\]
From the Euler formula \cite[p. 158]{Aig79}%
\[
\psi\left(  N\right)  =\varphi\left(  N\right)  /N=\prod_{j=1}^{n}\left(
1-\frac{1}{p_{j}}\right)
\]
we deduce that%
\begin{align*}
\psi\left(  N\right)   &  =\prod_{j=1}^{n}\left(  1-\frac{1}{p_{j}}\right)  \\
&  =\sum\limits_{b\in\left\{  0,1\right\}  ^{n}}\prod_{j=1}^{n}\left(
1+k_{j}-b_{j}-\frac{k_{j}-b_{j}}{p_{j}}\right)  \left(  -1\right)  ^{b_{j}}\\
&  =\sum\limits_{b\in\left\{  0,1\right\}  ^{n}}R\left(  \prod_{j=1}^{n}%
p_{j}^{k_{j}-b_{j}}\right)  \left(  -1\right)  ^{\sum_{j=1}^{n}b_{j}}\\
&  =\sum_{d|N}R\left(  \frac{N}{d}\right)  \mu\left(  d\right)  \text{ ,}%
\end{align*}
which proves (\ref{Psi}), and, in consequence, (\ref{PsiMobinv}).

Now, formula (\ref{eq_WS_main_formula}) folllows from (\ref{PsiMobinv}) and
(\ref{defect1}).

%
%


\begin{thebibliography}{99}

\bibitem{MO79}
      A.W. Marshall and  I. Olkin,
      {\sl Inequalities: Theory of Majorization and its Applications},
      Academic Press, New York, 1979.

\bibitem{ZKSS03}
      K. {\.Z}yczkowski, M. Ku{\'s},
      W. S{\l}omczy{\'n}ski  and H.-J. Sommers,
      Random unistochastic matrices,
      {\sl J. Phys.} {\bf A 36}, 3425-3450 (2003).

\bibitem{AMM91}
      G. Auberson, A. Martin and G. Mennessier,
      Commun. Math. Phys. {\bf 140}, 523 (1991).

\bibitem{BEKTZ05}
      I. Bengtsson, A.~Ericsson, M.~Ku{\'s},
      W.~Tadej, and  K.~{\.Z}yczkowski,
      Birkhoff's polytope and unistochastic matrices, $N=3$ and $N=4$,
      {\sl Comm. Math. Phys.} {\bf 259},  307-324 (2005).

\bibitem{Cr91} R. Craigen,
    Equivalence classes of inverse orthogonal and unit Hadamard,
    Bull. Austral. Math. Soc. 44, 109-115 (1991).

\bibitem{Po83}
      S. Popa,
      Orthogonal pairs of *-subalgebras in finite von Neumann
      algebras,
      J. Operator Theory 9 (1983),
      pp 253 - 268.

\bibitem{MW92}
      A. Munemasa and Y. Watatani, Orthogonal pairs of
      $*$--subalgebras and association schemes,
      C.R. Acad. Sci. Paris {\bf 314}, 329-331 (1992).

\bibitem{Ha96}
      U. Haagerup,
      Orthogonal maximal abelian $*$-subalgebras of the $n \times n$
      matrices and cyclic $n$--roots,
      Operator Algebras and Quantum Field Theory (Rome),
      1996 (Cambridge, MA: International Press),
      pp 296-322.

\bibitem{BF91}
      G. Bj{\"o}rk and R. Fr{\"o}berg,
      A faster way to count the solutions of inhomogeneous systems of
      algebraic equations,   with applications to cyclin $n$--roots,
      J. Symbolic Comp. {\bf 12}, 329-336 (1991).

\bibitem{BS95}
      G. Bj{\"o}rck and B. Saffari,
      New classes of finite unimodular sequences with  unimodular
      Fourier transform.  Circulant Hadamard matrices with complex  entries,
      {\sl C. R. Acad. Sci., Paris} {\bf 320} 319-24  (1995).

\bibitem{GR05}
      C. D. Godsil and  A. Roy,
      Equiangular lines, mutually unbiased bases, and spin models
      preprint {\sl www.arxiv.org/abs/quant-ph/0511004}

\bibitem{Di04}
      P. Di{\c t}{\v a},
      Some results on the parametrization of complex Hadamard matrices,
      J. Phys. A: Math. Gen. {\bf 37}, 5355-5374 (2004),

\bibitem{We01}
      R.F. Werner,
      All teleportation and dense coding schemes,
      J. Phys. A: Math. Gen. {\bf 34} 7081-94  (2001),

\bibitem{WGC03}
      A. W{\'o}jcik, A. Grudka and R.W. Chhajlany,
      Generation of inequivalent generalized Bell bases,
      Quantum Information Processing, 2, 201 (2003).

\bibitem{TZ06}
      W. Tadej and K. \.Zyczkowski,
      A concise guide to complex Hadamard matrices,
      Open Systems \& Infor. Dynamics 13, 133-177 (2006).
      For an updated version of the catalog see also
      {\sl http://chaos.if.uj.edu.pl/${\sim}$karol/hadamard}

\bibitem{Ja85}
      C. Jarlskog,
      Commutator of the Quark Mass Matrices in the
      Standard Electroweak Model and a Measure of Maximal CP Nonconservation
      Phys. Rev. Lett. 55, 1039-1042 (1985).

\bibitem{BD87}
      J. D. Bjorken and  I. Dunietz,
      Rephasing-invariant parametrizations of
      generalized Kobayashi-Maskawa matrices
      Phys. Rev. D 36, 2109-2118 (1987).

\bibitem{Di05}
      P. Di{\c t}{\v a},
      Global fits to the Cabibbo-Kobayashi-Maskawa matrix:
      Unitarity condition method versus standard unitarity triangles approach,
      Mod. Phys. Lett. {\bf A 20},  1709-1721 (2005).

\bibitem{KS97}
      T. Kottos and U. Smilansky,
      Quantum chaos on graphs,
      Phys. Rev. Lett. {\bf 79}, 4794-4797 (1997).

\bibitem{Ta01}
      G. Tanner,
      Unitary-stochastic matrix ensembles and spectral statistics,
      {\sl J. Phys. A} {\bf 34}, 8485-8500 (2001).

\bibitem{PZK01}
      P. Pako{\'n}ski, K. {\.Z}yczkowski and  M. Ku{\'s},
      Classical 1D maps, quantum graphs and ensembles of unitary matrices,
      J. Phys. {\bf A 34}, 9303-9317 (2001).

\bibitem{Ka89} 
        A. Karabegov, 
    The reconstruction of a unitary matrix from the moduli 
   of its elements and symbols on finite phase space,
preprint YERPHI-1194 (71)-89, Yerevan (1989).

 \bibitem{TK07}
       W. Tadej, M. Ku{\'s} et al.
       Defect of unitary matrices of a tensor product structure,
       in preparation.

\bibitem{Nicoara04}
      R. Nicoara,
      A finiteness result for commuting squares of matrix algebras,
      preprint {\sl www.arxiv.org/abs/math.OA/0404301}

\bibitem{Petrescu97}
      M. Petrescu,
      Existence of continuous families of complex Hadamard matrices of
      certain  prime dimensions,
      Ph.D thesis, UCLA 1997.

\bibitem{Ta06} W. Tadej,
      Permutation equivalence classes of Kronecker products
      of unitary Fourier matrices,
      Linear Algebra Appl. {\bf 418}, 719-736  (2006)

\bibitem{BN06} K. Beauchamp and R. Nicoara,
      Orthogonal maximal abelian *-subalgebras of the 6x6 matrices,
      preprint {\sl www.arxiv.org/abs/math.OA/0609076}

\bibitem{Sz06} F. Sz\"{o}ll\H{o}si,
      Parametrising complex Hadamard matrices,
      preprint math.CO/0610297

\bibitem{MS07} M. Matolcsi and  F. Sz\"{o}ll\H{o}si, 
 Towards a classification of $6\times 6$ complex Hadamard matrices, 
preprint math/0702043

\bibitem{Bxxx07} I. Bengtsson, W. Bruzda, A. Ericsson,
      J.-A. Larsson, W. Tadej, and K. {\.Z}yczkowski,
      MUBs and Hadamards of order six,
      J. Math. Phys.    {\bf 48}, 052106 (2007).  

\bibitem{Moor01} G. E. Moorhouse,
      The 2-Transitive Complex Hadamard Matrices,
      preprint 2001,
      {\sl http://www.uwyo.edu/moorhouse/pub/}

\bibitem{Tao04} T. Tao,
      Fuglede's conjecture is false in 5 and higher dimensions,
      Math. Res. Letters 11, 251-258 (2004).

\bibitem{Aig79} M. Aigner,
      Combinatorial Theory, Springer, New York, 1979.

\bibitem{Preprint07}
      W. Tadej and  K.  {\.Z}yczkowski,
      Defect of a unitary matrix,
      preprint 2007,
      {\sl http://arxiv.org/abs/math/0702510}

\end{thebibliography}
\end{document}